
\ifx\shlhetal\undefinedcontrolsequence\let\shlhetal\relax\fi
\def\fmtname{AmS-TeX}

\def\fmtversion{2.1}
\catcode`\@=11
\ifx\amstexloaded@\relax\catcode`\@=\active
  \endinput\else\let\amstexloaded@\relax\fi
\newlinechar=`\^^J
\def\W@{\immediate\write\sixt@@n}
\def\CR@{\W@{^^J\fmtname - Version \fmtversion^^J%
COPYRIGHT 1985, 1990, 1991 - AMERICAN MATHEMATICAL SOCIETY^^J%
Use of this macro package is not restricted provided^^J%
each use is acknowledged upon publication.^^J}}
\CR@ \everyjob{\CR@}
\message{Loading definitions for}
\message{misc utility macros,}
\toksdef\toks@@=2
\long\def\rightappend@#1\to#2{\toks@{\\{#1}}\toks@@
 =\expandafter{#2}\xdef#2{\the\toks@@\the\toks@}\toks@{}\toks@@{}}
\def\alloclist@{}
\newif\ifalloc@
\def\showallocations{{\def\\{\immediate\write\m@ne}\alloclist@}\alloc@true}
\def\alloc@#1#2#3#4#5{\global\advance\count1#1by\@ne
 \ch@ck#1#4#2\allocationnumber=\count1#1
 \global#3#5=\allocationnumber
 \edef\next@{\string#5=\string#2\the\allocationnumber}%
 \expandafter\rightappend@\next@\to\alloclist@}
\newcount\count@@
\newcount\count@@@
\def\FN@{\futurelet\next}
\def\DN@{\def\next@}
\def\DNii@{\def\nextii@}
\def\RIfM@{\relax\ifmmode}
\def\RIfMIfI@{\relax\ifmmode\ifinner}
\def\setboxz@h{\setbox\z@\hbox}
\def\wdz@{\wd\z@}
\def\boxz@{\box\z@}
\def\setbox@ne{\setbox\@ne}
\def\wd@ne{\wd\@ne}
\def\iterate{\body\expandafter\iterate\else\fi}
\def\err@#1{\errmessage{AmS-TeX error: #1}}
\newhelp\defaulthelp@{Sorry, I already gave what help I could...^^J
Maybe you should try asking a human?^^J
An error might have occurred before I noticed any problems.^^J
``If all else fails, read the instructions.''}
\def\Err@{\errhelp\defaulthelp@\err@}
\def\eat@#1{}
\def\in@#1#2{\def\in@@##1#1##2##3\in@@{\ifx\in@##2\in@false\else\in@true\fi}%
 \in@@#2#1\in@\in@@}
\newif\ifin@
\def\space@.{\futurelet\space@\relax}
\space@. %
\newhelp\athelp@
{Only certain combinations beginning with @ make sense to me.^^J
Perhaps you wanted \string\@\space for a printed @?^^J
I've ignored the character or group after @.}
{\catcode`\~=\active 
 \lccode`\~=`\@ \lowercase{\gdef~{\FN@\at@}}}
\def\at@{\let\next@\at@@
 \ifcat\noexpand\next a\else\ifcat\noexpand\next0\else
 \ifcat\noexpand\next\relax\else
   \let\next\at@@@\fi\fi\fi
 \next@}
\def\at@@#1{\expandafter
 \ifx\csname\space @\string#1\endcsname\relax
  \expandafter\at@@@ \else
  \csname\space @\string#1\expandafter\endcsname\fi}
\def\at@@@#1{\errhelp\athelp@ \err@{\Invalid@@ @}}
\def\atdef@#1{\expandafter\def\csname\space @\string#1\endcsname}
\newhelp\defahelp@{If you typed \string\define\space cs instead of
\string\define\string\cs\space^^J
I've substituted an inaccessible control sequence so that your^^J
definition will be completed without mixing me up too badly.^^J
If you typed \string\define{\string\cs} the inaccessible control sequence^^J
was defined to be \string\cs, and the rest of your^^J
definition appears as input.}
\newhelp\defbhelp@{I've ignored your definition, because it might^^J
conflict with other uses that are important to me.}
\def\define{\FN@\define@}
\def\define@{\ifcat\noexpand\next\relax
 \expandafter\define@@\else\errhelp\defahelp@                               
 \err@{\string\define\space must be followed by a control
 sequence}\expandafter\def\expandafter\nextii@\fi}                          
\def\undefined@@@@@@@@@@{}
\def\preloaded@@@@@@@@@@{}
\def\next@@@@@@@@@@{}
\def\define@@#1{\ifx#1\relax\errhelp\defbhelp@                              
 \err@{\string#1\space is already defined}\DN@{\DNii@}\else
 \expandafter\ifx\csname\expandafter\eat@\string                            
 #1@@@@@@@@@@\endcsname\undefined@@@@@@@@@@\errhelp\defbhelp@
 \err@{\string#1\space can't be defined}\DN@{\DNii@}\else
 \expandafter\ifx\csname\expandafter\eat@\string#1\endcsname\relax          
 \global\let#1\undefined\DN@{\def#1}\else\errhelp\defbhelp@
 \err@{\string#1\space is already defined}\DN@{\DNii@}\fi
 \fi\fi\next@}

\def\predefine#1#2{\let#1#2}
\def\undefine#1{\let#1\undefined}
\message{page layout,}
\newdimen\captionwidth@
\captionwidth@\hsize
\advance\captionwidth@-1.5in
\def\pagewidth#1{\hsize#1\relax
 \captionwidth@\hsize\advance\captionwidth@-1.5in}

\def\hcorrection#1{\advance\hoffset#1\relax}
\def\vcorrection#1{\advance\voffset#1\relax}
\message{accents/punctuation,}

\let\graveaccent\`
\let\acuteaccent\'
\let\tildeaccent\~
\let\hataccent\^
\let\underscore\_
\let\B\=
\let\D\.
\let\ic@\/
\def\/{\unskip\ic@}
\def\textfonti{\the\textfont\@ne}
\def\t#1#2{{\edef\next@{\the\font}\textfonti\accent"7F \next@#1#2}}
\def~{\unskip\nobreak\ \ignorespaces}
\def\.{.\spacefactor\@m}
\atdef@;{\leavevmode\null;}
\atdef@:{\leavevmode\null:}
\atdef@?{\leavevmode\null?}
\edef\@{\string @}
\def\relaxnext@{\let\next\relax}
\atdef@-{\relaxnext@\leavevmode
 \DN@{\ifx\next-\DN@-{\FN@\nextii@}\else
  \DN@{\leavevmode\hbox{-}}\fi\next@}%
 \DNii@{\ifx\next-\DN@-{\leavevmode\hbox{---}}\else
  \DN@{\leavevmode\hbox{--}}\fi\next@}%
 \FN@\next@}
\def\srdr@{\kern.16667em}
\def\drsr@{\kern.02778em}
\def\sldl@{\drsr@}
\def\dlsl@{\srdr@}
\atdef@"{\unskip\relaxnext@
 \DN@{\ifx\next\space@\DN@. {\FN@\nextii@}\else
  \DN@.{\FN@\nextii@}\fi\next@.}%
 \DNii@{\ifx\next`\DN@`{\FN@\nextiii@}\else
  \ifx\next\lq\DN@\lq{\FN@\nextiii@}\else
  \DN@####1{\FN@\nextiv@}\fi\fi\next@}%
 \def\nextiii@{\ifx\next`\DN@`{\sldl@``}\else\ifx\next\lq
  \DN@\lq{\sldl@``}\else\DN@{\dlsl@`}\fi\fi\next@}%
 \def\nextiv@{\ifx\next'\DN@'{\srdr@''}\else
  \ifx\next\rq\DN@\rq{\srdr@''}\else\DN@{\drsr@'}\fi\fi\next@}%
 \FN@\next@}

\def\textfontii{\the\textfont\tw@}
\def\lbrace@{\delimiter"4266308 }
\def\rbrace@{\delimiter"5267309 }
\def\{{\RIfM@\lbrace@\else{\textfontii f}\spacefactor\@m\fi}
\def\}{\RIfM@\rbrace@\else
 \let\@sf\empty\ifhmode\edef\@sf{\spacefactor\the\spacefactor}\fi
 {\textfontii g}\@sf\relax\fi}
\let\lbrace\{
\let\rbrace\}
\def\AmSTeX{{\textfontii A\kern-.1667em%
  \lower.5ex\hbox{M}\kern-.125emS}-\TeX}
\message{line and page breaks,}
\def\vmodeerr@#1{\Err@{\string#1\space not allowed between paragraphs}}
\def\mathmodeerr@#1{\Err@{\string#1\space not allowed in math mode}}
\def\linebreak{\RIfM@\mathmodeerr@\linebreak\else
 \ifhmode\unskip\unkern\break\else\vmodeerr@\linebreak\fi\fi}

\newskip\saveskip@
\def\allowlinebreak{\RIfM@\mathmodeerr@\allowlinebreak\else
 \ifhmode\saveskip@\lastskip\unskip
 \allowbreak\ifdim\saveskip@>\z@\hskip\saveskip@\fi
 \else\vmodeerr@\allowlinebreak\fi\fi}
\def\nolinebreak{\RIfM@\mathmodeerr@\nolinebreak\else
 \ifhmode\saveskip@\lastskip\unskip
 \nobreak\ifdim\saveskip@>\z@\hskip\saveskip@\fi
 \else\vmodeerr@\nolinebreak\fi\fi}
\def\newline{\relaxnext@
 \DN@{\RIfM@\expandafter\mathmodeerr@\expandafter\newline\else
  \ifhmode\ifx\next\par\else
  \expandafter\unskip\expandafter\null\expandafter\hfill\expandafter\break\fi
  \else
  \expandafter\vmodeerr@\expandafter\newline\fi\fi}%
 \FN@\next@}
\def\dmatherr@#1{\Err@{\string#1\space not allowed in display math mode}}
\def\nondmatherr@#1{\Err@{\string#1\space not allowed in non-display math
 mode}}
\def\onlydmatherr@#1{\Err@{\string#1\space allowed only in display math mode}}
\def\nonmatherr@#1{\Err@{\string#1\space allowed only in math mode}}
\def\mathbreak{\RIfMIfI@\break\else
 \dmatherr@\mathbreak\fi\else\nonmatherr@\mathbreak\fi}
\def\nomathbreak{\RIfMIfI@\nobreak\else
 \dmatherr@\nomathbreak\fi\else\nonmatherr@\nomathbreak\fi}
\def\allowmathbreak{\RIfMIfI@\allowbreak\else
 \dmatherr@\allowmathbreak\fi\else\nonmatherr@\allowmathbreak\fi}
\def\pagebreak{\RIfM@
 \ifinner\nondmatherr@\pagebreak\else\postdisplaypenalty-\@M\fi
 \else\ifvmode\removelastskip\break\else\vadjust{\break}\fi\fi}
\def\nopagebreak{\RIfM@
 \ifinner\nondmatherr@\nopagebreak\else\postdisplaypenalty\@M\fi
 \else\ifvmode\nobreak\else\vadjust{\nobreak}\fi\fi}
\def\nonvmodeerr@#1{\Err@{\string#1\space not allowed within a paragraph
 or in math}}
\def\vnonvmode@#1#2{\relaxnext@\DNii@{\ifx\next\par\DN@{#1}\else
 \DN@{#2}\fi\next@}%
 \ifvmode\DN@{#1}\else
 \DN@{\FN@\nextii@}\fi\next@}
\def\newpage{\vnonvmode@{\vfill\break}{\nonvmodeerr@\newpage}}
\def\smallpagebreak{\vnonvmode@\smallbreak{\nonvmodeerr@\smallpagebreak}}
\def\medpagebreak{\vnonvmode@\medbreak{\nonvmodeerr@\medpagebreak}}
\def\bigpagebreak{\vnonvmode@\bigbreak{\nonvmodeerr@\bigpagebreak}}
\def\NoBlackBoxes{\global\overfullrule\z@}
\def\BlackBoxes{\global\overfullrule5\p@}
\def\Invalid@#1{\def#1{\Err@{\Invalid@@\string#1}}}
\def\Invalid@@{Invalid use of }
\message{figures,}
\Invalid@\caption
\Invalid@\captionwidth
\newdimen\smallcaptionwidth@
\def\topspace{\mid@false\ins@}
\def\midspace{\mid@true\ins@}
\newif\ifmid@
\def\captionfont@{}
\def\ins@#1{\relaxnext@\allowbreak
 \smallcaptionwidth@\captionwidth@\gdef\thespace@{#1}%
 \DN@{\ifx\next\space@\DN@. {\FN@\nextii@}\else
  \DN@.{\FN@\nextii@}\fi\next@.}%
 \DNii@{\ifx\next\caption\DN@\caption{\FN@\nextiii@}%
  \else\let\next@\nextiv@\fi\next@}%
 \def\nextiv@{\vnonvmode@
  {\ifmid@\expandafter\midinsert\else\expandafter\topinsert\fi
   \vbox to\thespace@{}\endinsert}
  {\ifmid@\nonvmodeerr@\midspace\else\nonvmodeerr@\topspace\fi}}%
 \def\nextiii@{\ifx\next\captionwidth\expandafter\nextv@
  \else\expandafter\nextvi@\fi}%
 \def\nextv@\captionwidth##1##2{\smallcaptionwidth@##1\relax\nextvi@{##2}}%
 \def\nextvi@##1{\def\thecaption@{\captionfont@##1}%
  \DN@{\ifx\next\space@\DN@. {\FN@\nextvii@}\else
   \DN@.{\FN@\nextvii@}\fi\next@.}%
  \FN@\next@}%
 \def\nextvii@{\vnonvmode@
  {\ifmid@\expandafter\midinsert\else
  \expandafter\topinsert\fi\vbox to\thespace@{}\nobreak\smallskip
  \setboxz@h{\noindent\ignorespaces\thecaption@\unskip}%
  \ifdim\wdz@>\smallcaptionwidth@\centerline{\vbox{\hsize\smallcaptionwidth@
   \noindent\ignorespaces\thecaption@\unskip}}%
  \else\centerline{\boxz@}\fi\endinsert}
  {\ifmid@\nonvmodeerr@\midspace
  \else\nonvmodeerr@\topspace\fi}}%
 \FN@\next@}
\message{comments,}
\def\newcodes@{\catcode`\\12\catcode`\{12\catcode`\}12\catcode`\#12%
 \catcode`\%12\relax}
\def\oldcodes@{\catcode`\\0\catcode`\{1\catcode`\}2\catcode`\#6%
 \catcode`\%14\relax}
\def\comment{\newcodes@\endlinechar=10 \comment@}
{\lccode`\0=`\\
\lowercase{\gdef\comment@#1^^J{\comment@@#10endcomment\comment@@@}%
\gdef\comment@@#10endcomment{\FN@\comment@@@}%
\gdef\comment@@@#1\comment@@@{\ifx\next\comment@@@\let\next\comment@
 \else\def\next{\oldcodes@\endlinechar=`\^^M\relax}%
 \fi\next}}}
\def\pr@m@s{\ifx'\next\DN@##1{\prim@s}\else\let\next@\egroup\fi\next@}
\def\prime{{\null\prime@\null}}
\mathchardef\prime@="0230
\let\dsize\displaystyle

\let\ssize\scriptstyle

\message{math spacing,}
\def\,{\RIfM@\mskip\thinmuskip\relax\else\kern.16667em\fi}
\def\!{\RIfM@\mskip-\thinmuskip\relax\else\kern-.16667em\fi}
\let\thinspace\,
\let\negthinspace\!
\def\medspace{\RIfM@\mskip\medmuskip\relax\else\kern.222222em\fi}
\def\negmedspace{\RIfM@\mskip-\medmuskip\relax\else\kern-.222222em\fi}
\def\thickspace{\RIfM@\mskip\thickmuskip\relax\else\kern.27777em\fi}
\let\;\thickspace
\def\negthickspace{\RIfM@\mskip-\thickmuskip\relax\else
 \kern-.27777em\fi}
\atdef@,{\RIfM@\mskip.1\thinmuskip\else\leavevmode\null,\fi}
\atdef@!{\RIfM@\mskip-.1\thinmuskip\else\leavevmode\null!\fi}
\atdef@.{\RIfM@&&\else\leavevmode.\spacefactor3000 \fi}
\def\and{\DOTSB\;\mathchar"3026 \;}

\message{fractions,}
\def\frac#1#2{{#1\over#2}}

\newdimen\ex@
\ex@.2326ex
\Invalid@\thickness
\def\thickfrac{\relaxnext@
 \DN@{\ifx\next\thickness\let\next@\nextii@\else
 \DN@{\nextii@\thickness1}\fi\next@}%
 \DNii@\thickness##1##2##3{{##2\above##1\ex@##3}}%
 \FN@\next@}

\def\thickfracwithdelims#1#2{\relaxnext@\def\ldelim@{#1}\def\rdelim@{#2}%
 \DN@{\ifx\next\thickness\let\next@\nextii@\else
 \DN@{\nextii@\thickness1}\fi\next@}%
 \DNii@\thickness##1##2##3{{##2\abovewithdelims
 \ldelim@\rdelim@##1\ex@##3}}%
 \FN@\next@}

\def\:{\nobreak\hskip.1111em\mathpunct{}\nonscript\mkern-\thinmuskip{:}\hskip
 .3333emplus.0555em\relax}
\def\snug{\unskip\kern-\mathsurround}
\message{smash commands,}
\def\topsmash{\top@true\bot@false\smash@}
\def\botsmash{\top@false\bot@true\smash@}
\newif\iftop@
\newif\ifbot@
\def\smash{\top@true\bot@true\smash@}
\def\smash@{\RIfM@\expandafter\mathpalette\expandafter\mathsm@sh\else
 \expandafter\makesm@sh\fi}
\def\finsm@sh{\iftop@\ht\z@\z@\fi\ifbot@\dp\z@\z@\fi\leavevmode\boxz@}
\message{large operator symbols,}
\def\LimitsOnSums{\global\let\slimits@\displaylimits}
\def\NoLimitsOnSums{\global\let\slimits@\nolimits}
\LimitsOnSums
\mathchardef\coprod@="1360       \def\coprod{\DOTSB\coprod@\slimits@}
\mathchardef\bigvee@="1357       \def\bigvee{\DOTSB\bigvee@\slimits@}
\mathchardef\bigwedge@="1356     \def\bigwedge{\DOTSB\bigwedge@\slimits@}
\mathchardef\biguplus@="1355     \def\biguplus{\DOTSB\biguplus@\slimits@}
\mathchardef\bigcap@="1354       \def\bigcap{\DOTSB\bigcap@\slimits@}
\mathchardef\bigcup@="1353       \def\bigcup{\DOTSB\bigcup@\slimits@}
\mathchardef\prod@="1351         \def\prod{\DOTSB\prod@\slimits@}
\mathchardef\sum@="1350          \def\sum{\DOTSB\sum@\slimits@}
\mathchardef\bigotimes@="134E    \def\bigotimes{\DOTSB\bigotimes@\slimits@}
\mathchardef\bigoplus@="134C     \def\bigoplus{\DOTSB\bigoplus@\slimits@}
\mathchardef\bigodot@="134A      \def\bigodot{\DOTSB\bigodot@\slimits@}
\mathchardef\bigsqcup@="1346     \def\bigsqcup{\DOTSB\bigsqcup@\slimits@}
\message{integrals,}
\def\LimitsOnInts{\global\let\ilimits@\displaylimits}
\def\NoLimitsOnInts{\global\let\ilimits@\nolimits}
\NoLimitsOnInts
\def\int{\DOTSI\intop\ilimits@}
\def\oint{\DOTSI\ointop\ilimits@}
\def\intic@{\mathchoice{\hskip.5em}{\hskip.4em}{\hskip.4em}{\hskip.4em}}
\def\negintic@{\mathchoice
 {\hskip-.5em}{\hskip-.4em}{\hskip-.4em}{\hskip-.4em}}
\def\intkern@{\mathchoice{\!\!\!}{\!\!}{\!\!}{\!\!}}
\def\intdots@{\mathchoice{\plaincdots@}
 {{\cdotp}\mkern1.5mu{\cdotp}\mkern1.5mu{\cdotp}}
 {{\cdotp}\mkern1mu{\cdotp}\mkern1mu{\cdotp}}
 {{\cdotp}\mkern1mu{\cdotp}\mkern1mu{\cdotp}}}
\newcount\intno@
\def\iint{\DOTSI\intno@\tw@\FN@\ints@}
\def\iiint{\DOTSI\intno@\thr@@\FN@\ints@}
\def\iiiint{\DOTSI\intno@4 \FN@\ints@}
\def\idotsint{\DOTSI\intno@\z@\FN@\ints@}
\def\ints@{\findlimits@\ints@@}
\newif\iflimtoken@
\newif\iflimits@
\def\findlimits@{\limtoken@true\ifx\next\limits\limits@true
 \else\ifx\next\nolimits\limits@false\else
 \limtoken@false\ifx\ilimits@\nolimits\limits@false\else
 \ifinner\limits@false\else\limits@true\fi\fi\fi\fi}
\def\multint@{\int\ifnum\intno@=\z@\intdots@                                
 \else\intkern@\fi                                                          
 \ifnum\intno@>\tw@\int\intkern@\fi                                         
 \ifnum\intno@>\thr@@\int\intkern@\fi                                       
 \int}                                                                      
\def\multintlimits@{\intop\ifnum\intno@=\z@\intdots@\else\intkern@\fi
 \ifnum\intno@>\tw@\intop\intkern@\fi
 \ifnum\intno@>\thr@@\intop\intkern@\fi\intop}
\def\ints@@{\iflimtoken@                                                    
 \def\ints@@@{\iflimits@\negintic@\mathop{\intic@\multintlimits@}\limits    
  \else\multint@\nolimits\fi                                                
  \eat@}                                                                    
 \else                                                                      
 \def\ints@@@{\iflimits@\negintic@
  \mathop{\intic@\multintlimits@}\limits\else
  \multint@\nolimits\fi}\fi\ints@@@}
\def\LimitsOnNames{\global\let\nlimits@\displaylimits}
\def\NoLimitsOnNames{\global\let\nlimits@\nolimits@}
\LimitsOnNames
\def\nolimits@{\relaxnext@
 \DN@{\ifx\next\limits\DN@\limits{\nolimits}\else
  \let\next@\nolimits\fi\next@}%
 \FN@\next@}
\message{operator names,}
\def\newmcodes@{\mathcode`\'"27\mathcode`\*"2A\mathcode`\."613A%
 \mathcode`\-"2D\mathcode`\/"2F\mathcode`\:"603A }
\def\operatorname#1{\mathop{\newmcodes@\kern\z@\fam\z@#1}\nolimits@}
\def\operatornamewithlimits#1{\mathop{\newmcodes@\kern\z@\fam\z@#1}\nlimits@}
\def\qopname@#1{\mathop{\fam\z@#1}\nolimits@}
\def\qopnamewl@#1{\mathop{\fam\z@#1}\nlimits@}
\def\arccos{\qopname@{arccos}}
\def\arcsin{\qopname@{arcsin}}
\def\arctan{\qopname@{arctan}}
\def\arg{\qopname@{arg}}
\def\cos{\qopname@{cos}}
\def\cosh{\qopname@{cosh}}
\def\cot{\qopname@{cot}}
\def\coth{\qopname@{coth}}
\def\csc{\qopname@{csc}}
\def\deg{\qopname@{deg}}
\def\det{\qopnamewl@{det}}
\def\dim{\qopname@{dim}}
\def\exp{\qopname@{exp}}
\def\gcd{\qopnamewl@{gcd}}
\def\hom{\qopname@{hom}}
\def\inf{\qopnamewl@{inf}}
\def\injlim{\qopnamewl@{inj\,lim}}
\def\ker{\qopname@{ker}}
\def\lg{\qopname@{lg}}
\def\lim{\qopnamewl@{lim}}
\def\liminf{\qopnamewl@{lim\,inf}}
\def\limsup{\qopnamewl@{lim\,sup}}
\def\ln{\qopname@{ln}}
\def\log{\qopname@{log}}
\def\max{\qopnamewl@{max}}
\def\min{\qopnamewl@{min}}
\def\Pr{\qopnamewl@{Pr}}
\def\projlim{\qopnamewl@{proj\,lim}}
\def\sec{\qopname@{sec}}
\def\sin{\qopname@{sin}}
\def\sinh{\qopname@{sinh}}
\def\sup{\qopnamewl@{sup}}
\def\tan{\qopname@{tan}}
\def\tanh{\qopname@{tanh}}
\def\varinjlim{\mathop{\vtop{\ialign{##\crcr
 \hfil\rm lim\hfil\crcr\noalign{\nointerlineskip}\rightarrowfill\crcr
 \noalign{\nointerlineskip\kern-\ex@}\crcr}}}}
\def\varprojlim{\mathop{\vtop{\ialign{##\crcr
 \hfil\rm lim\hfil\crcr\noalign{\nointerlineskip}\leftarrowfill\crcr
 \noalign{\nointerlineskip\kern-\ex@}\crcr}}}}
\def\varliminf{\mathop{\underline{\vrule height\z@ depth.2exwidth\z@
 \hbox{\rm lim}}}}

\newdimen\buffer@
\buffer@\fontdimen13 \tenex
\newdimen\buffer
\buffer\buffer@

\def\ResetBuffer{\fontdimen13 \tenex\buffer@\global\buffer\buffer@}
\def\shave#1{\mathop{\hbox{$\m@th\fontdimen13 \tenex\z@                     
 \displaystyle{#1}$}}\fontdimen13 \tenex\buffer}

\message{multilevel sub/superscripts,}
\Invalid@\\
\def\Let@{\relax\iffalse{\fi\let\\=\cr\iffalse}\fi}
\Invalid@\vspace
\def\vspace@{\def\vspace##1{\crcr\noalign{\vskip##1\relax}}}
\def\multilimits@{\bgroup\vspace@\Let@
 \baselineskip\fontdimen10 \scriptfont\tw@
 \advance\baselineskip\fontdimen12 \scriptfont\tw@
 \lineskip\thr@@\fontdimen8 \scriptfont\thr@@
 \lineskiplimit\lineskip
 \vbox\bgroup\ialign\bgroup\hfil$\m@th\scriptstyle{##}$\hfil\crcr}
\def\Sb{_\multilimits@}
\def\endSb{\crcr\egroup\egroup\egroup}
\def\Sp{^\multilimits@}

\def\spreadlines#1{\RIfMIfI@\onlydmatherr@\spreadlines\else
 \openup#1\relax\fi\else\onlydmatherr@\spreadlines\fi}
\def\Mathstrut@{\copy\Mathstrutbox@}
\newbox\Mathstrutbox@
\setbox\Mathstrutbox@\null
\setboxz@h{$\m@th($}
\ht\Mathstrutbox@\ht\z@
\dp\Mathstrutbox@\dp\z@
\message{matrices,}
\newdimen\spreadmlines@
\def\spreadmatrixlines#1{\RIfMIfI@
 \onlydmatherr@\spreadmatrixlines\else
 \spreadmlines@#1\relax\fi\else\onlydmatherr@\spreadmatrixlines\fi}
\def\matrix{\null\,\vcenter\bgroup\Let@\vspace@
 \normalbaselines\openup\spreadmlines@\ialign
 \bgroup\hfil$\m@th##$\hfil&&\quad\hfil$\m@th##$\hfil\crcr
 \Mathstrut@\crcr\noalign{\kern-\baselineskip}}
\def\endmatrix{\crcr\Mathstrut@\crcr\noalign{\kern-\baselineskip}\egroup
 \egroup\,}
\def\format{\crcr\egroup\iffalse{\fi\ifnum`}=0 \fi\format@}
\newtoks\hashtoks@
\hashtoks@{#}
\def\format@#1\\{\def\preamble@{#1}%
 \def\l{$\m@th\the\hashtoks@$\hfil}%
 \def\c{\hfil$\m@th\the\hashtoks@$\hfil}%
 \def\r{\hfil$\m@th\the\hashtoks@$}%
 \edef\preamble@@{\preamble@}\ifnum`{=0 \fi\iffalse}\fi
 \ialign\bgroup\span\preamble@@\crcr}
\def\smallmatrix{\null\,\vcenter\bgroup\vspace@\Let@
 \baselineskip9\ex@\lineskip\ex@
 \ialign\bgroup\hfil$\m@th\scriptstyle{##}$\hfil&&\thickspace\hfil
 $\m@th\scriptstyle{##}$\hfil\crcr}
\def\endsmallmatrix{\crcr\egroup\egroup\,}

\newmuskip\dotsspace@
\dotsspace@1.5mu
\def\strip@#1 {#1}
\def\spacehdots#1\for#2{\multispan{#2}\xleaders
 \hbox{$\m@th\mkern\strip@#1 \dotsspace@.\mkern\strip@#1 \dotsspace@$}\hfill}
\def\hdotsfor#1{\spacehdots\@ne\for{#1}}
\def\multispan@#1{\omit\mscount#1\unskip\loop\ifnum\mscount>\@ne\sp@n\repeat}
\def\spaceinnerhdots#1\for#2\after#3{\multispan@{\strip@#2 }#3\xleaders
 \hbox{$\m@th\mkern\strip@#1 \dotsspace@.\mkern\strip@#1 \dotsspace@$}\hfill}
\def\innerhdotsfor#1\after#2{\spaceinnerhdots\@ne\for#1\after{#2}}
\def\cases{\bgroup\spreadmlines@\jot\left\{\,\matrix\format\l&\quad\l\\}
\def\endcases{\endmatrix\right.\egroup}
\message{multiline displays,}
\newif\ifinany@
\newif\ifinalign@
\newif\ifingather@
\def\strut@{\copy\strutbox@}
\newbox\strutbox@
\setbox\strutbox@\hbox{\vrule height8\p@ depth3\p@ width\z@}
\def\topaligned{\null\,\vtop\aligned@}
\def\botaligned{\null\,\vbox\aligned@}
\def\aligned{\null\,\vcenter\aligned@}
\def\aligned@{\bgroup\vspace@\Let@
 \ifinany@\else\openup\jot\fi\ialign
 \bgroup\hfil\strut@$\m@th\displaystyle{##}$&
 $\m@th\displaystyle{{}##}$\hfil\crcr}
\def\endaligned{\crcr\egroup\egroup}

\def\alignedat#1{\null\,\vcenter\bgroup\doat@{#1}\vspace@\Let@
 \ifinany@\else\openup\jot\fi\ialign\bgroup\span\preamble@@\crcr}
\newcount\atcount@
\def\doat@#1{\toks@{\hfil\strut@$\m@th
 \displaystyle{\the\hashtoks@}$&$\m@th\displaystyle
 {{}\the\hashtoks@}$\hfil}
 \atcount@#1\relax\advance\atcount@\m@ne                                    
 \loop\ifnum\atcount@>\z@\toks@=\expandafter{\the\toks@&\hfil$\m@th
 \displaystyle{\the\hashtoks@}$&$\m@th
 \displaystyle{{}\the\hashtoks@}$\hfil}\advance
  \atcount@\m@ne\repeat                                                     
 \xdef\preamble@{\the\toks@}\xdef\preamble@@{\preamble@}}

\def\gathered{\null\,\vcenter\bgroup\vspace@\Let@
 \ifinany@\else\openup\jot\fi\ialign
 \bgroup\hfil\strut@$\m@th\displaystyle{##}$\hfil\crcr}
\def\endgathered{\crcr\egroup\egroup}
\newif\iftagsleft@
\def\TagsOnLeft{\global\tagsleft@true}
\def\TagsOnRight{\global\tagsleft@false}
\TagsOnLeft
\newif\ifmathtags@
\def\TagsAsMath{\global\mathtags@true}
\def\TagsAsText{\global\mathtags@false}
\TagsAsText
\def\tagform@#1{\hbox{\rm(\ignorespaces#1\unskip)}}
\def\thetag{\leavevmode\tagform@}
\def\tag#1$${\iftagsleft@\leqno\else\eqno\fi                                
 \maketag@#1\maketag@                                                       
 $$}                                                                        
\def\maketag@{\FN@\maketag@@}
\def\maketag@@{\ifx\next"\expandafter\maketag@@@\else\expandafter\maketag@@@@
 \fi}
\def\maketag@@@"#1"#2\maketag@{\hbox{\rm#1}}                                
\def\maketag@@@@#1\maketag@{\ifmathtags@\tagform@{$\m@th#1$}\else
 \tagform@{#1}\fi}
\interdisplaylinepenalty\@M
\def\allowdisplaybreaks{\RIfMIfI@
 \onlydmatherr@\allowdisplaybreaks\else
 \interdisplaylinepenalty\z@\fi\else\onlydmatherr@\allowdisplaybreaks\fi}
\Invalid@\allowdisplaybreak
\Invalid@\displaybreak
\Invalid@\intertext
\def\allowdisplaybreak@{\def\allowdisplaybreak{\crcr\noalign{\allowbreak}}}
\def\displaybreak@{\def\displaybreak{\crcr\noalign{\break}}}
\def\intertext@{\def\intertext##1{\crcr\noalign{%
 \penalty\postdisplaypenalty \vskip\belowdisplayskip
 \vbox{\normalbaselines\noindent##1}%
 \penalty\predisplaypenalty \vskip\abovedisplayskip}}}
\newskip\centering@
\centering@\z@ plus\@m\p@
\def\align{\relax\ifingather@\DN@{\csname align (in
  \string\gather)\endcsname}\else
 \ifmmode\ifinner\DN@{\onlydmatherr@\align}\else
  \let\next@\align@\fi
 \else\DN@{\onlydmatherr@\align}\fi\fi\next@}
\newhelp\andhelp@
{An extra & here is so disastrous that you should probably exit^^J
and fix things up.}
\newif\iftag@
\newcount\and@
\def\align@{\inalign@true\inany@true
 \vspace@\allowdisplaybreak@\displaybreak@\intertext@
 \def\tag{\global\tag@true\ifnum\and@=\z@\DN@{&&}\else
          \DN@{&}\fi\next@}%
 \iftagsleft@\DN@{\csname align \endcsname}\else
  \DN@{\csname align \space\endcsname}\fi\next@}
\def\Tag@{\iftag@\else\errhelp\andhelp@\err@{Extra & on this line}\fi}
\newdimen\lwidth@
\newdimen\rwidth@
\newdimen\maxlwidth@
\newdimen\maxrwidth@
\newdimen\totwidth@
\def\measure@#1\endalign{\lwidth@\z@\rwidth@\z@\maxlwidth@\z@\maxrwidth@\z@
 \global\and@\z@                                                            
 \setbox@ne\vbox                                                            
  {\everycr{\noalign{\global\tag@false\global\and@\z@}}\Let@                
  \halign{\setboxz@h{$\m@th\displaystyle{\@lign##}$}
   \global\lwidth@\wdz@                                                     
   \ifdim\lwidth@>\maxlwidth@\global\maxlwidth@\lwidth@\fi                  
   \global\advance\and@\@ne                                                 
   &\setboxz@h{$\m@th\displaystyle{{}\@lign##}$}\global\rwidth@\wdz@        
   \ifdim\rwidth@>\maxrwidth@\global\maxrwidth@\rwidth@\fi                  
   \global\advance\and@\@ne                                                
   &\Tag@
   \eat@{##}\crcr#1\crcr}}
 \totwidth@\maxlwidth@\advance\totwidth@\maxrwidth@}                       
\def\displ@y@{\global\dt@ptrue\openup\jot
 \everycr{\noalign{\global\tag@false\global\and@\z@\ifdt@p\global\dt@pfalse
 \vskip-\lineskiplimit\vskip\normallineskiplimit\else
 \penalty\interdisplaylinepenalty\fi}}}
\def\black@#1{\noalign{\ifdim#1>\displaywidth
 \dimen@\prevdepth\nointerlineskip                                          
 \vskip-\ht\strutbox@\vskip-\dp\strutbox@                                   
 \vbox{\noindent\hbox to#1{\strut@\hfill}}
 \prevdepth\dimen@                                                          
 \fi}}
\expandafter\def\csname align \space\endcsname#1\endalign
 {\measure@#1\endalign\global\and@\z@                                       
 \ifingather@\everycr{\noalign{\global\and@\z@}}\else\displ@y@\fi           
 \Let@\tabskip\centering@                                                   
 \halign to\displaywidth
  {\hfil\strut@\setboxz@h{$\m@th\displaystyle{\@lign##}$}
  \global\lwidth@\wdz@\boxz@\global\advance\and@\@ne                        
  \tabskip\z@skip                                                           
  &\setboxz@h{$\m@th\displaystyle{{}\@lign##}$}
  \global\rwidth@\wdz@\boxz@\hfill\global\advance\and@\@ne                  
  \tabskip\centering@                                                       
  &\setboxz@h{\@lign\strut@\maketag@##\maketag@}
  \dimen@\displaywidth\advance\dimen@-\totwidth@
  \divide\dimen@\tw@\advance\dimen@\maxrwidth@\advance\dimen@-\rwidth@     
  \ifdim\dimen@<\tw@\wdz@\llap{\vtop{\normalbaselines\null\boxz@}}
  \else\llap{\boxz@}\fi                                                    
  \tabskip\z@skip                                                          
  \crcr#1\crcr                                                             
  \black@\totwidth@}}                                                      
\newdimen\lineht@
\expandafter\def\csname align \endcsname#1\endalign{\measure@#1\endalign
 \global\and@\z@
 \ifdim\totwidth@>\displaywidth\let\displaywidth@\totwidth@\else
  \let\displaywidth@\displaywidth\fi                                        
 \ifingather@\everycr{\noalign{\global\and@\z@}}\else\displ@y@\fi
 \Let@\tabskip\centering@\halign to\displaywidth
  {\hfil\strut@\setboxz@h{$\m@th\displaystyle{\@lign##}$}%
  \global\lwidth@\wdz@\global\lineht@\ht\z@                                 
  \boxz@\global\advance\and@\@ne
  \tabskip\z@skip&\setboxz@h{$\m@th\displaystyle{{}\@lign##}$}%
  \global\rwidth@\wdz@\ifdim\ht\z@>\lineht@\global\lineht@\ht\z@\fi         
  \boxz@\hfil\global\advance\and@\@ne
  \tabskip\centering@&\kern-\displaywidth@                                  
  \setboxz@h{\@lign\strut@\maketag@##\maketag@}%
  \dimen@\displaywidth\advance\dimen@-\totwidth@
  \divide\dimen@\tw@\advance\dimen@\maxlwidth@\advance\dimen@-\lwidth@
  \ifdim\dimen@<\tw@\wdz@
   \rlap{\vbox{\normalbaselines\boxz@\vbox to\lineht@{}}}\else
   \rlap{\boxz@}\fi
  \tabskip\displaywidth@\crcr#1\crcr\black@\totwidth@}}
\expandafter\def\csname align (in \string\gather)\endcsname
  #1\endalign{\vcenter{\align@#1\endalign}}
\Invalid@\endalign
\newif\ifxat@
\def\alignat{\RIfMIfI@\DN@{\onlydmatherr@\alignat}\else
 \DN@{\csname alignat \endcsname}\fi\else
 \DN@{\onlydmatherr@\alignat}\fi\next@}
\newif\ifmeasuring@
\newbox\savealignat@
\expandafter\def\csname alignat \endcsname#1#2\endalignat                   
 {\inany@true\xat@false
 \def\tag{\global\tag@true\count@#1\relax\multiply\count@\tw@
  \xdef\tag@{}\loop\ifnum\count@>\and@\xdef\tag@{&\tag@}\advance\count@\m@ne
  \repeat\tag@}%
 \vspace@\allowdisplaybreak@\displaybreak@\intertext@
 \displ@y@\measuring@true                                                   
 \setbox\savealignat@\hbox{$\m@th\displaystyle\Let@
  \attag@{#1}
  \vbox{\halign{\span\preamble@@\crcr#2\crcr}}$}%
 \measuring@false                                                           
 \Let@\attag@{#1}
 \tabskip\centering@\halign to\displaywidth
  {\span\preamble@@\crcr#2\crcr                                             
  \black@{\wd\savealignat@}}}                                               
\Invalid@\endalignat
\def\xalignat{\RIfMIfI@
 \DN@{\onlydmatherr@\xalignat}\else
 \DN@{\csname xalignat \endcsname}\fi\else
 \DN@{\onlydmatherr@\xalignat}\fi\next@}
\expandafter\def\csname xalignat \endcsname#1#2\endxalignat
 {\inany@true\xat@true
 \def\tag{\global\tag@true\def\tag@{}\count@#1\relax\multiply\count@\tw@
  \loop\ifnum\count@>\and@\xdef\tag@{&\tag@}\advance\count@\m@ne\repeat\tag@}%
 \vspace@\allowdisplaybreak@\displaybreak@\intertext@
 \displ@y@\measuring@true\setbox\savealignat@\hbox{$\m@th\displaystyle\Let@
 \attag@{#1}\vbox{\halign{\span\preamble@@\crcr#2\crcr}}$}%
 \measuring@false\Let@
 \attag@{#1}\tabskip\centering@\halign to\displaywidth
 {\span\preamble@@\crcr#2\crcr\black@{\wd\savealignat@}}}
\def\attag@#1{\let\Maketag@\maketag@\let\TAG@\Tag@                          
 \let\Tag@=0\let\maketag@=0
 \ifmeasuring@\def\llap@##1{\setboxz@h{##1}\hbox to\tw@\wdz@{}}%
  \def\rlap@##1{\setboxz@h{##1}\hbox to\tw@\wdz@{}}\else
  \let\llap@\llap\let\rlap@\rlap\fi                                         
 \toks@{\hfil\strut@$\m@th\displaystyle{\@lign\the\hashtoks@}$\tabskip\z@skip
  \global\advance\and@\@ne&$\m@th\displaystyle{{}\@lign\the\hashtoks@}$\hfil
  \ifxat@\tabskip\centering@\fi\global\advance\and@\@ne}
 \iftagsleft@
  \toks@@{\tabskip\centering@&\Tag@\kern-\displaywidth
   \rlap@{\@lign\maketag@\the\hashtoks@\maketag@}%
   \global\advance\and@\@ne\tabskip\displaywidth}\else
  \toks@@{\tabskip\centering@&\Tag@\llap@{\@lign\maketag@
   \the\hashtoks@\maketag@}\global\advance\and@\@ne\tabskip\z@skip}\fi      
 \atcount@#1\relax\advance\atcount@\m@ne
 \loop\ifnum\atcount@>\z@
 \toks@=\expandafter{\the\toks@&\hfil$\m@th\displaystyle{\@lign
  \the\hashtoks@}$\global\advance\and@\@ne
  \tabskip\z@skip&$\m@th\displaystyle{{}\@lign\the\hashtoks@}$\hfil\ifxat@
  \tabskip\centering@\fi\global\advance\and@\@ne}\advance\atcount@\m@ne
 \repeat                                                                    
 \xdef\preamble@{\the\toks@\the\toks@@}
 \xdef\preamble@@{\preamble@}
 \let\maketag@\Maketag@\let\Tag@\TAG@}                                      
\Invalid@\endxalignat
\def\xxalignat{\RIfMIfI@
 \DN@{\onlydmatherr@\xxalignat}\else\DN@{\csname xxalignat
  \endcsname}\fi\else
 \DN@{\onlydmatherr@\xxalignat}\fi\next@}
\expandafter\def\csname xxalignat \endcsname#1#2\endxxalignat{\inany@true
 \vspace@\allowdisplaybreak@\displaybreak@\intertext@
 \displ@y\setbox\savealignat@\hbox{$\m@th\displaystyle\Let@
 \xxattag@{#1}\vbox{\halign{\span\preamble@@\crcr#2\crcr}}$}%
 \Let@\xxattag@{#1}\tabskip\z@skip\halign to\displaywidth
 {\span\preamble@@\crcr#2\crcr\black@{\wd\savealignat@}}}
\def\xxattag@#1{\toks@{\tabskip\z@skip\hfil\strut@
 $\m@th\displaystyle{\the\hashtoks@}$&%
 $\m@th\displaystyle{{}\the\hashtoks@}$\hfil\tabskip\centering@&}%
 \atcount@#1\relax\advance\atcount@\m@ne\loop\ifnum\atcount@>\z@
 \toks@=\expandafter{\the\toks@&\hfil$\m@th\displaystyle{\the\hashtoks@}$%
  \tabskip\z@skip&$\m@th\displaystyle{{}\the\hashtoks@}$\hfil
  \tabskip\centering@}\advance\atcount@\m@ne\repeat
 \xdef\preamble@{\the\toks@\tabskip\z@skip}\xdef\preamble@@{\preamble@}}
\Invalid@\endxxalignat
\newdimen\gwidth@
\newdimen\gmaxwidth@
\def\gmeasure@#1\endgather{\gwidth@\z@\gmaxwidth@\z@\setbox@ne\vbox{\Let@
 \halign{\setboxz@h{$\m@th\displaystyle{##}$}\global\gwidth@\wdz@
 \ifdim\gwidth@>\gmaxwidth@\global\gmaxwidth@\gwidth@\fi
 &\eat@{##}\crcr#1\crcr}}}
\def\gather{\RIfMIfI@\DN@{\onlydmatherr@\gather}\else
 \ingather@true\inany@true\def\tag{&}%
 \vspace@\allowdisplaybreak@\displaybreak@\intertext@
 \displ@y\Let@
 \iftagsleft@\DN@{\csname gather \endcsname}\else
  \DN@{\csname gather \space\endcsname}\fi\fi
 \else\DN@{\onlydmatherr@\gather}\fi\next@}
\expandafter\def\csname gather \space\endcsname#1\endgather
 {\gmeasure@#1\endgather\tabskip\centering@
 \halign to\displaywidth{\hfil\strut@\setboxz@h{$\m@th\displaystyle{##}$}%
 \global\gwidth@\wdz@\boxz@\hfil&
 \setboxz@h{\strut@{\maketag@##\maketag@}}%
 \dimen@\displaywidth\advance\dimen@-\gwidth@
 \ifdim\dimen@>\tw@\wdz@\llap{\boxz@}\else
 \llap{\vtop{\normalbaselines\null\boxz@}}\fi
 \tabskip\z@skip\crcr#1\crcr\black@\gmaxwidth@}}
\newdimen\glineht@
\expandafter\def\csname gather \endcsname#1\endgather{\gmeasure@#1\endgather
 \ifdim\gmaxwidth@>\displaywidth\let\gdisplaywidth@\gmaxwidth@\else
 \let\gdisplaywidth@\displaywidth\fi\tabskip\centering@\halign to\displaywidth
 {\hfil\strut@\setboxz@h{$\m@th\displaystyle{##}$}%
 \global\gwidth@\wdz@\global\glineht@\ht\z@\boxz@\hfil&\kern-\gdisplaywidth@
 \setboxz@h{\strut@{\maketag@##\maketag@}}%
 \dimen@\displaywidth\advance\dimen@-\gwidth@
 \ifdim\dimen@>\tw@\wdz@\rlap{\boxz@}\else
 \rlap{\vbox{\normalbaselines\boxz@\vbox to\glineht@{}}}\fi
 \tabskip\gdisplaywidth@\crcr#1\crcr\black@\gmaxwidth@}}
\newif\ifctagsplit@
\def\CenteredTagsOnSplits{\global\ctagsplit@true}
\def\TopOrBottomTagsOnSplits{\global\ctagsplit@false}
\TopOrBottomTagsOnSplits
\def\split{\relax\ifinany@\let\next@\insplit@\else
 \ifmmode\ifinner\def\next@{\onlydmatherr@\split}\else
 \let\next@\outsplit@\fi\else
 \def\next@{\onlydmatherr@\split}\fi\fi\next@}
\def\insplit@{\global\setbox\z@\vbox\bgroup\vspace@\Let@\ialign\bgroup
 \hfil\strut@$\m@th\displaystyle{##}$&$\m@th\displaystyle{{}##}$\hfill\crcr}
\def\endsplit{\crcr\egroup\egroup\iftagsleft@\expandafter\lendsplit@\else
 \expandafter\rendsplit@\fi}
\def\rendsplit@{\global\setbox9 \vbox
 {\unvcopy\z@\global\setbox8 \lastbox\unskip}
 \setbox@ne\hbox{\unhcopy8 \unskip\global\setbox\tw@\lastbox
 \unskip\global\setbox\thr@@\lastbox}
 \global\setbox7 \hbox{\unhbox\tw@\unskip}
 \ifinalign@\ifctagsplit@                                                   
  \gdef\split@{\hbox to\wd\thr@@{}&
   \vcenter{\vbox{\moveleft\wd\thr@@\boxz@}}}
 \else\gdef\split@{&\vbox{\moveleft\wd\thr@@\box9}\crcr
  \box\thr@@&\box7}\fi                                                      
 \else                                                                      
  \ifctagsplit@\gdef\split@{\vcenter{\boxz@}}\else
  \gdef\split@{\box9\crcr\hbox{\box\thr@@\box7}}\fi
 \fi
 \split@}                                                                   
\def\lendsplit@{\global\setbox9\vtop{\unvcopy\z@}
 \setbox@ne\vbox{\unvcopy\z@\global\setbox8\lastbox}
 \setbox@ne\hbox{\unhcopy8\unskip\setbox\tw@\lastbox
  \unskip\global\setbox\thr@@\lastbox}
 \ifinalign@\ifctagsplit@                                                   
  \gdef\split@{\hbox to\wd\thr@@{}&
  \vcenter{\vbox{\moveleft\wd\thr@@\box9}}}
  \else                                                                     
  \gdef\split@{\hbox to\wd\thr@@{}&\vbox{\moveleft\wd\thr@@\box9}}\fi
 \else
  \ifctagsplit@\gdef\split@{\vcenter{\box9}}\else
  \gdef\split@{\box9}\fi
 \fi\split@}
\def\outsplit@#1$${\align\insplit@#1\endalign$$}
\newdimen\multlinegap@
\multlinegap@1em
\newdimen\multlinetaggap@
\multlinetaggap@1em
\def\MultlineGap#1{\global\multlinegap@#1\relax}
\def\multlinegap#1{\RIfMIfI@\onlydmatherr@\multlinegap\else
 \multlinegap@#1\relax\fi\else\onlydmatherr@\multlinegap\fi}
\def\nomultlinegap{\multlinegap{\z@}}
\def\multline{\RIfMIfI@
 \DN@{\onlydmatherr@\multline}\else
 \DN@{\multline@}\fi\else
 \DN@{\onlydmatherr@\multline}\fi\next@}
\newif\iftagin@
\def\tagin@#1{\tagin@false\in@\tag{#1}\ifin@\tagin@true\fi}
\def\multline@#1$${\inany@true\vspace@\allowdisplaybreak@\displaybreak@
 \tagin@{#1}\iftagsleft@\DN@{\multline@l#1$$}\else
 \DN@{\multline@r#1$$}\fi\next@}
\newdimen\mwidth@
\def\rmmeasure@#1\endmultline{%
 \def\shoveleft##1{##1}\def\shoveright##1{##1}
 \setbox@ne\vbox{\Let@\halign{\setboxz@h
  {$\m@th\@lign\displaystyle{}##$}\global\mwidth@\wdz@
  \crcr#1\crcr}}}
\newdimen\mlineht@
\newif\ifzerocr@
\newif\ifonecr@
\def\lmmeasure@#1\endmultline{\global\zerocr@true\global\onecr@false
 \everycr{\noalign{\ifonecr@\global\onecr@false\fi
  \ifzerocr@\global\zerocr@false\global\onecr@true\fi}}
  \def\shoveleft##1{##1}\def\shoveright##1{##1}%
 \setbox@ne\vbox{\Let@\halign{\setboxz@h
  {$\m@th\@lign\displaystyle{}##$}\ifonecr@\global\mwidth@\wdz@
  \global\mlineht@\ht\z@\fi\crcr#1\crcr}}}
\newbox\mtagbox@
\newdimen\ltwidth@
\newdimen\rtwidth@
\def\multline@l#1$${\iftagin@\DN@{\lmultline@@#1$$}\else
 \DN@{\setbox\mtagbox@\null\ltwidth@\z@\rtwidth@\z@
  \lmultline@@@#1$$}\fi\next@}
\def\lmultline@@#1\endmultline\tag#2$${%
 \setbox\mtagbox@\hbox{\maketag@#2\maketag@}
 \lmmeasure@#1\endmultline\dimen@\mwidth@\advance\dimen@\wd\mtagbox@
 \advance\dimen@\multlinetaggap@                                            
 \ifdim\dimen@>\displaywidth\ltwidth@\z@\else\ltwidth@\wd\mtagbox@\fi       
 \lmultline@@@#1\endmultline$$}
\def\lmultline@@@{\displ@y
 \def\shoveright##1{##1\hfilneg\hskip\multlinegap@}%
 \def\shoveleft##1{\setboxz@h{$\m@th\displaystyle{}##1$}%
  \setbox@ne\hbox{$\m@th\displaystyle##1$}%
  \hfilneg
  \iftagin@
   \ifdim\ltwidth@>\z@\hskip\ltwidth@\hskip\multlinetaggap@\fi
  \else\hskip\multlinegap@\fi\hskip.5\wd@ne\hskip-.5\wdz@##1}
  \halign\bgroup\Let@\hbox to\displaywidth
   {\strut@$\m@th\displaystyle\hfil{}##\hfil$}\crcr
   \hfilneg                                                                 
   \iftagin@                                                                
    \ifdim\ltwidth@>\z@                                                     
     \box\mtagbox@\hskip\multlinetaggap@                                    
    \else
     \rlap{\vbox{\normalbaselines\hbox{\strut@\box\mtagbox@}%
     \vbox to\mlineht@{}}}\fi                                               
   \else\hskip\multlinegap@\fi}                                             
\def\multline@r#1$${\iftagin@\DN@{\rmultline@@#1$$}\else
 \DN@{\setbox\mtagbox@\null\ltwidth@\z@\rtwidth@\z@
  \rmultline@@@#1$$}\fi\next@}
\def\rmultline@@#1\endmultline\tag#2$${\ltwidth@\z@
 \setbox\mtagbox@\hbox{\maketag@#2\maketag@}%
 \rmmeasure@#1\endmultline\dimen@\mwidth@\advance\dimen@\wd\mtagbox@
 \advance\dimen@\multlinetaggap@
 \ifdim\dimen@>\displaywidth\rtwidth@\z@\else\rtwidth@\wd\mtagbox@\fi
 \rmultline@@@#1\endmultline$$}
\def\rmultline@@@{\displ@y
 \def\shoveright##1{##1\hfilneg\iftagin@\ifdim\rtwidth@>\z@
  \hskip\rtwidth@\hskip\multlinetaggap@\fi\else\hskip\multlinegap@\fi}%
 \def\shoveleft##1{\setboxz@h{$\m@th\displaystyle{}##1$}%
  \setbox@ne\hbox{$\m@th\displaystyle##1$}%
  \hfilneg\hskip\multlinegap@\hskip.5\wd@ne\hskip-.5\wdz@##1}%
 \halign\bgroup\Let@\hbox to\displaywidth
  {\strut@$\m@th\displaystyle\hfil{}##\hfil$}\crcr
 \hfilneg\hskip\multlinegap@}
\def\endmultline{\iftagsleft@\expandafter\lendmultline@\else
 \expandafter\rendmultline@\fi}
\def\lendmultline@{\hfilneg\hskip\multlinegap@\crcr\egroup}
\def\rendmultline@{\iftagin@                                                
 \ifdim\rtwidth@>\z@                                                        
  \hskip\multlinetaggap@\box\mtagbox@                                       
 \else\llap{\vtop{\normalbaselines\null\hbox{\strut@\box\mtagbox@}}}\fi     
 \else\hskip\multlinegap@\fi                                                
 \hfilneg\crcr\egroup}
\def\bmod{\mskip-\medmuskip\mkern5mu\mathbin{\fam\z@ mod}\penalty900
 \mkern5mu\mskip-\medmuskip}
\def\pmod#1{\allowbreak\ifinner\mkern8mu\else\mkern18mu\fi
 ({\fam\z@ mod}\,\,#1)}
\def\pod#1{\allowbreak\ifinner\mkern8mu\else\mkern18mu\fi(#1)}
\def\mod#1{\allowbreak\ifinner\mkern12mu\else\mkern18mu\fi{\fam\z@ mod}\,\,#1}
\message{continued fractions,}
\newcount\cfraccount@
\def\cfrac{\bgroup\bgroup\advance\cfraccount@\@ne\strut
 \iffalse{\fi\def\\{\over\displaystyle}\iffalse}\fi}
\def\lcfrac{\bgroup\bgroup\advance\cfraccount@\@ne\strut
 \iffalse{\fi\def\\{\hfill\over\displaystyle}\iffalse}\fi}
\def\rcfrac{\bgroup\bgroup\advance\cfraccount@\@ne\strut\hfill
 \iffalse{\fi\def\\{\over\displaystyle}\iffalse}\fi}
\def\gloop@#1\repeat{\gdef\body{#1}\iterate}
\def\endcfrac{\gloop@\ifnum\cfraccount@>\z@\global\advance\cfraccount@\m@ne
 \egroup\hskip-\nulldelimiterspace\egroup\repeat}
\message{compound symbols,}
\def\binrel@#1{\setboxz@h{\thinmuskip0mu
  \medmuskip\m@ne mu\thickmuskip\@ne mu$#1\m@th$}%
 \setbox@ne\hbox{\thinmuskip0mu\medmuskip\m@ne mu\thickmuskip
  \@ne mu${}#1{}\m@th$}%
 \setbox\tw@\hbox{\hskip\wd@ne\hskip-\wdz@}}
\def\overset#1\to#2{\binrel@{#2}\ifdim\wd\tw@<\z@
 \mathbin{\mathop{\kern\z@#2}\limits^{#1}}\else\ifdim\wd\tw@>\z@
 \mathrel{\mathop{\kern\z@#2}\limits^{#1}}\else
 {\mathop{\kern\z@#2}\limits^{#1}}{}\fi\fi}
\def\underset#1\to#2{\binrel@{#2}\ifdim\wd\tw@<\z@
 \mathbin{\mathop{\kern\z@#2}\limits_{#1}}\else\ifdim\wd\tw@>\z@
 \mathrel{\mathop{\kern\z@#2}\limits_{#1}}\else
 {\mathop{\kern\z@#2}\limits_{#1}}{}\fi\fi}
\def\oversetbrace#1\to#2{\overbrace{#2}^{#1}}
\def\undersetbrace#1\to#2{\underbrace{#2}_{#1}}
\def\sideset#1\and#2\to#3{%
 \setbox@ne\hbox{$\dsize{\vphantom{#3}}#1{#3}\m@th$}%
 \setbox\tw@\hbox{$\dsize{#3}#2\m@th$}%
 \hskip\wd@ne\hskip-\wd\tw@\mathop{\hskip\wd\tw@\hskip-\wd@ne
  {\vphantom{#3}}#1{#3}#2}}
\def\rightarrowfill@#1{\setboxz@h{$#1-\m@th$}\ht\z@\z@
  $#1\m@th\copy\z@\mkern-6mu\cleaders
  \hbox{$#1\mkern-2mu\box\z@\mkern-2mu$}\hfill
  \mkern-6mu\mathord\rightarrow$}
\def\leftarrowfill@#1{\setboxz@h{$#1-\m@th$}\ht\z@\z@
  $#1\m@th\mathord\leftarrow\mkern-6mu\cleaders
  \hbox{$#1\mkern-2mu\copy\z@\mkern-2mu$}\hfill
  \mkern-6mu\box\z@$}
\def\leftrightarrowfill@#1{\setboxz@h{$#1-\m@th$}\ht\z@\z@
  $#1\m@th\mathord\leftarrow\mkern-6mu\cleaders
  \hbox{$#1\mkern-2mu\box\z@\mkern-2mu$}\hfill
  \mkern-6mu\mathord\rightarrow$}
\def\overrightarrow{\mathpalette\overrightarrow@}
\def\overrightarrow@#1#2{\vbox{\ialign{##\crcr\rightarrowfill@#1\crcr
 \noalign{\kern-\ex@\nointerlineskip}$\m@th\hfil#1#2\hfil$\crcr}}}

\def\overleftarrow{\mathpalette\overleftarrow@}
\def\overleftarrow@#1#2{\vbox{\ialign{##\crcr\leftarrowfill@#1\crcr
 \noalign{\kern-\ex@\nointerlineskip}$\m@th\hfil#1#2\hfil$\crcr}}}
\def\overleftrightarrow{\mathpalette\overleftrightarrow@}
\def\overleftrightarrow@#1#2{\vbox{\ialign{##\crcr\leftrightarrowfill@#1\crcr
 \noalign{\kern-\ex@\nointerlineskip}$\m@th\hfil#1#2\hfil$\crcr}}}
\def\underrightarrow{\mathpalette\underrightarrow@}
\def\underrightarrow@#1#2{\vtop{\ialign{##\crcr$\m@th\hfil#1#2\hfil$\crcr
 \noalign{\nointerlineskip}\rightarrowfill@#1\crcr}}}

\def\underleftarrow{\mathpalette\underleftarrow@}
\def\underleftarrow@#1#2{\vtop{\ialign{##\crcr$\m@th\hfil#1#2\hfil$\crcr
 \noalign{\nointerlineskip}\leftarrowfill@#1\crcr}}}
\def\underleftrightarrow{\mathpalette\underleftrightarrow@}
\def\underleftrightarrow@#1#2{\vtop{\ialign{##\crcr$\m@th\hfil#1#2\hfil$\crcr
 \noalign{\nointerlineskip}\leftrightarrowfill@#1\crcr}}}
\message{various kinds of dots,}
\let\DOTSI\relax
\let\DOTSB\relax

\newif\ifmath@
{\uccode`7=`\\ \uccode`8=`m \uccode`9=`a \uccode`0=`t \uccode`!=`h
 \uppercase{\gdef\math@#1#2#3#4#5#6\math@{\global\math@false\ifx 7#1\ifx 8#2%
 \ifx 9#3\ifx 0#4\ifx !#5\xdef\meaning@{#6}\global\math@true\fi\fi\fi\fi\fi}}}
\newif\ifmathch@
{\uccode`7=`c \uccode`8=`h \uccode`9=`\"
 \uppercase{\gdef\mathch@#1#2#3#4#5#6\mathch@{\global\mathch@false
  \ifx 7#1\ifx 8#2\ifx 9#5\global\mathch@true\xdef\meaning@{9#6}\fi\fi\fi}}}
\newcount\classnum@
\def\getmathch@#1.#2\getmathch@{\classnum@#1 \divide\classnum@4096
 \ifcase\number\classnum@\or\or\gdef\thedots@{\dotsb@}\or
 \gdef\thedots@{\dotsb@}\fi}
\newif\ifmathbin@
{\uccode`4=`b \uccode`5=`i \uccode`6=`n
 \uppercase{\gdef\mathbin@#1#2#3{\relaxnext@
  \DNii@##1\mathbin@{\ifx\space@\next\global\mathbin@true\fi}%
 \global\mathbin@false\DN@##1\mathbin@{}%
 \ifx 4#1\ifx 5#2\ifx 6#3\DN@{\FN@\nextii@}\fi\fi\fi\next@}}}
\newif\ifmathrel@
{\uccode`4=`r \uccode`5=`e \uccode`6=`l
 \uppercase{\gdef\mathrel@#1#2#3{\relaxnext@
  \DNii@##1\mathrel@{\ifx\space@\next\global\mathrel@true\fi}%
 \global\mathrel@false\DN@##1\mathrel@{}%
 \ifx 4#1\ifx 5#2\ifx 6#3\DN@{\FN@\nextii@}\fi\fi\fi\next@}}}
\newif\ifmacro@
{\uccode`5=`m \uccode`6=`a \uccode`7=`c
 \uppercase{\gdef\macro@#1#2#3#4\macro@{\global\macro@false
  \ifx 5#1\ifx 6#2\ifx 7#3\global\macro@true
  \xdef\meaning@{\macro@@#4\macro@@}\fi\fi\fi}}}
\def\macro@@#1->#2\macro@@{#2}
\newif\ifDOTS@
\newcount\DOTSCASE@
{\uccode`6=`\\ \uccode`7=`D \uccode`8=`O \uccode`9=`T \uccode`0=`S
 \uppercase{\gdef\DOTS@#1#2#3#4#5{\global\DOTS@false\DN@##1\DOTS@{}%
  \ifx 6#1\ifx 7#2\ifx 8#3\ifx 9#4\ifx 0#5\let\next@\DOTS@@\fi\fi\fi\fi\fi
  \next@}}}
{\uccode`3=`B \uccode`4=`I \uccode`5=`X
 \uppercase{\gdef\DOTS@@#1{\relaxnext@
  \DNii@##1\DOTS@{\ifx\space@\next\global\DOTS@true\fi}%
  \DN@{\FN@\nextii@}%
  \ifx 3#1\global\DOTSCASE@\z@\else
  \ifx 4#1\global\DOTSCASE@\@ne\else
  \ifx 5#1\global\DOTSCASE@\tw@\else\DN@##1\DOTS@{}%
  \fi\fi\fi\next@}}}
\newif\ifnot@
{\uccode`5=`\\ \uccode`6=`n \uccode`7=`o \uccode`8=`t
 \uppercase{\gdef\not@#1#2#3#4{\relaxnext@
  \DNii@##1\not@{\ifx\space@\next\global\not@true\fi}%
 \global\not@false\DN@##1\not@{}%
 \ifx 5#1\ifx 6#2\ifx 7#3\ifx 8#4\DN@{\FN@\nextii@}\fi\fi\fi
 \fi\next@}}}
\newif\ifkeybin@
\def\keybin@{\keybin@true
 \ifx\next+\else\ifx\next=\else\ifx\next<\else\ifx\next>\else\ifx\next-\else
 \ifx\next*\else\ifx\next:\else\keybin@false\fi\fi\fi\fi\fi\fi\fi}
\def\dots{\RIfM@\expandafter\mdots@\else\expandafter\tdots@\fi}
\def\tdots@{\unskip\relaxnext@
 \DN@{$\m@th\mathinner{\ldotp\ldotp\ldotp}\,
   \ifx\next,\,$\else\ifx\next.\,$\else\ifx\next;\,$\else\ifx\next:\,$\else
   \ifx\next?\,$\else\ifx\next!\,$\else$ \fi\fi\fi\fi\fi\fi}%
 \ \FN@\next@}
\def\mdots@{\FN@\mdots@@}
\def\mdots@@{\gdef\thedots@{\dotso@}
 \ifx\next\boldkey\gdef\thedots@\boldkey{\boldkeydots@}\else                
 \ifx\next\boldsymbol\gdef\thedots@\boldsymbol{\boldsymboldots@}\else       
 \ifx,\next\gdef\thedots@{\dotsc}
 \else\ifx\not\next\gdef\thedots@{\dotsb@}
 \else\keybin@
 \ifkeybin@\gdef\thedots@{\dotsb@}
 \else\xdef\meaning@{\meaning\next..........}\xdef\meaning@@{\meaning@}
  \expandafter\math@\meaning@\math@
  \ifmath@
   \expandafter\mathch@\meaning@\mathch@
   \ifmathch@\expandafter\getmathch@\meaning@\getmathch@\fi                 
  \else\expandafter\macro@\meaning@@\macro@                                 
  \ifmacro@                                                                
   \expandafter\not@\meaning@\not@\ifnot@\gdef\thedots@{\dotsb@}
  \else\expandafter\DOTS@\meaning@\DOTS@
  \ifDOTS@
   \ifcase\number\DOTSCASE@\gdef\thedots@{\dotsb@}%
    \or\gdef\thedots@{\dotsi}\else\fi                                      
  \else\expandafter\math@\meaning@\math@                                   
  \ifmath@\expandafter\mathbin@\meaning@\mathbin@
  \ifmathbin@\gdef\thedots@{\dotsb@}
  \else\expandafter\mathrel@\meaning@\mathrel@
  \ifmathrel@\gdef\thedots@{\dotsb@}
  \fi\fi\fi\fi\fi\fi\fi\fi\fi\fi\fi\fi
 \thedots@}
\def\plainldots@{\mathinner{\ldotp\ldotp\ldotp}}
\def\plaincdots@{\mathinner{\cdotp\cdotp\cdotp}}
\def\dotsi{\!\plaincdots@}
\let\dotsb@\plaincdots@
\newif\ifextra@
\newif\ifrightdelim@
\def\rightdelim@{\global\rightdelim@true                                    
 \ifx\next)\else                                                            
 \ifx\next]\else
 \ifx\next\rbrack\else
 \ifx\next\}\else
 \ifx\next\rbrace\else
 \ifx\next\rangle\else
 \ifx\next\rceil\else
 \ifx\next\rfloor\else
 \ifx\next\rgroup\else
 \ifx\next\rmoustache\else
 \ifx\next\right\else
 \ifx\next\bigr\else
 \ifx\next\biggr\else
 \ifx\next\Bigr\else                                                        
 \ifx\next\Biggr\else\global\rightdelim@false
 \fi\fi\fi\fi\fi\fi\fi\fi\fi\fi\fi\fi\fi\fi\fi}
\def\extra@{%
 \global\extra@false\rightdelim@\ifrightdelim@\global\extra@true            
 \else\ifx\next$\global\extra@true                                          
 \else\xdef\meaning@{\meaning\next..........}
 \expandafter\macro@\meaning@\macro@\ifmacro@                               
 \expandafter\DOTS@\meaning@\DOTS@
 \ifDOTS@
 \ifnum\DOTSCASE@=\tw@\global\extra@true                                    
 \fi\fi\fi\fi\fi}
\newif\ifbold@
\def\dotso@{\relaxnext@
 \ifbold@
  \let\next\delayed@
  \DNii@{\extra@\plainldots@\ifextra@\,\fi}%
 \else
  \DNii@{\DN@{\extra@\plainldots@\ifextra@\,\fi}\FN@\next@}%
 \fi
 \nextii@}
\def\extrap@#1{%
 \ifx\next,\DN@{#1\,}\else
 \ifx\next;\DN@{#1\,}\else
 \ifx\next.\DN@{#1\,}\else\extra@
 \ifextra@\DN@{#1\,}\else
 \let\next@#1\fi\fi\fi\fi\next@}
\def\ldots{\DN@{\extrap@\plainldots@}%
 \FN@\next@}
\def\cdots{\DN@{\extrap@\plaincdots@}%
 \FN@\next@}

\def\dotsc{\relaxnext@
 \DN@{\ifx\next;\plainldots@\,\else
  \ifx\next.\plainldots@\,\else\extra@\plainldots@
  \ifextra@\,\fi\fi\fi}%
 \FN@\next@}
\def\cdot{\mathchar"2201 }

\def\mapsto{\DOTSB\mapstochar\rightarrow}

\message{special superscripts,}
\def\dddot#1{{\mathop{#1}\limits^{\vbox to-1.4\ex@{\kern-\tw@\ex@
 \hbox{\rm...}\vss}}}}
\def\ddddot#1{{\mathop{#1}\limits^{\vbox to-1.4\ex@{\kern-\tw@\ex@
 \hbox{\rm....}\vss}}}}
\def\sphat{^{\mathchoice{}{}%
 {\,\,\botsmash{\hbox{\lower4\ex@\hbox{$\m@th\widehat{\null}$}}}}%
 {\,\botsmash{\hbox{\lower3\ex@\hbox{$\m@th\hat{\null}$}}}}}}

\def\spacute{^{\!\botsmash{\hbox{\lower\@ne ex\hbox{\'{}}}}}}
\def\spgrave{^{\mathchoice{}{}{}{\!}%
 \botsmash{\hbox{\lower\@ne ex\hbox{\`{}}}}}}
\def\spdot{^{\hbox{\raise\ex@\hbox{\rm.}}}}
\def\spddot{^{\hbox{\raise\ex@\hbox{\rm..}}}}
\def\spdddot{^{\hbox{\raise\ex@\hbox{\rm...}}}}
\def\spddddot{^{\hbox{\raise\ex@\hbox{\rm....}}}}
\def\spbreve{^{\!\botsmash{\hbox{\lower4\ex@\hbox{\u{}}}}}}

\message{\string\text,}
\def\textonlyfont@#1#2{\def#1{\RIfM@
 \Err@{Use \string#1\space only in text}\else#2\fi}}
\textonlyfont@\rm\tenrm
\textonlyfont@\it\tenit
\textonlyfont@\sl\tensl
\textonlyfont@\bf\tenbf
\def\oldnos#1{\RIfM@{\mathcode`\,="013B \fam\@ne#1}\else
 \leavevmode\hbox{$\m@th\mathcode`\,="013B \fam\@ne#1$}\fi}
\def\text{\RIfM@\expandafter\text@\else\expandafter\text@@\fi}
\def\text@@#1{\leavevmode\hbox{#1}}
\def\mathhexbox@#1#2#3{\text{$\m@th\mathchar"#1#2#3$}}
\def\dag{{\mathhexbox@279}}
\def\ddag{{\mathhexbox@27A}}
\def\S{{\mathhexbox@278}}
\def\P{{\mathhexbox@27B}}
\newif\iffirstchoice@
\firstchoice@true
\def\text@#1{\mathchoice
 {\hbox{\everymath{\displaystyle}\def\textfonti{\the\textfont\@ne}%
  \def\textfontii{\the\textfont\tw@}\textdef@@ T#1}}
 {\hbox{\firstchoice@false
  \everymath{\textstyle}\def\textfonti{\the\textfont\@ne}%
  \def\textfontii{\the\textfont\tw@}\textdef@@ T#1}}
 {\hbox{\firstchoice@false
  \everymath{\scriptstyle}\def\textfonti{\the\scriptfont\@ne}%
  \def\textfontii{\the\scriptfont\tw@}\textdef@@ S\rm#1}}
 {\hbox{\firstchoice@false
  \everymath{\scriptscriptstyle}\def\textfonti
  {\the\scriptscriptfont\@ne}%
  \def\textfontii{\the\scriptscriptfont\tw@}\textdef@@ s\rm#1}}}
\def\textdef@@#1{\textdef@#1\rm\textdef@#1\bf\textdef@#1\sl\textdef@#1\it}
\def\rmfam{0}
\def\textdef@#1#2{%
 \DN@{\csname\expandafter\eat@\string#2fam\endcsname}%
 \if S#1\edef#2{\the\scriptfont\next@\relax}%
 \else\if s#1\edef#2{\the\scriptscriptfont\next@\relax}%
 \else\edef#2{\the\textfont\next@\relax}\fi\fi}
\scriptfont\itfam\tenit \scriptscriptfont\itfam\tenit
\scriptfont\slfam\tensl \scriptscriptfont\slfam\tensl
\newif\iftopfolded@
\newif\ifbotfolded@
\def\topfoldedtext{\topfolded@true\botfolded@false\foldedtext@}
\def\botfoldedtext{\botfolded@true\topfolded@false\foldedtext@}
\def\foldedtext{\topfolded@false\botfolded@false\foldedtext@}
\Invalid@\foldedwidth
\def\foldedtext@{\relaxnext@
 \DN@{\ifx\next\foldedwidth\let\next@\nextii@\else
  \DN@{\nextii@\foldedwidth{.3\hsize}}\fi\next@}%
 \DNii@\foldedwidth##1##2{\setbox\z@\vbox
  {\normalbaselines\hsize##1\relax
  \tolerance1600 \noindent\ignorespaces##2}\ifbotfolded@\boxz@\else
  \iftopfolded@\vtop{\unvbox\z@}\else\vcenter{\boxz@}\fi\fi}%
 \FN@\next@}
\message{math font commands,}
\def\bold{\RIfM@\expandafter\bold@\else
 \expandafter\nonmatherr@\expandafter\bold\fi}
\def\bold@#1{{\bold@@{#1}}}
\def\bold@@#1{\fam\bffam\relax#1}
\def\slanted{\RIfM@\expandafter\slanted@\else
 \expandafter\nonmatherr@\expandafter\slanted\fi}
\def\slanted@#1{{\slanted@@{#1}}}
\def\slanted@@#1{\fam\slfam\relax#1}
\def\roman{\RIfM@\expandafter\roman@\else
 \expandafter\nonmatherr@\expandafter\roman\fi}
\def\roman@#1{{\roman@@{#1}}}
\def\roman@@#1{\fam\rmfam\relax#1}
\def\italic{\RIfM@\expandafter\italic@\else
 \expandafter\nonmatherr@\expandafter\italic\fi}
\def\italic@#1{{\italic@@{#1}}}
\def\italic@@#1{\fam\itfam\relax#1}
\def\Cal{\RIfM@\expandafter\Cal@\else
 \expandafter\nonmatherr@\expandafter\Cal\fi}
\def\Cal@#1{{\Cal@@{#1}}}
\def\Cal@@#1{\noaccents@\fam\tw@#1}
\mathchardef\Gamma="0000
\mathchardef\Delta="0001
\mathchardef\Theta="0002
\mathchardef\Lambda="0003
\mathchardef\Xi="0004
\mathchardef\Pi="0005
\mathchardef\Sigma="0006
\mathchardef\Upsilon="0007
\mathchardef\Phi="0008
\mathchardef\Psi="0009
\mathchardef\Omega="000A
\mathchardef\varGamma="0100
\mathchardef\varDelta="0101
\mathchardef\varTheta="0102
\mathchardef\varLambda="0103
\mathchardef\varXi="0104
\mathchardef\varPi="0105
\mathchardef\varSigma="0106
\mathchardef\varUpsilon="0107
\mathchardef\varPhi="0108
\mathchardef\varPsi="0109
\mathchardef\varOmega="010A
\let\alloc@@\alloc@
\def\hexnumber@#1{\ifcase#1 0\or 1\or 2\or 3\or 4\or 5\or 6\or 7\or 8\or
 9\or A\or B\or C\or D\or E\or F\fi}
\def\loadmsam{%
 \font@\tenmsa=msam10
 \font@\sevenmsa=msam7
 \font@\fivemsa=msam5
 \alloc@@8\fam\chardef\sixt@@n\msafam
 \textfont\msafam=\tenmsa
 \scriptfont\msafam=\sevenmsa
 \scriptscriptfont\msafam=\fivemsa
 \edef\next{\hexnumber@\msafam}%
 \mathchardef\dabar@"0\next39
 \edef\dashrightarrow{\mathrel{\dabar@\dabar@\mathchar"0\next4B}}%
 \edef\dashleftarrow{\mathrel{\mathchar"0\next4C\dabar@\dabar@}}%
 \let\dasharrow\dashrightarrow
 \edef\ulcorner{\delimiter"4\next70\next70 }%
 \edef\urcorner{\delimiter"5\next71\next71 }%
 \edef\llcorner{\delimiter"4\next78\next78 }%
 \edef\lrcorner{\delimiter"5\next79\next79 }%
 \edef\yen{{\noexpand\mathhexbox@\next55}}%
 \edef\checkmark{{\noexpand\mathhexbox@\next58}}%
 \edef\circledR{{\noexpand\mathhexbox@\next72}}%
 \edef\maltese{{\noexpand\mathhexbox@\next7A}}%
 \global\let\loadmsam\empty}%
\def\loadmsbm{%
 \font@\tenmsb=msbm10 \font@\sevenmsb=msbm7 \font@\fivemsb=msbm5
 \alloc@@8\fam\chardef\sixt@@n\msbfam
 \textfont\msbfam=\tenmsb
 \scriptfont\msbfam=\sevenmsb \scriptscriptfont\msbfam=\fivemsb
 \global\let\loadmsbm\empty
 }
\def\widehat#1{\ifx\undefined\msbfam \DN@{362}%
  \else \setboxz@h{$\m@th#1$}%
    \edef\next@{\ifdim\wdz@>\tw@ em%
        \hexnumber@\msbfam 5B%
      \else 362\fi}\fi
  \mathaccent"0\next@{#1}}
\def\widetilde#1{\ifx\undefined\msbfam \DN@{365}%
  \else \setboxz@h{$\m@th#1$}%
    \edef\next@{\ifdim\wdz@>\tw@ em%
        \hexnumber@\msbfam 5D%
      \else 365\fi}\fi
  \mathaccent"0\next@{#1}}
\message{\string\newsymbol,}
\def\newsymbol#1#2#3#4#5{\define#1{}%
  \count@#2\relax \advance\count@\m@ne 
 \ifcase\count@
   \ifx\undefined\msafam\loadmsam\fi \let\next@\msafam
 \or \ifx\undefined\msbfam\loadmsbm\fi \let\next@\msbfam
 \else  \Err@{\Invalid@@\string\newsymbol}\let\next@\tw@\fi
 \mathchardef#1="#3\hexnumber@\next@#4#5\space}

\def\Bbb{\RIfM@\expandafter\Bbb@\else
 \expandafter\nonmatherr@\expandafter\Bbb\fi}
\def\Bbb@#1{{\Bbb@@{#1}}}
\def\Bbb@@#1{\noaccents@\fam\msbfam\relax#1}
\message{bold Greek and bold symbols,}
\def\loadbold{%
 \font@\tencmmib=cmmib10 \font@\sevencmmib=cmmib7 \font@\fivecmmib=cmmib5
 \skewchar\tencmmib'177 \skewchar\sevencmmib'177 \skewchar\fivecmmib'177
 \alloc@@8\fam\chardef\sixt@@n\cmmibfam
 \textfont\cmmibfam\tencmmib
 \scriptfont\cmmibfam\sevencmmib \scriptscriptfont\cmmibfam\fivecmmib
 \font@\tencmbsy=cmbsy10 \font@\sevencmbsy=cmbsy7 \font@\fivecmbsy=cmbsy5
 \skewchar\tencmbsy'60 \skewchar\sevencmbsy'60 \skewchar\fivecmbsy'60
 \alloc@@8\fam\chardef\sixt@@n\cmbsyfam
 \textfont\cmbsyfam\tencmbsy
 \scriptfont\cmbsyfam\sevencmbsy \scriptscriptfont\cmbsyfam\fivecmbsy
 \let\loadbold\empty
}
\def\boldnotloaded#1{\Err@{\ifcase#1\or First\else Second\fi
       bold symbol font not loaded}}
\def\mathchari@#1#2#3{\ifx\undefined\cmmibfam
    \boldnotloaded@\@ne
  \else\mathchar"#1\hexnumber@\cmmibfam#2#3\space \fi}
\def\mathcharii@#1#2#3{\ifx\undefined\cmbsyfam
    \boldnotloaded\tw@
  \else \mathchar"#1\hexnumber@\cmbsyfam#2#3\space\fi}
\edef\bffam@{\hexnumber@\bffam}
\def\boldkey#1{\ifcat\noexpand#1A%
  \ifx\undefined\cmmibfam \boldnotloaded\@ne
  \else {\fam\cmmibfam#1}\fi
 \else
 \ifx#1!\mathchar"5\bffam@21 \else
 \ifx#1(\mathchar"4\bffam@28 \else\ifx#1)\mathchar"5\bffam@29 \else
 \ifx#1+\mathchar"2\bffam@2B \else\ifx#1:\mathchar"3\bffam@3A \else
 \ifx#1;\mathchar"6\bffam@3B \else\ifx#1=\mathchar"3\bffam@3D \else
 \ifx#1?\mathchar"5\bffam@3F \else\ifx#1[\mathchar"4\bffam@5B \else
 \ifx#1]\mathchar"5\bffam@5D \else
 \ifx#1,\mathchari@63B \else
 \ifx#1-\mathcharii@200 \else
 \ifx#1.\mathchari@03A \else
 \ifx#1/\mathchari@03D \else
 \ifx#1<\mathchari@33C \else
 \ifx#1>\mathchari@33E \else
 \ifx#1*\mathcharii@203 \else
 \ifx#1|\mathcharii@06A \else
 \ifx#10\bold0\else\ifx#11\bold1\else\ifx#12\bold2\else\ifx#13\bold3\else
 \ifx#14\bold4\else\ifx#15\bold5\else\ifx#16\bold6\else\ifx#17\bold7\else
 \ifx#18\bold8\else\ifx#19\bold9\else
  \Err@{\string\boldkey\space can't be used with #1}%
 \fi\fi\fi\fi\fi\fi\fi\fi\fi\fi\fi\fi\fi\fi\fi
 \fi\fi\fi\fi\fi\fi\fi\fi\fi\fi\fi\fi\fi\fi}
\def\boldsymbol#1{%
 \DN@{\Err@{You can't use \string\boldsymbol\space with \string#1}#1}%
 \ifcat\noexpand#1A%
   \let\next@\relax
   \ifx\undefined\cmmibfam \boldnotloaded\@ne
   \else {\fam\cmmibfam#1}\fi
 \else
  \xdef\meaning@{\meaning#1.........}%
  \expandafter\math@\meaning@\math@
  \ifmath@
   \expandafter\mathch@\meaning@\mathch@
   \ifmathch@
    \expandafter\boldsymbol@@\meaning@\boldsymbol@@
   \fi
  \else
   \expandafter\macro@\meaning@\macro@
   \expandafter\delim@\meaning@\delim@
   \ifdelim@
    \expandafter\delim@@\meaning@\delim@@
   \else
    \boldsymbol@{#1}%
   \fi
  \fi
 \fi
 \next@}
\def\mathhexboxii@#1#2{\ifx\undefined\cmbsyfam
    \boldnotloaded\tw@
  \else \mathhexbox@{\hexnumber@\cmbsyfam}{#1}{#2}\fi}
\def\boldsymbol@#1{\let\next@\relax\let\next#1%
 \ifx\next\cdot\mathcharii@201 \else
 \ifx\next\prime{{\null\mathcharii@030 \null}}\else
 \ifx\next\lbrack\mathchar"4\bffam@5B \else
 \ifx\next\rbrack\mathchar"5\bffam@5D \else
 \ifx\next\{\mathcharii@466 \else
 \ifx\next\lbrace\mathcharii@466 \else
 \ifx\next\}\mathcharii@567 \else
 \ifx\next\rbrace\mathcharii@567 \else
 \ifx\next\surd{{\mathcharii@170}}\else
 \ifx\next\S{{\mathhexboxii@78}}\else
 \ifx\next\P{{\mathhexboxii@7B}}\else
 \ifx\next\dag{{\mathhexboxii@79}}\else
 \ifx\next\ddag{{\mathhexboxii@7A}}\else
 \DN@{\Err@{You can't use \string\boldsymbol\space with \string#1}#1}%
 \fi\fi\fi\fi\fi\fi\fi\fi\fi\fi\fi\fi\fi}
\def\boldsymbol@@#1.#2\boldsymbol@@{\classnum@#1 \count@@@\classnum@        
 \divide\classnum@4096 \count@\classnum@                                    
 \multiply\count@4096 \advance\count@@@-\count@ \count@@\count@@@           
 \divide\count@@@\@cclvi \count@\count@@                                    
 \multiply\count@@@\@cclvi \advance\count@@-\count@@@                       
 \divide\count@@@\@cclvi                                                    
 \multiply\classnum@4096 \advance\classnum@\count@@                         
 \ifnum\count@@@=\z@                                                        
  \count@"\bffam@ \multiply\count@\@cclvi
  \advance\classnum@\count@
  \DN@{\mathchar\number\classnum@}%
 \else
  \ifnum\count@@@=\@ne                                                      
   \ifx\undefined\cmmibfam \DN@{\boldnotloaded\@ne}%
   \else \count@\cmmibfam \multiply\count@\@cclvi
     \advance\classnum@\count@
     \DN@{\mathchar\number\classnum@}\fi
  \else
   \ifnum\count@@@=\tw@                                                    
     \ifx\undefined\cmbsyfam
       \DN@{\boldnotloaded\tw@}%
     \else
       \count@\cmbsyfam \multiply\count@\@cclvi
       \advance\classnum@\count@
       \DN@{\mathchar\number\classnum@}%
     \fi
  \fi
 \fi
\fi}
\newif\ifdelim@
\newcount\delimcount@
{\uccode`6=`\\ \uccode`7=`d \uccode`8=`e \uccode`9=`l
 \uppercase{\gdef\delim@#1#2#3#4#5\delim@
  {\delim@false\ifx 6#1\ifx 7#2\ifx 8#3\ifx 9#4\delim@true
   \xdef\meaning@{#5}\fi\fi\fi\fi}}}
\def\delim@@#1"#2#3#4#5#6\delim@@{\if#32%
\let\next@\relax
 \ifx\undefined\cmbsyfam \boldnotloaded\@ne
 \else \mathcharii@#2#4#5\space \fi\fi}
\def\vert{\delimiter"026A30C }
\def\Vert{\delimiter"026B30D }
\let\|\Vert
\def\backslash{\delimiter"026E30F }
\def\boldkeydots@#1{\bold@true\let\next=#1\let\delayed@=#1\mdots@@
 \boldkey#1\bold@false}  
\def\boldsymboldots@#1{\bold@true\let\next#1\let\delayed@#1\mdots@@
 \boldsymbol#1\bold@false}
\message{Euler fonts,}

\def\frak{\mathfont@\frak}

\def\loadmathfont#1{%
   \expandafter\font@\csname ten#1\endcsname=#110
   \expandafter\font@\csname seven#1\endcsname=#17
   \expandafter\font@\csname five#1\endcsname=#15
   \edef\next{\noexpand\alloc@@8\fam\chardef\sixt@@n
     \expandafter\noexpand\csname#1fam\endcsname}%
   \next
   \textfont\csname#1fam\endcsname \csname ten#1\endcsname
   \scriptfont\csname#1fam\endcsname \csname seven#1\endcsname
   \scriptscriptfont\csname#1fam\endcsname \csname five#1\endcsname
   \expandafter\def\csname #1\expandafter\endcsname\expandafter{%
      \expandafter\mathfont@\csname#1\endcsname}%
 \expandafter\gdef\csname load#1\endcsname{}%
}
\def\mathfont@#1{\RIfM@\expandafter\mathfont@@\expandafter#1\else
  \expandafter\nonmatherr@\expandafter#1\fi}
\def\mathfont@@#1#2{{\mathfont@@@#1{#2}}}
\def\mathfont@@@#1#2{\noaccents@
   \fam\csname\expandafter\eat@\string#1fam\endcsname
   \relax#2}
\message{math accents,}
\def\accentclass@{7}
\def\noaccents@{\def\accentclass@{0}}
\def\makeacc@#1#2{\def#1{\mathaccent"\accentclass@#2 }}
\makeacc@\hat{05E}
\makeacc@\check{014}
\makeacc@\tilde{07E}
\makeacc@\acute{013}
\makeacc@\grave{012}
\makeacc@\dot{05F}
\makeacc@\ddot{07F}
\makeacc@\breve{015}
\makeacc@\bar{016}

\newcount\skewcharcount@
\newcount\familycount@
\def\theskewchar@{\familycount@\@ne
 \global\skewcharcount@\the\skewchar\textfont\@ne                           
 \ifnum\fam>\m@ne\ifnum\fam<16
  \global\familycount@\the\fam\relax
  \global\skewcharcount@\the\skewchar\textfont\the\fam\relax\fi\fi          
 \ifnum\skewcharcount@>\m@ne
  \ifnum\skewcharcount@<128
  \multiply\familycount@256
  \global\advance\skewcharcount@\familycount@
  \global\advance\skewcharcount@28672
  \mathchar\skewcharcount@\else
  \global\skewcharcount@\m@ne\fi\else
 \global\skewcharcount@\m@ne\fi}                                            
\newcount\pointcount@
\def\getpoints@#1.#2\getpoints@{\pointcount@#1 }
\newdimen\accentdimen@
\newcount\accentmu@
\def\dimentomu@{\multiply\accentdimen@ 100
 \expandafter\getpoints@\the\accentdimen@\getpoints@
 \multiply\pointcount@18
 \divide\pointcount@\@m
 \global\accentmu@\pointcount@}
\def\Makeacc@#1#2{\def#1{\RIfM@\DN@{\mathaccent@
 {"\accentclass@#2 }}\else\DN@{\nonmatherr@{#1}}\fi\next@}}
\def\unbracefonts@{\let\Cal@\Cal@@\let\roman@\roman@@\let\bold@\bold@@
 \let\slanted@\slanted@@}
\def\mathaccent@#1#2{\ifnum\fam=\m@ne\xdef\thefam@{1}\else
 \xdef\thefam@{\the\fam}\fi                                                 
 \accentdimen@\z@                                                           
 \setboxz@h{\unbracefonts@$\m@th\fam\thefam@\relax#2$}
 \ifdim\accentdimen@=\z@\DN@{\mathaccent#1{#2}}
  \setbox@ne\hbox{\unbracefonts@$\m@th\fam\thefam@\relax#2\theskewchar@$}
  \setbox\tw@\hbox{$\m@th\ifnum\skewcharcount@=\m@ne\else
   \mathchar\skewcharcount@\fi$}
  \global\accentdimen@\wd@ne\global\advance\accentdimen@-\wdz@
  \global\advance\accentdimen@-\wd\tw@                                     
  \global\multiply\accentdimen@\tw@
  \dimentomu@\global\advance\accentmu@\@ne                                 
 \else\DN@{{\mathaccent#1{#2\mkern\accentmu@ mu}%
    \mkern-\accentmu@ mu}{}}\fi                                             
 \next@}\Makeacc@\Hat{05E}
\Makeacc@\Check{014}
\Makeacc@\Tilde{07E}
\Makeacc@\Acute{013}
\Makeacc@\Grave{012}
\Makeacc@\Dot{05F}
\Makeacc@\Ddot{07F}
\Makeacc@\Breve{015}
\Makeacc@\Bar{016}
\def\Vec{\RIfM@\DN@{\mathaccent@{"017E }}\else
 \DN@{\nonmatherr@\Vec}\fi\next@}
\def\accentedsymbol#1#2{\csname newbox\expandafter\endcsname
  \csname\expandafter\eat@\string#1@box\endcsname
 \expandafter\setbox\csname\expandafter\eat@
  \string#1@box\endcsname\hbox{$\m@th#2$}\define
  #1{\copy\csname\expandafter\eat@\string#1@box\endcsname{}}}
\message{roots,}
\def\sqrt#1{\radical"270370 {#1}}
\let\underline@\underline
\let\overline@\overline
\def\underline#1{\underline@{#1}}
\def\overline#1{\overline@{#1}}
\Invalid@\leftroot
\Invalid@\uproot
\newcount\uproot@
\newcount\leftroot@
\def\root{\relaxnext@
  \DN@{\ifx\next\uproot\let\next@\nextii@\else
   \ifx\next\leftroot\let\next@\nextiii@\else
   \let\next@\plainroot@\fi\fi\next@}%
  \DNii@\uproot##1{\uproot@##1\relax\FN@\nextiv@}%
  \def\nextiv@{\ifx\next\space@\DN@. {\FN@\nextv@}\else
   \DN@.{\FN@\nextv@}\fi\next@.}%
  \def\nextv@{\ifx\next\leftroot\let\next@\nextvi@\else
   \let\next@\plainroot@\fi\next@}%
  \def\nextvi@\leftroot##1{\leftroot@##1\relax\plainroot@}%
   \def\nextiii@\leftroot##1{\leftroot@##1\relax\FN@\nextvii@}%
  \def\nextvii@{\ifx\next\space@
   \DN@. {\FN@\nextviii@}\else
   \DN@.{\FN@\nextviii@}\fi\next@.}%
  \def\nextviii@{\ifx\next\uproot\let\next@\nextix@\else
   \let\next@\plainroot@\fi\next@}%
  \def\nextix@\uproot##1{\uproot@##1\relax\plainroot@}%
  \bgroup\uproot@\z@\leftroot@\z@\FN@\next@}
\def\plainroot@#1\of#2{\setbox\rootbox\hbox{$\m@th\scriptscriptstyle{#1}$}%
 \mathchoice{\r@@t\displaystyle{#2}}{\r@@t\textstyle{#2}}
 {\r@@t\scriptstyle{#2}}{\r@@t\scriptscriptstyle{#2}}\egroup}
\def\r@@t#1#2{\setboxz@h{$\m@th#1\sqrt{#2}$}%
 \dimen@\ht\z@\advance\dimen@-\dp\z@
 \setbox@ne\hbox{$\m@th#1\mskip\uproot@ mu$}\advance\dimen@ 1.667\wd@ne
 \mkern-\leftroot@ mu\mkern5mu\raise.6\dimen@\copy\rootbox
 \mkern-10mu\mkern\leftroot@ mu\boxz@}
\def\boxed#1{\setboxz@h{$\m@th\displaystyle{#1}$}\dimen@.4\ex@
 \advance\dimen@3\ex@\advance\dimen@\dp\z@
 \hbox{\lower\dimen@\hbox{%
 \vbox{\hrule height.4\ex@
 \hbox{\vrule width.4\ex@\hskip3\ex@\vbox{\vskip3\ex@\boxz@\vskip3\ex@}%
 \hskip3\ex@\vrule width.4\ex@}\hrule height.4\ex@}%
 }}}
\message{commutative diagrams,}
\let\ampersand@\relax
\newdimen\minaw@
\minaw@11.11128\ex@
\newdimen\minCDaw@
\minCDaw@2.5pc
\def\minCDarrowwidth#1{\RIfMIfI@\onlydmatherr@\minCDarrowwidth
 \else\minCDaw@#1\relax\fi\else\onlydmatherr@\minCDarrowwidth\fi}
\newif\ifCD@
\def\CD{\bgroup\vspace@\relax\let\ampersand@&\iffalse}\fi
 \CD@true\vcenter\bgroup\Let@\tabskip\z@skip\baselineskip20\ex@
 \lineskip3\ex@\lineskiplimit3\ex@\halign\bgroup
 &\hfill$\m@th##$\hfill\crcr}
\def\endCD{\crcr\egroup\egroup\egroup}
\newdimen\bigaw@
\atdef@>#1>#2>{\ampersand@                                                  
 \setboxz@h{$\m@th\ssize\;{#1}\;\;$}
 \setbox@ne\hbox{$\m@th\ssize\;{#2}\;\;$}
 \setbox\tw@\hbox{$\m@th#2$}
 \ifCD@\global\bigaw@\minCDaw@\else\global\bigaw@\minaw@\fi                 
 \ifdim\wdz@>\bigaw@\global\bigaw@\wdz@\fi
 \ifdim\wd@ne>\bigaw@\global\bigaw@\wd@ne\fi                                
 \ifCD@\enskip\fi                                                           
 \ifdim\wd\tw@>\z@
  \mathrel{\mathop{\hbox to\bigaw@{\rightarrowfill@\displaystyle}}%
    \limits^{#1}_{#2}}
 \else\mathrel{\mathop{\hbox to\bigaw@{\rightarrowfill@\displaystyle}}%
    \limits^{#1}}\fi                                                        
 \ifCD@\enskip\fi                                                          
 \ampersand@}                                                              
\atdef@<#1<#2<{\ampersand@\setboxz@h{$\m@th\ssize\;\;{#1}\;$}%
 \setbox@ne\hbox{$\m@th\ssize\;\;{#2}\;$}\setbox\tw@\hbox{$\m@th#2$}%
 \ifCD@\global\bigaw@\minCDaw@\else\global\bigaw@\minaw@\fi
 \ifdim\wdz@>\bigaw@\global\bigaw@\wdz@\fi
 \ifdim\wd@ne>\bigaw@\global\bigaw@\wd@ne\fi
 \ifCD@\enskip\fi
 \ifdim\wd\tw@>\z@
  \mathrel{\mathop{\hbox to\bigaw@{\leftarrowfill@\displaystyle}}%
       \limits^{#1}_{#2}}\else
  \mathrel{\mathop{\hbox to\bigaw@{\leftarrowfill@\displaystyle}}%
       \limits^{#1}}\fi
 \ifCD@\enskip\fi\ampersand@}
\begingroup
 \catcode`\~=\active \lccode`\~=`\@
 \lowercase{%
  \global\atdef@)#1)#2){~>#1>#2>}
  \global\atdef@(#1(#2({~<#1<#2<}}
\endgroup
\atdef@ A#1A#2A{\llap{$\m@th\vcenter{\hbox
 {$\ssize#1$}}$}\Big\uparrow\rlap{$\m@th\vcenter{\hbox{$\ssize#2$}}$}&&}
\atdef@ V#1V#2V{\llap{$\m@th\vcenter{\hbox
 {$\ssize#1$}}$}\Big\downarrow\rlap{$\m@th\vcenter{\hbox{$\ssize#2$}}$}&&}
\atdef@={&\enskip\mathrel
 {\vbox{\hrule width\minCDaw@\vskip3\ex@\hrule width
 \minCDaw@}}\enskip&}
\atdef@|{\Big\Vert&&}
\atdef@\vert{\Big\Vert&&}
\def\pretend#1\haswidth#2{\setboxz@h{$\m@th\scriptstyle{#2}$}\hbox
 to\wdz@{\hfill$\m@th\scriptstyle{#1}$\hfill}}
\message{poor man's bold,}
\def\pmb{\RIfM@\expandafter\mathpalette\expandafter\pmb@\else
 \expandafter\pmb@@\fi}
\def\pmb@@#1{\leavevmode\setboxz@h{#1}%
   \dimen@-\wdz@
   \kern-.5\ex@\copy\z@
   \kern\dimen@\kern.25\ex@\raise.4\ex@\copy\z@
   \kern\dimen@\kern.25\ex@\box\z@
}
\def\binrel@@#1{\ifdim\wd2<\z@\mathbin{#1}\else\ifdim\wd\tw@>\z@
 \mathrel{#1}\else{#1}\fi\fi}
\newdimen\pmbraise@
\def\pmb@#1#2{\setbox\thr@@\hbox{$\m@th#1{#2}$}%
 \setbox4\hbox{$\m@th#1\mkern.5mu$}\pmbraise@\wd4\relax
 \binrel@{#2}%
 \dimen@-\wd\thr@@
   \binrel@@{%
   \mkern-.8mu\copy\thr@@
   \kern\dimen@\mkern.4mu\raise\pmbraise@\copy\thr@@
   \kern\dimen@\mkern.4mu\box\thr@@
}}
\def\documentstyle#1{\W@{}\input #1.sty\relax}
\message{syntax check,}
\font\dummyft@=dummy
\fontdimen1 \dummyft@=\z@
\fontdimen2 \dummyft@=\z@
\fontdimen3 \dummyft@=\z@
\fontdimen4 \dummyft@=\z@
\fontdimen5 \dummyft@=\z@
\fontdimen6 \dummyft@=\z@
\fontdimen7 \dummyft@=\z@
\fontdimen8 \dummyft@=\z@
\fontdimen9 \dummyft@=\z@
\fontdimen10 \dummyft@=\z@
\fontdimen11 \dummyft@=\z@
\fontdimen12 \dummyft@=\z@
\fontdimen13 \dummyft@=\z@
\fontdimen14 \dummyft@=\z@
\fontdimen15 \dummyft@=\z@
\fontdimen16 \dummyft@=\z@
\fontdimen17 \dummyft@=\z@
\fontdimen18 \dummyft@=\z@
\fontdimen19 \dummyft@=\z@
\fontdimen20 \dummyft@=\z@
\fontdimen21 \dummyft@=\z@
\fontdimen22 \dummyft@=\z@
\def\fontlist@{\\{\tenrm}\\{\sevenrm}\\{\fiverm}\\{\teni}\\{\seveni}%
 \\{\fivei}\\{\tensy}\\{\sevensy}\\{\fivesy}\\{\tenex}\\{\tenbf}\\{\sevenbf}%
 \\{\fivebf}\\{\tensl}\\{\tenit}}
\def\font@#1=#2 {\rightappend@#1\to\fontlist@\font#1=#2 }
\def\dodummy@{{\def\\##1{\global\let##1\dummyft@}\fontlist@}}
\def\nopages@{\output{\setbox\z@\box\@cclv \deadcycles\z@}%
 \alloc@5\toks\toksdef\@cclvi\output}
\let\galleys\nopages@
\newif\ifsyntax@
\newcount\countxviii@
\def\syntax{\syntax@true\dodummy@\countxviii@\count18
 \loop\ifnum\countxviii@>\m@ne\textfont\countxviii@=\dummyft@
 \scriptfont\countxviii@=\dummyft@\scriptscriptfont\countxviii@=\dummyft@
 \advance\countxviii@\m@ne\repeat                                           
 \dummyft@\tracinglostchars\z@\nopages@\frenchspacing\hbadness\@M}
\def\first@#1#2\end{#1}
\def\printoptions{\W@{Do you want S(yntax check),
  G(alleys) or P(ages)?}%
 \message{Type S, G or P, followed by <return>: }%
 \begingroup 
 \endlinechar\m@ne 
 \read\m@ne to\ans@
 \edef\ans@{\uppercase{\def\noexpand\ans@{%
   \expandafter\first@\ans@ P\end}}}%
 \expandafter\endgroup\ans@
 \if\ans@ P
 \else \if\ans@ S\syntax
 \else \if\ans@ G\galleys
 \else\message{? Unknown option: \ans@; using the `pages' option.}%
 \fi\fi\fi}
\def\alloc@#1#2#3#4#5{\global\advance\count1#1by\@ne
 \ch@ck#1#4#2\allocationnumber=\count1#1
 \global#3#5=\allocationnumber
 \ifalloc@\wlog{\string#5=\string#2\the\allocationnumber}\fi}
\def\document{\def\alloclist@{}\def\fontlist@{}}
\let\enddocument\bye

\let\proclaim\undefined
\let\footnote\undefined
\let\=\undefined
\let\>\undefined

\catcode`\@=\active
\message{... finished}

\expandafter\ifx\csname mathdefs.tex\endcsname\relax
  \expandafter\gdef\csname mathdefs.tex\endcsname{}
\else \message{Hey!  Apparently you were trying to
  \string\input{mathdefs.tex} twice.   This does not make sense.} 
\errmessage{Please edit your file (probably \jobname.tex) and remove
any duplicate ``\string\input'' lines}\endinput\fi




\catcode`\X=12\catcode`\@=11

\def\n@wcount{\alloc@0\count\countdef\insc@unt}
\def\n@wwrite{\alloc@7\write\chardef\sixt@@n}
\def\n@wread{\alloc@6\read\chardef\sixt@@n}
\def\r@s@t{\relax}\def\v@idline{\par}\def\@mputate#1/{#1}
\def\l@c@l#1X{\firstpart.#1}\def\gl@b@l#1X{#1}\def\t@d@l#1X{{}}

\def\crossrefs#1{\ifx\all#1\let\tr@ce=\all\else\def\tr@ce{#1,}\fi
   \n@wwrite\cit@tionsout\openout\cit@tionsout=\jobname.cit 
   \write\cit@tionsout{\tr@ce}\expandafter\setfl@gs\tr@ce,}
\def\setfl@gs#1,{\def\@{#1}\ifx\@\empty\let\next=\relax
   \else\let\next=\setfl@gs\expandafter\xdef
   \csname#1tr@cetrue\endcsname{}\fi\next}
\def\m@ketag#1#2{\expandafter\n@wcount\csname#2tagno\endcsname
     \csname#2tagno\endcsname=0\let\tail=\all\xdef\all{\tail#2,}
   \ifx#1\l@c@l\let\tail=\r@s@t\xdef\r@s@t{\csname#2tagno\endcsname=0\tail}\fi
   \expandafter\gdef\csname#2cite\endcsname##1{\expandafter
     \ifx\csname#2tag##1\endcsname\relax?\else\csname#2tag##1\endcsname\fi
     \expandafter\ifx\csname#2tr@cetrue\endcsname\relax\else
     \write\cit@tionsout{#2tag ##1 cited on page \folio.}\fi}
   \expandafter\gdef\csname#2page\endcsname##1{\expandafter
     \ifx\csname#2page##1\endcsname\relax?\else\csname#2page##1\endcsname\fi
     \expandafter\ifx\csname#2tr@cetrue\endcsname\relax\else
     \write\cit@tionsout{#2tag ##1 cited on page \folio.}\fi}
   \expandafter\gdef\csname#2tag\endcsname##1{\expandafter
      \ifx\csname#2check##1\endcsname\relax
      \expandafter\xdef\csname#2check##1\endcsname{}%
      \else\immediate\write16{Warning: #2tag ##1 used more than once.}\fi
      \multit@g{#1}{#2}##1/X%
      \write\t@gsout{#2tag ##1 assigned number \csname#2tag##1\endcsname\space
      on page \number\count0.}%
   \csname#2tag##1\endcsname}}

\def\multit@g#1#2#3/#4X{\def\t@mp{#4}\ifx\t@mp\empty%
      \global\advance\csname#2tagno\endcsname by 1 
      \expandafter\xdef\csname#2tag#3\endcsname
      {#1\number\csname#2tagno\endcsnameX}%
   \else\expandafter\ifx\csname#2last#3\endcsname\relax
      \expandafter\n@wcount\csname#2last#3\endcsname
      \global\advance\csname#2tagno\endcsname by 1 
      \expandafter\xdef\csname#2tag#3\endcsname
      {#1\number\csname#2tagno\endcsnameX}
      \write\t@gsout{#2tag #3 assigned number \csname#2tag#3\endcsname\space
      on page \number\count0.}\fi
   \global\advance\csname#2last#3\endcsname by 1
   \def\t@mp{\expandafter\xdef\csname#2tag#3/}%
   \expandafter\t@mp\@mputate#4\endcsname
   {\csname#2tag#3\endcsname\lastpart{\csname#2last#3\endcsname}}\fi}
\def\t@gs#1{\def\all{}\m@ketag#1e\m@ketag#1s\m@ketag\t@d@l p
\let\realscite\scite
\let\realstag\stag
   \m@ketag\gl@b@l r \n@wread\t@gsin
   \openin\t@gsin=\jobname.tgs \re@der \closein\t@gsin
   \n@wwrite\t@gsout\openout\t@gsout=\jobname.tgs }
\outer\def\localtags{\t@gs\l@c@l}
\outer\def\globaltags{\t@gs\gl@b@l}
\outer\def\newlocaltag#1{\m@ketag\l@c@l{#1}}
\outer\def\newglobaltag#1{\m@ketag\gl@b@l{#1}}

\newif\ifpr@ 
\def\m@kecs #1tag #2 assigned number #3 on page #4.%
   {\expandafter\gdef\csname#1tag#2\endcsname{#3}
   \expandafter\gdef\csname#1page#2\endcsname{#4}
   \ifpr@\expandafter\xdef\csname#1check#2\endcsname{}\fi}
\def\re@der{\ifeof\t@gsin\let\next=\relax\else
   \read\t@gsin to\t@gline\ifx\t@gline\v@idline\else
   \expandafter\m@kecs \t@gline\fi\let \next=\re@der\fi\next}
\def\pretags#1{\pr@true\pret@gs#1,,}
\def\pret@gs#1,{\def\@{#1}\ifx\@\empty\let\n@xtfile=\relax
   \else\let\n@xtfile=\pret@gs \openin\t@gsin=#1.tgs \message{#1} \re@der 
   \closein\t@gsin\fi \n@xtfile}

\newcount\sectno\sectno=0\newcount\subsectno\subsectno=0
\newif\ifultr@local \def\ultralocal{\ultr@localtrue}
\def\firstpart{\number\sectno}
\def\lastpart#1{\ifcase#1 \or a\or b\or c\or d\or e\or f\or g\or h\or 
   i\or k\or l\or m\or n\or o\or p\or q\or r\or s\or t\or u\or v\or w\or 
   x\or y\or z \fi}

\def\resetall{\global\advance\sectno by 1\subsectno=0
   \gdef\firstpart{\number\sectno}\r@s@t}
\def\resetsub{\global\advance\subsectno by 1
   \gdef\firstpart{\number\sectno.\number\subsectno}\r@s@t}
\def\newsection#1\par{\resetall\vskip0pt plus.3\vsize\penalty-250
   \vskip0pt plus-.3\vsize\bigskip\bigskip
   \message{#1}\leftline{\bf#1}\nobreak\bigskip}
\def\subsection#1\par{\ifultr@local\resetsub\fi
   \vskip0pt plus.2\vsize\penalty-250\vskip0pt plus-.2\vsize
   \bigskip\smallskip\message{#1}\leftline{\bf#1}\nobreak\medskip}


\newdimen\marginshift

\newdimen\margindelta
\newdimen\marginmax
\newdimen\marginmin

\def\margininit{       
\marginmax=3 true cm                  
				      
\margindelta=0.1 true cm              
\marginmin=0.1true cm                 
\marginshift=\marginmin
}    

\def\t@gsjj#1,{\def\@{#1}\ifx\@\empty\let\next=\relax\else\let\next=\t@gsjj
   \def\@@{p}\ifx\@\@@\else
   \expandafter\gdef\csname#1cite\endcsname##1{\citejj{##1}}
   \expandafter\gdef\csname#1page\endcsname##1{?}
   \expandafter\gdef\csname#1tag\endcsname##1{\tagjj{##1}}\fi\fi\next}
\newif\ifshowstuffinmargin
\showstuffinmarginfalse
\def\jjtags{\ifx\shlhetal\relax 
  \else
\ifx\shlhetal\undefinedcontrolseq
\else
\showstuffinmargintrue
\ifx\all\relax\else\expandafter\t@gsjj\all,\fi\fi \fi
}

\def\tagjj#1{\realstag{#1}\mginpar{\zeigen{#1}}}
\def\citejj#1{\rechnen{#1}\mginpar{\zeigen{#1}}}     

\def\rechnen#1{\expandafter\ifx\csname stag#1\endcsname\relax ??\else
                           \csname stag#1\endcsname\fi}

\newdimen\theight

\def\marginfont{\sevenrm}

\def\trymarginbox#1{\setbox0=\hbox{\marginfont\hskip\marginshift #1}%
		\global\marginshift\wd0 
		\global\advance\marginshift\margindelta}

\def \mginpar#1{%
\ifvmode\setbox0\hbox to \hsize{\hfill\rlap{\marginfont\quad#1}}%
\ht0 0cm
\dp0 0cm
\box0\vskip-\baselineskip
\else 
             \vadjust{\trymarginbox{#1}%
		\ifdim\marginshift>\marginmax \global\marginshift\marginmin
			\trymarginbox{#1}%
                \fi
             \theight=\ht0
             \advance\theight by \dp0    \advance\theight by \lineskip
             \kern -\theight \vbox to \theight{\rightline{\rlap{\box0}}%
\vss}}\fi}


\def\t@gsoff#1,{\def\@{#1}\ifx\@\empty\let\next=\relax\else\let\next=\t@gsoff
   \def\@@{p}\ifx\@\@@\else
   \expandafter\gdef\csname#1cite\endcsname##1{\zeigen{##1}}
   \expandafter\gdef\csname#1page\endcsname##1{?}
   \expandafter\gdef\csname#1tag\endcsname##1{\zeigen{##1}}\fi\fi\next}
\def\verbatimtags{\showstuffinmarginfalse
\ifx\all\relax\else\expandafter\t@gsoff\all,\fi}
\def\zeigen#1{\hbox{$\langle$}#1\hbox{$\rangle$}}

\def\margincite#1{\ifshowstuffinmargin\mginpar{\zeigen{#1}}\fi}

\def\(#1){\edef\dot@g{\ifmmode\ifinner(\hbox{\noexpand\etag{#1}})
   \else\noexpand\eqno(\hbox{\noexpand\etag{#1}})\fi
   \else(\noexpand\ecite{#1})\fi}\dot@g}

\newif\ifbr@ck
\def\eat#1{}
\def\[#1]{\br@cktrue[\br@cket#1'X]}
\def\br@cket#1'#2X{\def\temp{#2}\ifx\temp\empty\let\next\eat
   \else\let\next\br@cket\fi
   \ifbr@ck\br@ckfalse\br@ck@t#1,X\else\br@cktrue#1\fi\next#2X}
\def\br@ck@t#1,#2X{\def\temp{#2}\ifx\temp\empty\let\neext\eat
   \else\let\neext\br@ck@t\def\temp{,}\fi
   \def\teemp{#1}\ifx\teemp\empty\else\rcite{#1}\fi\temp\neext#2X}
\def\resetbr@cket{\gdef\[##1]{[\rtag{##1}]}}
\def\references{\resetbr@cket\newsection References\par}

\newtoks\symb@ls\newtoks\s@mb@ls\newtoks\p@gelist\n@wcount\ftn@mber
    \ftn@mber=1\newif\ifftn@mbers\ftn@mbersfalse\newif\ifbyp@ge\byp@gefalse
\def\defm@rk{\ifftn@mbers\n@mberm@rk\else\symb@lm@rk\fi}
\def\n@mberm@rk{\xdef\m@rk{{\the\ftn@mber}}%
    \global\advance\ftn@mber by 1 }
\def\rot@te#1{\let\temp=#1\global#1=\expandafter\r@t@te\the\temp,X}
\def\r@t@te#1,#2X{{#2#1}\xdef\m@rk{{#1}}}
\def\b@@st#1{{$^{#1}$}}\def\str@p#1{#1}
\def\symb@lm@rk{\ifbyp@ge\rot@te\p@gelist\ifnum\expandafter\str@p\m@rk=1 
    \s@mb@ls=\symb@ls\fi\write\f@nsout{\number\count0}\fi \rot@te\s@mb@ls}
\def\byp@ge{\byp@getrue\n@wwrite\f@nsin\openin\f@nsin=\jobname.fns 
    \n@wcount\currentp@ge\currentp@ge=0\p@gelist={0}
    \re@dfns\closein\f@nsin\rot@te\p@gelist
    \n@wread\f@nsout\openout\f@nsout=\jobname.fns }
\def\m@kelist#1X#2{{#1,#2}}
\def\re@dfns{\ifeof\f@nsin\let\next=\relax\else\read\f@nsin to \f@nline
    \ifx\f@nline\v@idline\else\let\t@mplist=\p@gelist
    \ifnum\currentp@ge=\f@nline
    \global\p@gelist=\expandafter\m@kelist\the\t@mplistX0
    \else\currentp@ge=\f@nline
    \global\p@gelist=\expandafter\m@kelist\the\t@mplistX1\fi\fi
    \let\next=\re@dfns\fi\next}
\def\symbols#1{\symb@ls={#1}\s@mb@ls=\symb@ls} 
\def\bigsymbol{\textstyle}
\symbols{\bigsymbol\ast,\dagger,\ddagger,\sharp,\flat,\natural,\star}
\def\ftnumbers{\ftn@mberstrue} \def\ftsymbols{\ftn@mbersfalse}
\def\paginal{\byp@ge} \def\resetftnumbers{\ftn@mber=1}
\def\ftnote#1{\defm@rk\expandafter\expandafter\expandafter\footnote
    \expandafter\b@@st\m@rk{#1}}

\long\def\jump#1\endjump{}
\def\ssum{\mathop{\lower .1em\hbox{$\textstyle\Sigma$}}\nolimits}

\def\qed{\nobreak\kern 1em \vrule height .5em width .5em depth 0em}
\def\newneq{\hbox{\rlap{\hbox to 1\wd9{\hss$=$\hss}}\raise .1em 
   \hbox to 1\wd9{\hss$\scriptscriptstyle/$\hss}}}
\def\subsetne{\setbox9 = \hbox{$\subset$}\mathrel{\hbox{\rlap
   {\lower .4em \newneq}\raise .13em \hbox{$\subset$}}}}
\def\supsetne{\setbox9 = \hbox{$\subset$}\mathrel{\hbox{\rlap
   {\lower .4em \newneq}\raise .13em \hbox{$\supset$}}}}

\def\vbar{\mathchoice{\vrule height6.3ptdepth-.5ptwidth.8pt\kern-.8pt}
   {\vrule height6.3ptdepth-.5ptwidth.8pt\kern-.8pt}
   {\vrule height4.1ptdepth-.35ptwidth.6pt\kern-.6pt}
   {\vrule height3.1ptdepth-.25ptwidth.5pt\kern-.5pt}}
\def\f@dge{\mathchoice{}{}{\mkern.5mu}{\mkern.8mu}}
\def\b@c#1#2{{\rm \mkern#2mu\vbar\mkern-#2mu#1}}
\def\b@b#1{{\rm I\mkern-3.5mu #1}}
\def\b@a#1#2{{\rm #1\mkern-#2mu\f@dge #1}}
\def\bb#1{{\count4=`#1 \advance\count4by-64 \ifcase\count4\or\b@a A{11.5}\or
   \b@b B\or\b@c C{5}\or\b@b D\or\b@b E\or\b@b F \or\b@c G{5}\or\b@b H\or
   \b@b I\or\b@c J{3}\or\b@b K\or\b@b L \or\b@b M\or\b@b N\or\b@c O{5} \or
   \b@b P\or\b@c Q{5}\or\b@b R\or\b@a S{8}\or\b@a T{10.5}\or\b@c U{5}\or
   \b@a V{12}\or\b@a W{16.5}\or\b@a X{11}\or\b@a Y{11.7}\or\b@a Z{7.5}\fi}}

\catcode`\X=11 \catcode`\@=12




\let\thischap\jobname

\def\partof#1{\csname returnthe#1part\endcsname}
\def\chapof#1{\csname returnthe#1chap\endcsname}

\def\setchapter#1,#2,#3.{%
  \expandafter\def\csname returnthe#1part\endcsname{#2}%
  \expandafter\def\csname returnthe#1chap\endcsname{#3}%
}

\setchapter 300a,A,I.
\setchapter 300b,A,II.
\setchapter 300c,A,III.
\setchapter 300d,A,IV.
\setchapter 300e,A,V.
\setchapter 300f,A,VI.
\setchapter 300g,A,VII.
\setchapter F604,B,0.
\setchapter  88r,B,I.
\setchapter  600,B,II.
\setchapter  705,B,III.

\def\cprefix#1{
\edef\theotherpart{\partof{#1}}\edef\theotherchap{\chapof{#1}}%
\ifx\theotherpart\thispart
   \ifx\theotherchap\thischap 
    \else 
     \theotherchap%
    \fi
   \else 
     \theotherpart.\theotherchap\fi}

\def\sectioncite[#1]#2{%
     \cprefix{#2}#1}

\edef\thispart{\partof{\thischap}}
\edef\thischap{\chapof{\thischap}}



\expandafter\ifx\csname citeadd.tex\endcsname\relax
\expandafter\gdef\csname citeadd.tex\endcsname{}
\else \message{Hey!  Apparently you were trying to
\string\input{citeadd.tex} twice.   This does not make sense.} 
\errmessage{Please edit your file (probably \jobname.tex) and remove
any duplicate ``\string\input'' lines}\endinput\fi

\def\sciteu{\sciteerror{undefined}}

\def\sciteerror#1#2{{\mathortextbf{\scite{#2}}}\complainaboutcitation{#1}{#2}}
\def\mathortextbf#1{\hbox{\bf #1}}
\def\complainaboutcitation#1#2{%
\vadjust{\line{\llap{---$\!\!>$ }\qquad scite$\{$#2$\}$ #1\hfil}}}

\sectno=-1   
\localtags
\jjtags
\NoBlackBoxes
\define\mr{\medskip\roster}
\define\sn{\smallskip\noindent}
\define\mn{\medskip\noindent}
\define\bn{\bigskip\noindent}
\define\ub{\underbar}
\define\wilog{\text{without loss of generality}}
\define\ermn{\endroster\medskip\noindent}
\define\dbca{\dsize\bigcap}
\define\dbcu{\dsize\bigcup}
\define \nl{\newline}
\magnification=\magstep 1
\documentstyle{amsppt}

{    
\catcode`@11

\ifx\alicetwothousandloaded@\relax
  \endinput\else\global\let\alicetwothousandloaded@\relax\fi

\gdef\subjclass{\let\savedef@\subjclass
 \def\subjclass##1\endsubjclass{\let\subjclass\savedef@
   \toks@{\def\usualspace{{\rm\enspace}}\eightpoint}%
   \toks@@{##1\unskip.}%
   \edef\thesubjclass@{\the\toks@
     \frills@{{\noexpand\rm2000 {\noexpand\it Mathematics Subject
       Classification}.\noexpand\enspace}}%
     \the\toks@@}}%
  \nofrillscheck\subjclass}
} 


\expandafter\ifx\csname alice2jlem.tex\endcsname\relax
  \expandafter\xdef\csname alice2jlem.tex\endcsname{\the\catcode`@}
\else \message{Hey!  Apparently you were trying to
\string\input{alice2jlem.tex}  twice.   This does not make sense.}
\errmessage{Please edit your file (probably \jobname.tex) and remove
any duplicate ``\string\input'' lines}\endinput\fi

\expandafter\ifx\csname bib4plain.tex\endcsname\relax
  \expandafter\gdef\csname bib4plain.tex\endcsname{}
\else \message{Hey!  Apparently you were trying to \string\input
  bib4plain.tex twice.   This does not make sense.}
\errmessage{Please edit your file (probably \jobname.tex) and remove
any duplicate ``\string\input'' lines}\endinput\fi

\def\renewcommand{\newcommand}	       
\edef\cite{\the\catcode`@}%
\catcode`@ = 11
\let\@oldatcatcode = \cite
\chardef\@letter = 11
\chardef\@other = 12
%
%
%
%
\def\@innerdef#1#2{\edef#1{\expandafter\noexpand\csname #2\endcsname}}%
%
%
\@innerdef\@innernewcount{newcount}%
\@innerdef\@innernewdimen{newdimen}%
\@innerdef\@innernewif{newif}%
\@innerdef\@innernewwrite{newwrite}%
%
%
%
\def\@gobble#1{}%
%
%
%
\ifx\inputlineno\@undefined
   \let\@linenumber = \empty 
\else
   \def\@linenumber{\the\inputlineno:\space}%
\fi
%
%
%
\def\@futurenonspacelet#1{\def\cs{#1}%
   \afterassignment\@stepone\let\@nexttoken=
}%
\begingroup 
\def\\{\global\let\@stoken= }%
\\ 
\endgroup
\def\@stepone{\expandafter\futurelet\cs\@steptwo}%
\def\@steptwo{\expandafter\ifx\cs\@stoken\let\@@next=\@stepthree
   \else\let\@@next=\@nexttoken\fi \@@next}%
\def\@stepthree{\afterassignment\@stepone\let\@@next= }%
%
%
%
\def\@getoptionalarg#1{%
   \let\@optionaltemp = #1%
   \let\@optionalnext = \relax
   \@futurenonspacelet\@optionalnext\@bracketcheck
}%
%
%
\def\@bracketcheck{%
   \ifx [\@optionalnext
      \expandafter\@@getoptionalarg
   \else
      \let\@optionalarg = \empty
      \expandafter\@optionaltemp
   \fi
}%
\def\@@getoptionalarg[#1]{%
   \def\@optionalarg{#1}%
   \@optionaltemp
}%
%
%
%
\def\@nnil{\@nil}%
\def\@fornoop#1\@@#2#3{}%
\def\@for#1:=#2\do#3{%
   \edef\@fortmp{#2}%
   \ifx\@fortmp\empty \else
      \expandafter\@forloop#2,\@nil,\@nil\@@#1{#3}%
   \fi
}%
\def\@forloop#1,#2,#3\@@#4#5{\def#4{#1}\ifx #4\@nnil \else
       #5\def#4{#2}\ifx #4\@nnil \else#5\@iforloop #3\@@#4{#5}\fi\fi
}%
\def\@iforloop#1,#2\@@#3#4{\def#3{#1}\ifx #3\@nnil
       \let\@nextwhile=\@fornoop \else
      #4\relax\let\@nextwhile=\@iforloop\fi\@nextwhile#2\@@#3{#4}%
}%
%
%
%
\@innernewif\if@fileexists
\def\@testfileexistence{\@getoptionalarg\@finishtestfileexistence}%
\def\@finishtestfileexistence#1{%
   \begingroup
      \def\extension{#1}%
      \immediate\openin0 =
         \ifx\@optionalarg\empty\jobname\else\@optionalarg\fi
         \ifx\extension\empty \else .#1\fi
         \space
      \ifeof 0
         \global\@fileexistsfalse
      \else
         \global\@fileexiststrue
      \fi
      \immediate\closein0
   \endgroup
}%
%
%
%
%
\def\bibliographystyle#1{%
   \@readauxfile
   \@writeaux{\string\bibstyle{#1}}%
}%
\let\bibstyle = \@gobble
%
%
\let\bblfilebasename = \jobname
\def\bibliography#1{%
   \@readauxfile
   \@writeaux{\string\bibdata{#1}}%
   \@testfileexistence[\bblfilebasename]{bbl}%
   \if@fileexists
      \nobreak
      \@readbblfile
   \fi
}%
\let\bibdata = \@gobble
%
%
\def\nocite#1{%
   \@readauxfile
   \@writeaux{\string\citation{#1}}%
}%
\@innernewif\if@notfirstcitation
%
%
\def\cite{\@getoptionalarg\@cite}%
%
%
\def\@cite#1{%
   \let\@citenotetext = \@optionalarg
   \printcitestart
   \nocite{#1}%
   \@notfirstcitationfalse
   \@for \@citation :=#1\do
   {%
      \expandafter\@onecitation\@citation\@@
   }%
   \ifx\empty\@citenotetext\else
      \printcitenote{\@citenotetext}%
   \fi
   \printcitefinish
}%
\newif\ifweareinprivate
\weareinprivatetrue
\ifx\shlhetal\undefinedcontrolseq\weareinprivatefalse\fi
\ifx\shlhetal\relax\weareinprivatefalse\fi
\def\@onecitation#1\@@{%
   \if@notfirstcitation
      \printbetweencitations
   \fi
   \expandafter \ifx \csname\@citelabel{#1}\endcsname \relax
      \if@citewarning
         \message{\@linenumber Undefined citation `#1'.}%
      \fi
     \ifweareinprivate
      \expandafter\gdef\csname\@citelabel{#1}\endcsname{%
\strut 
\vadjust{\vskip-\dp\strutbox
\vbox to 0pt{\vss\parindent0cm \leftskip=\hsize 
\advance\leftskip3mm
\advance\hsize 4cm\strut\openup-4pt 
\rightskip 0cm plus 1cm minus 0.5cm ?  #1 ?\strut}}
         {\tt
            \escapechar = -1
            \nobreak\hskip0pt\pfeilsw
            \expandafter\string\csname#1\endcsname
             \pfeilso
            \nobreak\hskip0pt
         }%
      }%
     \else  
      \expandafter\gdef\csname\@citelabel{#1}\endcsname{%
            {\tt\expandafter\string\csname#1\endcsname}
      }%
     \fi  
   \fi
   \csname\@citelabel{#1}\endcsname
   \@notfirstcitationtrue
}%
%
%
\def\@citelabel#1{b@#1}%
%
%
\def\@citedef#1#2{\expandafter\gdef\csname\@citelabel{#1}\endcsname{#2}}%
%
%
%
\def\@readbblfile{%
   \ifx\@itemnum\@undefined
      \@innernewcount\@itemnum
   \fi
   \begingroup
      \def\begin##1##2{%
         \setbox0 = \hbox{\biblabelcontents{##2}}%
         \biblabelwidth = \wd0
      }%
      \def\end##1{}
      %
      %
      \@itemnum = 0
      \def\bibitem{\@getoptionalarg\@bibitem}%
      \def\@bibitem{%
         \ifx\@optionalarg\empty
            \expandafter\@numberedbibitem
         \else
            \expandafter\@alphabibitem
         \fi
      }%
      \def\@alphabibitem##1{%
         \expandafter \xdef\csname\@citelabel{##1}\endcsname {\@optionalarg}%
         \ifx\biblabelprecontents\@undefined
            \let\biblabelprecontents = \relax
         \fi
         \ifx\biblabelpostcontents\@undefined
            \let\biblabelpostcontents = \hss
         \fi
         \@finishbibitem{##1}%
      }%
      \def\@numberedbibitem##1{%
         \advance\@itemnum by 1
         \expandafter \xdef\csname\@citelabel{##1}\endcsname{\number\@itemnum}%
         \ifx\biblabelprecontents\@undefined
            \let\biblabelprecontents = \hss
         \fi
         \ifx\biblabelpostcontents\@undefined
            \let\biblabelpostcontents = \relax
         \fi
         \@finishbibitem{##1}%
      }%
      \def\@finishbibitem##1{%
         \biblabelprint{\csname\@citelabel{##1}\endcsname}%
         \@writeaux{\string\@citedef{##1}{\csname\@citelabel{##1}\endcsname}}%
         \ignorespaces
      }%
      %
      %
      \let\em = \bblem
      \let\newblock = \bblnewblock
      \let\sc = \bblsc
      \frenchspacing
      \clubpenalty = 4000 \widowpenalty = 4000
      \tolerance = 10000 \hfuzz = .5pt
      \everypar = {\hangindent = \biblabelwidth
                      \advance\hangindent by \biblabelextraspace}%
      \bblrm
      \parskip = 1.5ex plus .5ex minus .5ex
      \biblabelextraspace = .5em
      \bblhook
      \input \bblfilebasename.bbl
   \endgroup
}%
%
%
\@innernewdimen\biblabelwidth
\@innernewdimen\biblabelextraspace
%
%
%
\def\biblabelprint#1{%
   \noindent
   \hbox to \biblabelwidth{%
      \biblabelprecontents
      \biblabelcontents{#1}%
      \biblabelpostcontents
   }%
   \kern\biblabelextraspace
}%
%
%
%
\def\biblabelcontents#1{{\bblrm [#1]}}%
%
%
\def\bblrm{\rm}%
%
%
\def\bblem{\it}%
%
%
\def\bblsc{\ifx\@scfont\@undefined
              \font\@scfont = cmcsc10
           \fi
           \@scfont
}%
%
%
\def\bblnewblock{\hskip .11em plus .33em minus .07em }%
%
%
\let\bblhook = \empty
%
%
%
\def\printcitestart{[}
\def\printcitefinish{]}
\def\printbetweencitations{, }
\def\printcitenote#1{, #1}
%
%
%
\let\citation = \@gobble
%
%
%
\@innernewcount\@numparams
%
%
\def\newcommand#1{%
   \def\@commandname{#1}%
   \@getoptionalarg\@continuenewcommand
}%
%
%
\def\@continuenewcommand{%
   \@numparams = \ifx\@optionalarg\empty 0\else\@optionalarg \fi \relax
   \@newcommand
}%
%
%
\def\@newcommand#1{%
   \def\@startdef{\expandafter\edef\@commandname}%
   \ifnum\@numparams=0
      \let\@paramdef = \empty
   \else
      \ifnum\@numparams>9
         \errmessage{\the\@numparams\space is too many parameters}%
      \else
         \ifnum\@numparams<0
            \errmessage{\the\@numparams\space is too few parameters}%
         \else
            \edef\@paramdef{%
               \ifcase\@numparams
                  \empty  No arguments.
               \or ####1%
               \or ####1####2%
               \or ####1####2####3%
               \or ####1####2####3####4%
               \or ####1####2####3####4####5%
               \or ####1####2####3####4####5####6%
               \or ####1####2####3####4####5####6####7%
               \or ####1####2####3####4####5####6####7####8%
               \or ####1####2####3####4####5####6####7####8####9%
               \fi
            }%
         \fi
      \fi
   \fi
   \expandafter\@startdef\@paramdef{#1}%
}%
%
%
%
%
\def\@readauxfile{%
   \if@auxfiledone \else 
      \global\@auxfiledonetrue
      \@testfileexistence{aux}%
      \if@fileexists
         \begingroup
            \endlinechar = -1
            \catcode`@ = 11
            \input \jobname.aux
         \endgroup
      \else
         \message{\@undefinedmessage}%
         \global\@citewarningfalse
      \fi
      \immediate\openout\@auxfile = \jobname.aux
   \fi
}%
%
%
\newif\if@auxfiledone
\ifx\noauxfile\@undefined \else \@auxfiledonetrue\fi
%
%
%
%
\@innernewwrite\@auxfile
\def\@writeaux#1{\ifx\noauxfile\@undefined \write\@auxfile{#1}\fi}%
%
%
%
\ifx\@undefinedmessage\@undefined
   \def\@undefinedmessage{No .aux file; I won't give you warnings about
                          undefined citations.}%
\fi
%
%
\@innernewif\if@citewarning
\ifx\noauxfile\@undefined \@citewarningtrue\fi
%
%
%
\catcode`@ = \@oldatcatcode

\def\pfeilso{\leavevmode
            \vrule width 1pt height9pt depth 0pt\relax
           \vrule width 1pt height8.7pt depth 0pt\relax
           \vrule width 1pt height8.3pt depth 0pt\relax
           \vrule width 1pt height8.0pt depth 0pt\relax
           \vrule width 1pt height7.7pt depth 0pt\relax
            \vrule width 1pt height7.3pt depth 0pt\relax
            \vrule width 1pt height7.0pt depth 0pt\relax
            \vrule width 1pt height6.7pt depth 0pt\relax
            \vrule width 1pt height6.3pt depth 0pt\relax
            \vrule width 1pt height6.0pt depth 0pt\relax
            \vrule width 1pt height5.7pt depth 0pt\relax
            \vrule width 1pt height5.3pt depth 0pt\relax
            \vrule width 1pt height5.0pt depth 0pt\relax
            \vrule width 1pt height4.7pt depth 0pt\relax
            \vrule width 1pt height4.3pt depth 0pt\relax
            \vrule width 1pt height4.0pt depth 0pt\relax
            \vrule width 1pt height3.7pt depth 0pt\relax
            \vrule width 1pt height3.3pt depth 0pt\relax
            \vrule width 1pt height3.0pt depth 0pt\relax
            \vrule width 1pt height2.7pt depth 0pt\relax
            \vrule width 1pt height2.3pt depth 0pt\relax
            \vrule width 1pt height2.0pt depth 0pt\relax
            \vrule width 1pt height1.7pt depth 0pt\relax
            \vrule width 1pt height1.3pt depth 0pt\relax
            \vrule width 1pt height1.0pt depth 0pt\relax
            \vrule width 1pt height0.7pt depth 0pt\relax
            \vrule width 1pt height0.3pt depth 0pt\relax}

\def\pfeilsw{ \leavevmode 
            \vrule width 1pt height0.3pt depth 0pt\relax
            \vrule width 1pt height0.7pt depth 0pt\relax
            \vrule width 1pt height1.0pt depth 0pt\relax
            \vrule width 1pt height1.3pt depth 0pt\relax
            \vrule width 1pt height1.7pt depth 0pt\relax
            \vrule width 1pt height2.0pt depth 0pt\relax
            \vrule width 1pt height2.3pt depth 0pt\relax
            \vrule width 1pt height2.7pt depth 0pt\relax
            \vrule width 1pt height3.0pt depth 0pt\relax
            \vrule width 1pt height3.3pt depth 0pt\relax
            \vrule width 1pt height3.7pt depth 0pt\relax
            \vrule width 1pt height4.0pt depth 0pt\relax
            \vrule width 1pt height4.3pt depth 0pt\relax
            \vrule width 1pt height4.7pt depth 0pt\relax
            \vrule width 1pt height5.0pt depth 0pt\relax
            \vrule width 1pt height5.3pt depth 0pt\relax
            \vrule width 1pt height5.7pt depth 0pt\relax
            \vrule width 1pt height6.0pt depth 0pt\relax
            \vrule width 1pt height6.3pt depth 0pt\relax
            \vrule width 1pt height6.7pt depth 0pt\relax
            \vrule width 1pt height7.0pt depth 0pt\relax
            \vrule width 1pt height7.3pt depth 0pt\relax
            \vrule width 1pt height7.7pt depth 0pt\relax
            \vrule width 1pt height8.0pt depth 0pt\relax
            \vrule width 1pt height8.3pt depth 0pt\relax
            \vrule width 1pt height8.7pt depth 0pt\relax
            \vrule width 1pt height9pt depth 0pt\relax
      }


\def\widestnumber#1#2{}

\def\citewarning#1{\ifx\shlhetal\relax 
    \else
    \par{#1}\par
    \fi
}

\def\rm{\fam0 \tenrm}

\def\fakesubhead#1\endsubhead{\bigskip\noindent{\bf#1}\par}



%
%
%

%

\font\textrsfs=rsfs10
\font\scriptrsfs=rsfs7
\font\scriptscriptrsfs=rsfs5

\newfam\rsfsfam
\textfont\rsfsfam=\textrsfs
\scriptfont\rsfsfam=\scriptrsfs
\scriptscriptfont\rsfsfam=\scriptscriptrsfs

\edef\oldcatcodeofat{\the\catcode`\@}
\catcode`\@11

\def\Cal@@#1{\noaccents@ \fam \rsfsfam #1}

\catcode`\@\oldcatcodeofat


\expandafter\ifx \csname margininit\endcsname \relax\else\margininit\fi

\long\def\red#1\endred{}
\long\def\green#1\endgreen{}
\long\def\blue#1\endblue{}

\def\endred{ \unmatched endred! }
\def\endgreen{ \unmatched endgreen! }
\def\endblue{ \unmatched endblue! }

\ifx\latexcolors\undefinedcs\def\latexcolors{}\fi

\def\emptycs{}
\def\evaluatelatexcolors{%
        \ifx\latexcolors\emptycs\else
        \expandafter\xxevaluate\latexcolors\xxfertig\evaluatelatexcolors\fi}
\def\xxevaluate#1,#2\xxfertig{\setupthiscolor{#1}%
        \def\latexcolors{#2}}

\font\smallfont=cmsl7
\def\rutgerscolor{\ifmmode\else\endgraf\fi\smallfont
\advance\leftskip0.5cm\relax}
\def\setupthiscolor#1{\edef\tmptmpcs{\noexpand\bgroup\noexpand\rutgerscolor
\noexpand\def\noexpand\currentcolor{#1}%
\noexpand}%
\expandafter\let\csname#1\endcsname\tmptmpcs
\def\tmptmpcs{\checkColorUnmatched{#1}\popthecolor}
\expandafter\let\csname end#1\endcsname\tmptmpcs}

\def\checkColorUnmatched#1{\def\expectcolor{#1}%
    \ifx\expectcolor\currentcolor   
    \else \edef\failhere{\noexpand\tryingToClose '\currentcolor' with end\expectcolor}\failhere\fi}

\def\currentcolor{???}

\def\popthecolor{\ifmmode\else\endgraf\fi\egroup}

\expandafter\def\csname#1\endcsname{}

\evaluatelatexcolors

 \let\outerhead\head
 \def\head{\innerhead}
 \let\innerhead\outerhead

 \let\outersubhead\subhead
 \def\subhead{\innersubhead}
 \let\innersubhead\outersubhead

 \let\outersubsubhead\subsubhead
 \def\subsubhead{\innersubsubhead}
 \let\innersubsubhead\outersubsubhead

 \def\proclaim{\innerproclaim}
 \let\innerproclaim\outerproclaim

 %
 %
 %
 %

\def\demo#1{\medskip\noindent{\it #1.\/}}
\def\enddemo{\smallskip}

\def\remark#1{\medskip\noindent{\it #1.\/}}
\def\endremark{\smallskip}

\topmatter
\title{MORE ON CARDINAL ARITHMETIC \\
 Sh410} \endtitle
\author {Saharon Shelah \thanks {Done 1990 \null\newline
Partially supported by the Basic Research Fund, Israel Academy of Science
\null\newline
I thank Alice Leonhardt for the beautiful typing \null\newline
Publ. 410 \null\newline   
 Latest Revision -  02/March/12} \endthanks} \endauthor
\affil{Institute of Mathematics\\
 The Hebrew University\\
 Jerusalem, Israel
 \medskip
 Rutgers University\\
 Mathematics Department\\
 New Brunswick, NJ  USA} \endaffil
\endtopmatter
\document  
 
\newpage

\head {An Annotated Table of Contents} \endhead  \resetall 
\bn
\S1 $\quad$ On Normal Filters 
\sn
In \cite[\S4]{Sh:355} and \cite[\S3,\S5]{Sh:400} we have computed  
cov$(\lambda,\lambda,\theta,\sigma)$ when $\theta > \sigma = \text{cf }
\sigma > \aleph_0$, using tcf $\dsize \prod_{i < \kappa} \lambda_i/J$, for 
$\sigma$-complete ideals $J$ and $\sigma \leq \kappa < \theta$.  In 
\cite[\S4]{Sh:371} we deal with a similar theorem where we restrict ourselves 
to normal ideals, namely  prc,  but its computation, using  
pp's,  did not always yield exact values (i.e. the upper and lower bounds 
tend not to match).  Here we give reasonably exact values for  
prc$_J(f,\bar \mu)$,  using the true cofinalities of $\dsize \prod_{i <
\kappa} \mu'_i/J_1$, where $J_1$ is a normal filter on $\kappa$ extending
$J$ and for $i < \kappa,\mu'_i$ is a regular ordinal satisfying  
$\mu_i \leq \mu'_i \leq f(i)$.  We also give a sufficient condition for the
existence of normal ideal $J$ on $\kappa$ such that for some sequence  
$\langle \lambda_i:i < \kappa \rangle$ of regulars,  we have $\lambda = 
\text{ tlim}\langle \lambda_i:i < \kappa \rangle,\mu = \text{ tcf}
\dsize \prod_{i < \kappa} \lambda_i/J$.
\bn
\S2 $\quad$  On measures of the size of $[\lambda]^{< \kappa}$
\sn
We mainly investigate cardinals like 

$$
\align
\text{Min} \bigl\{ |{\Cal P}|:&{\Cal P} \subseteq [\lambda]^{< \theta}
(\lambda) \text{ and for every }  Z \in {\Cal S}_{\leq \kappa}(\lambda)
\text{ there is a sequence} \\
  &\langle Z_n:n < \omega \rangle \text{ of subsets of } Z \text{ such that }
Z = \dbcu_{n < \omega} Z_n \text{ and} \\
  &(\forall n < \omega)(\forall y \in {\Cal S}_{< \theta}(Z_n))
  (\exists z \in {\Cal P})[y \subseteq z] \bigr\}.
\endalign
$$ 
\mn
We also give sufficient 
conditions for the strong covering to hold for a pair $(W,V)$ of universes.
\bn
\S3 $\quad$ pcf - Inaccessibility and characterizing the existence of
non $<_J$-decreasing sequences (for topology)
\sn
We restate various results using  pcf  inaccessibility and present more 
consequences of the proofs in \cite[\S2,\S4]{Sh:400}.  We characterize those  
$\kappa < \sigma < \theta$ for which there is a sequence $\langle f_\alpha:
\alpha < \theta \rangle$ of members of ${}^\kappa \sigma$ such that 
$\alpha < \beta \Rightarrow f_\alpha \nleq f_\beta$;  
answering a question of Gerlits, Hajnal and Szentmiklossy.  [See more
in \cite[\S6]{Sh:513}].
\bn
\S4 $\quad$  Entangled Orders - Narrow Order Boolean Algebras Revisited 
\sn
We show that for a class of cardinals $\lambda$ there is an entangled linear
order of cardinality $\lambda^+$.  This holds for $\lambda$ if there is a 
$\kappa$ such that $\kappa^{+4} \leq \text{ cf}(\lambda) < \lambda \leq  
2^\kappa$.  [See more in \cite{Sh:462} and \cite{Sh:666}.]
\bn
\S5 $\quad$ prd: Measuring $\dsize \prod_{i < \kappa} f(i)$ by a family of
ideals and a family of sequences $\langle B_i:i < \kappa \rangle,|B_i| <
\mu_i$
\sn
This generalizes Section 1, replacing normality by an abstract property; we 
also present a generalization of the concept of a normal filter, and deduce  
prd$_J(\bar f,\bar \mu) \leq \text{prd}_J(\bar f,\bar \mu)^+$ and  
prd$_J(\aleph_{\bar f},\bar \mu) \leq \aleph_{(\text{prd}_J(\bar f,
\bar \mu))^+}$ under suitable conditions.
\bn
\S6 $\quad$  The Existence of Strongly Almost Disjoint Families 
\sn
We characterize such existence questions by pp's.  An example is the 
question of the existence of a family of $\lambda^+$ subsets of $\lambda > 
\kappa^{\aleph_0}$, each of cardinality $\kappa \,(> \aleph_0)$ such that
the intersection of any two is finite.
\newpage

\head {\S1  On Normal Filters} \endhead  \resetall \sectno=1
\bn
The following Lemma \scite{1.1} is similar to \cite[5.4]{Sh:355}, 
\cite[3.5]{Sh:400}, but deal with normal ideals (see \cite[\S4]{Sh:371}, in 
particular \cite{Sh:374}, Definition \scite{4.1}, Claim \sciteu{4.6}).  
Remember  prc  is defined in \cite[\S4]{Sh:371} as:
\definition{\stag{1.0} Definition}  1) For a regular uncountable cardinal  
$\kappa$,  normal ideal $J$ on $\kappa,\bar \mu$ a $\kappa$-sequence of 
cardinals $> \kappa$, and $f \in {}^\kappa$ Ord,  we define: 

$$
\align
\text{prc}_J(f,\bar \mu) = \text{ Min} \biggl\{|{\Cal P}|:&{\Cal P}
\text{ is a family of } \kappa\text{-sequences of sets of ordinals}, \\
  &\bar B = \langle B_i:i < \kappa \rangle,|B_i| < \mu_i 
  \text{ or at least} \\
  &\{i < \kappa:|B_i| \geqq \mu_i\} \in J, \text{ such that: for every }
g \in {}^\kappa\text{Ord}, \\
  &g \leqq_J f \text{ there is a sequence } \langle \bar A^\zeta:\zeta < 
\kappa \rangle \text{ of members} \\
  &{\Cal P} \text{ satisfying } \{i < \kappa:g(i) \notin
\dbcu_{\zeta < i} A^\zeta_i\} \in J \biggr\}.
\endalign
$$
\mn
2)  We may write $f$ as a sequence of ordinals say $\langle \lambda_i:i < 
\kappa \rangle$, and if $\lambda_i = \lambda$ for each $i$, we may write  
$\lambda$. \nl
3)  prc$'_J(f,\mu)$ is defined similarly but from
$\bar B = \langle B_i:i < \kappa \rangle$ we demand this time $|B_i| < \mu$.
\enddefinition
\bigskip

\remark{Remark}  See there (\cite[4.2,4.3]{Sh:371}) for some basic properties.
But \scite{1.1} below substantially improves \cite[Claim 4.6]{Sh:371} there.
\endremark
\bigskip

\proclaim{\stag{1.1} Lemma}  1) Let $\kappa$ be a regular uncountable 
cardinal, $f:\kappa \rightarrow$ ordinals,  $J$ a normal ideal on $\kappa$,  
and $\bar \mu  = \langle \mu_i:i < \kappa \rangle$ a sequence of regular 
cardinals.  \ub{Then}

$$
\align
{\text{\rm prc\/}}_J(f,\bar \mu) = {\text{\rm sup\/}} \biggl\{
{\text{\rm tcf\/}}
[\dsize \prod_{i < \kappa} \mu'_i/J_1]:&J_1 \text{ a normal 
ideal on } \kappa \text{ extending }  J, \\
  &\text{such that the tcf is well defined and} \\
  &\{i < \kappa:\text{ not } ``\mu_i \leqq \mu'_i = {\text{\rm cf\/}}
(\mu'_i) \leqq f(i)"\} \in J \biggr\}
\endalign
$$
\mn
provided that: 
\mr
\item "{$(\alpha)$}"   $\mu_i = \mu > \kappa$.
\ermn
2)  We can replace assumption $(\alpha)$ by $(\beta)$ below, and  
$\dsize \bigwedge_i \mu_i > \kappa$. 
\mr
\item "{$(\beta)$}"   $\mu_i$ strictly increasing and for limit $i$: \nl
if $i = \dsize \sum_{j<i} \mu_j$ is regular then $i = \mu_i$, otherwise  
$\mu_i = (\dsize \sum_{j<i} \mu_j)^+$.
\endroster
\endproclaim
\bigskip

\remark{\stag{1.1A} Remark}  1) On getting $=^+$ see \cite[6.1(C)]{Sh:420} and 
\cite[\S4]{Sh:430}.  The problem is when pcf$({\frak a})$ 
has an accumulation point which is inaccessible. \nl
2)  In the case $(\beta)$ holds, if $\bar \mu^* = \langle \mu^*_i:i < \kappa 
\rangle$ is (strictly) increasing continuous, 
$\underset {i < \kappa} {}\to \sup \mu^*_i = 
\underset {i < \kappa} {}\to \sup
\mu_i$ then prc$(f,\bar \mu) + (\dsize \sum_{i<\kappa} \mu_i)^+ = 
\text{ prc}(f,\bar \mu^*) + (\dsize \sum_{i<\kappa} \mu_i)^+$, by 
\cite[4.10(2)]{Sh:371} + \cite[2.1]{Sh:355}. \nl
3)  If in $(\beta)$ we place $``\mu_i > \kappa"$ by 
$\underset {i < \kappa} {}\to \sup \mu_i = \kappa,
\dsize \bigwedge_{i<\kappa^*} \mu_i < \kappa$ (so $\kappa$ is inaccessible 
we can get:

$$
\text{prc}_J(f,\bar \mu) = \sup \bigl\{ \text{nor-cf}_i
\dsize \prod_{i < \kappa} \mu'_i/J:\{i:\mu_i \leq \mu'_i = \text{ cf}(\mu'_i)
\leq f(i)\} \in J \bigr\}.
$$
\endremark
\bigskip

\demo{Proof}  \ub{The inequality} $\geq$: 

Same proof as that of $``\lambda(1) \leq \lambda(2)"$ in the proof of 
\cite[4.6]{Sh:371}.
\mn
\ub{The inequality} $\leq$: 

Let $\lambda^*$ be the successor of the sup.
\enddemo
\bigskip
\demo{\stag{1.1B} Fact}  There is a family ${\Cal P}^*$ such that: 
\mr
\item "{$(i)$}"   members of ${\Cal P}^*$ are of the form 
$$
\langle B_{i,\zeta }:i < \kappa,\zeta < \zeta_i \rangle \text{ or }
\langle <B_{i,\zeta }:\zeta < \zeta_i>:i < \kappa \rangle
$$
\ermn
where $\zeta_i < \mu_i$ and each $B_{i,\zeta}$ is a non-empty 
subset of $f(i) + 1$
\mr
\item "{$(ii)$}"   $|{\Cal P}^*| < \lambda^*$ 
\sn
\item "{$(iii)$}"  if $\langle B_{i,\zeta}:i < \kappa,\zeta < \zeta_i
\rangle \in {\Cal P}^*,g \in \dsize \prod_{i< \kappa} (f(i)+1),A \subseteq  
\kappa$ and for $i \in A$ we have 
$\xi_i < \zeta_i$ satisfying $g(i) \in B_{i,\xi_i}$ \ub{then}
there are $E,\langle A_j:j < \kappa \rangle$ and for $j < \kappa$
sequences $\bar B^j = \langle B^j_{i,\zeta}:
i < \kappa,\zeta < \zeta^j_i \rangle$ and $\langle \xi^j_i:
i \in A_j \rangle$ such that:  
{\roster
\itemitem{ $(a)$ }   $E \subseteq \kappa$ and $\kappa \backslash E \in J$;
\sn
\itemitem{ $(b)$ }  $A \cap E = \{i < \kappa:i \in \dbcu_{j<i} A_j\}$;
\sn
\itemitem{ $(c)$ }  $\bar B^j \in {\Cal P}^*$ for $j < \kappa$;
\sn
\itemitem{ $(d)$ }  for $j < \kappa$ and $i \in A_j$ we have:  
$\xi^j_i < \zeta^j_i$ and  
$B^j_{i,(\xi^j_i)} \subseteq B_{i,\xi_i}$ and 
$$
[|B^i_{i,(\xi^j_i)}| \geq \mu_i \Rightarrow |B^j_{i,(\xi^j_i)}| < 
|B_{i,(\xi^j_i)}|]
$$  
\sn
\itemitem{ $(e)$ }  for $i \in A_j$ we have $g(i) \in B^j_{i,(\xi^j_i)}$ 
\endroster}
\item "{$(iv)$}"  $\langle <f(i)+1>:i < \kappa \rangle \in {\Cal P}^*$ 
(so $\zeta_i = 1$ here).
\endroster
\enddemo
\bigskip

\demo{Proof of the Inequality from the Fact}  
Let us define a family ${\Cal P}'$:  
$$
{\Cal P}' = \bigl\{ \langle \cup \{B_{i,\zeta}:\zeta < \zeta_i \text{ and }
|B_{i,\zeta}| < \mu_i\}:i < \kappa \rangle:\langle B_{i,\zeta}:i < \kappa,
\zeta < \zeta_i \rangle \in {\Cal P}^* \bigr\}.
$$
\mn
Now each member of ${\Cal P}'$ has the right form as each $\mu_i$ is 
regular and $\langle B_{i,\zeta}:i < \kappa,\zeta < \zeta_i \rangle \in  
{\Cal P}^*$ implies $\zeta_i < \mu_i$.  Also the cardinality of  
${\Cal P}'$ is $< \lambda^*$ (by (ii) of \scite{1.1A}).

Let $g \leq f$ and it is enough to find $\langle B^\epsilon_i:
i < \kappa \rangle \in {\Cal P}'$ for $\epsilon < \kappa$ such that 
$\{i:g(i) \notin \dbcu_{\epsilon <i} B^\varepsilon_i\} \in J$. 

We choose by induction on $n$, for every $\eta \in {}^n \kappa$ the 
following:  $\bar B^\eta = \langle B^\eta_{i,\zeta}:i < \kappa,\zeta < 
\zeta^\eta_i \rangle \in {\Cal P}^*,\langle \xi^\eta_i:i \in A_\eta \rangle$  
with $\xi^\eta_i < \zeta^\eta_i$ and $A_\eta \subseteq \kappa$ such that:
\mr  
\item "{$(\alpha)$}"   $\bar B^{<>} = \langle <f(i)+1>:i < \kappa \rangle,  
\zeta^{<>}_i = 1,\xi^{<>}_i = 0,A_{<>} = \kappa$,
\sn
\item "{$(\beta)$}"   $[i \in A_{\eta \char94 <j>} \Rightarrow g(i) \in  
B^{\eta \char94 <j>}_{i,(\xi^{\eta \char94 <j>}_i)} \subseteq 
B^\eta_{i,\xi^\eta_i}]$,
\sn
\item "{$(\gamma)$}"  $[i \in A_{\eta \char94 <j>} \& |B^\eta_{i,\xi^\eta_i}| 
\geq \mu_i \Rightarrow |B^{\eta \char94 <j>}_{i,(\xi^{\eta \char94 <j>}_i)}| 
< |B^\eta_{i,\xi^\eta_i}|]$ and 
\sn
\item "{$(\delta)$}"   $\{i \in A_\eta:i \notin \dbcu_{j <i} 
A_{\eta \char94 <j>}\} \in J$.
\ermn
The induction step is by (iii) of the fact; in the end let for $\eta \in  
{}^{\omega >}\kappa$ and $i < \kappa:B^*_{\eta,i} =: \cup \{B^\eta_{i,\zeta}:
\zeta < \zeta^\eta_i \text{ and } |B^\eta_{i,\zeta}| < \mu _i\}$.

Clearly for each $\eta \in {}^{\omega >}\kappa$  we have  
$\langle B^*_{\eta,i}:i < \kappa \rangle \in {\Cal P}'$;  let us 
enumerate ${}^{\omega >}\kappa$ as $\{\rho_\epsilon:\epsilon < \kappa\}$ 
such that $[\rho \triangleleft \rho_\epsilon \Rightarrow \rho \in 
\{\rho_\zeta:\zeta < \epsilon\}]$, and let us define $\bar B^\epsilon = 
\langle B^*_{\rho_\epsilon,i}:i < \kappa \rangle \in {\Cal P}'$,  hence
by Definition \scite{1.0} it is enough to show that $E =: \{i < \kappa:g(i) 
\notin \dbcu_{\epsilon <i} B^*_{\rho_\epsilon,i}\}$ belongs to $J$.  
We know that for every $\eta \in {}^{\omega >}\kappa$ the set $X_\eta =: 
\{i \in A_\eta:i \notin \dbcu_{j<i} A_{\eta \char94 <j>}\}$ belongs to $J$.  
Also the sets 

$$
Y =: \{\delta < \kappa:\delta \text{ limit and } \neg(\forall \epsilon)
[\epsilon < \delta \Rightarrow \rho_\epsilon \in {}^{\omega >}\delta]\}
$$

$$
X =: \{i < \kappa:\text{if } i \notin Y \text{ then for some } \eta \in  
{}^{\omega >}i \text{ we have } i \in X_\eta\}
$$
\mn
belong to $J$.   It suffices to show 
\mr
\item "{$(*)$}"  for every $\delta \in \kappa \backslash X$,  
for some $\eta \in {}^{\omega >}\delta$ we have $\delta \in B^*_{\eta,i}$.
\ermn
Why $(*)$ holds?  Choose by induction on $n < \omega,\rho_n \in {}^n\delta$  
such that:  $\delta \in A_{\rho_n}$.  For $n = 0$ remember $A_{<>} = \kappa$.
For $n + 1$, as $\delta \in A_{\rho_n}$ and $\delta \notin X$ clearly  
$\delta \notin X_{\rho_n}$ so necessarily $\delta \in \dbcu_{j<i} 
A_{\rho_{n \char94 <j>}}$.
\mn
Now $\langle|B^{\rho_n}_{\delta,\xi^{\rho_n}_\delta}|:n < \omega \rangle$ 
is non-increasing hence (by $(\beta)$ above) is eventually constant hence 
(by $(\gamma)$ above) for some
$n,|B^{\rho_n}_{\delta,\xi^{\rho_n}_\delta}| < 
\mu_\delta$,  hence $g(\delta) \in B^*_{\rho_n,\delta}$.  
So $(*)$ holds and we have finished proving that ${\Cal P}'$ exemplify
the inequality from the fact.
\enddemo
\bigskip

\demo{Proof of the Fact 1.1A}  It suffices to prove that for any $\bar B^* = 
\langle B^*_{i,\zeta}:i < \kappa,\zeta < \zeta_i \rangle$ satisfying 
the requirements in (i), we can find ${\Cal P}_{{\bar B}^*}$ satisfying $(i) 
+ (ii)$ and $(iii)$ for the given $\bar B^*$ and any $g,\langle \xi_i:
i < \kappa \rangle$ as there.  Let  $Y^i_0 = \{\zeta < \zeta_i:|B^*_{i,\zeta}|
< \mu_i\},Y^i_1 = \{\zeta < \zeta_i:|B^*_{\zeta,i}| \geq \mu_i > \text{ cf}
|B^*_{\zeta,i}|\}$ and $Y^i_2 = \{\zeta < \zeta_i:\text{ cf}|B^*_{i,\zeta}| 
\geq \mu_i\}$ (for each $i < \kappa,\zeta < \zeta_i)$.

Clearly $\langle Y^i_0,Y^i_1,Y^i_2 \rangle$ is a partition of $\{\zeta:\zeta <
\zeta_i\}$, now for $\zeta \in Y^i_1 \cup Y^i_2$, let $\mu^i_\zeta =: 
\text{ cf}|B^*_{\zeta,i}|$, and $\langle B_{i,\zeta,\epsilon}:\epsilon < 
\mu^i_\zeta \rangle$  be an increasing continuous sequence of subsets of  
$B^*_{i,\zeta}$ of cardinality $< |B^*_{\zeta,i}|$ and $\dbcu_{\epsilon  <
\mu^i_\zeta} B_{i,\zeta,\epsilon} = B^*_{i,\zeta}$.

Now let ${\frak a} =: 
\{\mu^i_\zeta:i < \kappa,\zeta \in Y^i_2\}$, so ${\frak a}$ is a set 
of regular cardinals.
\bn
\ub{Case 1}:    Assume assumption $(\alpha)$ of \scite{1.1}; \nl
so ${\frak a}$ is a set of $< \mu + \kappa^+ = \mu$ regular cardinals, 
each $\geq \mu$, so the pcf analysis of \cite[\S2]{Sh:371} apply.  Let us get 
$\left< \langle f^\theta_\alpha:\alpha < \theta \rangle:\theta 
\in \text{ pcf} ({\frak a}) \right>$.

Now for each $\theta \in \text{ pcf}({\frak a}) 
\cup \{1\}$ which is $< \lambda^*$ 
and $\alpha < \theta$,  we choose  $\bar B^{\theta,\alpha} = 
\langle B^{\theta,\alpha}_{i,\zeta}:
i < \kappa,\zeta < \zeta^{\theta,\alpha}_i \rangle$ such that:  

$$
\align
\{B^{\theta,\alpha}_{i,\zeta}:\zeta < \zeta^{\theta,\alpha}_i\} = &\{B^*_{i,
\zeta}:\zeta \in Y^i_0\} \cup \{B^*_{i,\zeta,\epsilon}:\epsilon < 
\mu^i_\zeta,\zeta \in Y^i_1\} \\
  &\cup \{B_{i,\zeta,f^\theta_\alpha(\mu^i_\zeta)}:\zeta \in Y^i_2 
\text{ and } \theta \in \text{ pcf}({\frak a}) 
(\text{i.e. } \theta \ne 1)\}.
\endalign
$$
\mn
Let ${\Cal P}_{{\bar B}^*} =: \{\bar B^{\theta,\alpha}:\theta \in \lambda^*  
\cap \text{ pcf}({\frak a}) \text{ and } \alpha < \theta$ or $\theta =
1,\alpha = 0\}$.  Now it is as required, in particular 
$|{\Cal P}_{{\bar B}^*}| < \lambda^*$ because $\lambda^* > \sup(\lambda^* \cap
\text{ pcf}({\frak a})$ holds as $\lambda^*$ is a successor cardinal.
\bn
\ub{Case 2}:  Assume assumption.  Here we partition ${\frak a}$ to $\kappa$  
sets diagonally; i.e. 
without loss of generality each $\zeta_i$ is a cardinal hence by
clause $(\beta)$ we have
$\zeta_i \leq \dbcu_{j <i} \mu_{j+1} + \mu_0$.  So for every limit  
$i < \kappa$ we have 
$\zeta_i \leq \dbcu_{j<i} \mu_j$.  Remember $\kappa < \dbcu_{i <
\kappa} \mu_i$ and even $\kappa < \mu_0$ and let for $j < \kappa$: \nl
${\frak a}_j =: \{\mu^i_\zeta:i \in (j,\kappa)$ and 
$\zeta \in Y^i_2$ and $\zeta < \mu_j\}$.  
So $|{\frak a}_j| \leq \kappa + \mu_j < \text{ Min}({\frak a}_j)$
because and relevant $\mu^i_\zeta$ is $\ge \mu_i > \mu_j \ge \kappa$
and if $\langle \xi_i:
i < \kappa \rangle \in \dsize \prod_{i < \kappa} (1 + \zeta_i)$ then we can 
define $h:\kappa \rightarrow \kappa,h(i) < 1 + i$ such that:  
$[i < \kappa \and i \text{ limit } \Rightarrow \mu^i_{\xi_i} \in
{\frak a}_{h(i)}]$.  \hfil$\square_{\scite{1.1}}$\margincite{1.1}
\enddemo
\bn
The following lemma generalizes \cite[1.5]{Sh:371}.
\proclaim{\stag{1.2} Lemma}  Suppose $\sigma_1 \leq \sigma_2 \leq \kappa < 
\theta < \lambda$ are cardinals, $\sigma_1,\sigma_2,\kappa$ are regular, 
$\lambda > \text{ cf}(\lambda) = \kappa > \aleph_0,\lambda < \mu = \text{ cf}
(\mu) < \text{ pp}^+_{\Gamma(\theta,\sigma_2)}(\lambda)$, and for every 
large enough $\lambda' < \lambda,[\sigma_1 \leq \text{ cf}(\lambda') < \theta
\Rightarrow \text{ pp}_{\Gamma(\theta,\sigma_1)}(\lambda') < \lambda]$.

\ub{Then} there is 
an increasing sequence $\langle \mu_i:i < \kappa \rangle$ of
regular cardinals $< \lambda,\lambda = \underset {i < \kappa} {}\to \sup 
\mu_i$ and an ideal $J$ on $\kappa$ satisfying
$\lambda = \text{ tlim}_J \mu_i$ and $\mu = 
\text{ tcf}(\dsize \prod_{i < \kappa} \mu_i/J)$ such that:
\mr
\item "{$(a)$}"   $J$ is $\sigma_2$-complete and extend  
$J^{\text{bd}}_\kappa$
\sn
\item "{$(b)$}"  if $\kappa \geq \sigma^+_1$ and $(\forall \alpha < \kappa)
[\text{cov}(|\alpha|,\sigma_2,\sigma_1,2) < \kappa]$ then $J$ is normal;
\sn
\item "{$(c)$}"  if $\sigma_1 = \aleph_0$ then $J = J^{\text{bd}}_\kappa$.
\endroster
\endproclaim
\bigskip

\demo{Proof}  For $(c)$ see \cite[1.6]{Sh:371}, so we can assume  
$\sigma_1 > \aleph_0$ as otherwise we have there gotten a conclusion 
stronger than (a) + (b) + (c). 

Let ${\frak a} \subseteq \text{ Reg } \cap \lambda$ have 
cardinality $< \theta$, be 
unbounded in $\lambda,I$ a $\sigma_2$-complete ideal on $a$, and $(\forall 
\lambda' < \lambda)[{\frak a} \cap \lambda' \in I]$ and 
$\mu = \text{ tcf}(\Pi {\frak a}/I)$. 
As cf$(\lambda) = \kappa < \theta < \lambda$ without loss of generality  
$\theta < \text{ Min}({\frak a})$, and 
let $\langle \lambda_i:i < \kappa \rangle$  
be increasing continuous with limit $\lambda,\aleph_0 \leq \text{ cf}
(\lambda_i) < \kappa$ (remember $\kappa > \aleph_0$); without loss of 
generality $\theta < \lambda_0 < \text{ Min}){\frak a})$ and 
pp$_{\Gamma(\theta,
\sigma_1)}(\lambda_i) < \lambda_{i+1}$ and, if $i > 0$ is a limit ordinal 
then ${\frak a} \cap \lambda_i$ is unbounded in $\lambda_i$.  
Also without loss of generality for every $i < \kappa$: 
\mr
\item "{$(*)_1$}"  $\lambda_0 \leq \lambda' < 
\lambda_i \and \sigma_1 \leq \text{ cf}(\lambda') 
< \theta \Rightarrow \text{ pp}_{\Gamma(\theta,\sigma_1)}(\lambda') < 
\lambda_{i+1}$
\ermn
hence
\mr
\item "{$(*)_2$}"  ${\frak b} \subseteq (\lambda_0,\lambda_i) 
\and |{\frak b}| < \theta \Rightarrow 
\text{ sup pcf}_{\sigma_1 \text{-complete}}({\frak b}) < \lambda_{i+1}$.
\ermn
Let $\langle {\frak b}_\sigma [{\frak a}]:
\sigma \in \text{ pcf}({\frak a}) \rangle$ be a generating 
sequence (exists by \cite[2.6]{Sh:371}) and without loss of generality 
$\mu = \text{ max pcf}({\frak a})$.  
By \cite[3.6]{Sh:345a} (and \scite{3.1}(7)) 
without loss of generality $\sigma \in {\frak b}_\chi[{\frak a}] 
\Rightarrow {\frak b}_\sigma[{\frak a}] \subseteq {\frak
b}_\chi[{\frak a}]$.   Let $i < \kappa$ satisfy cf$(i) \geq \sigma_1$, as  
$|{\frak a}| < \theta$ and $(*)_1$ we have
$\text{sup pcf}_{\sigma_1 \text{-complete}}(\lambda_i \cap {\frak a})
< \lambda_{i+1}$ hence for some ${\frak c}_i$ we have:  
${\frak c}_i \subseteq \lambda_{i+1} \cap \text{ pcf}({\frak a} 
\cap \lambda_i),|{\frak c}_i| < \sigma_1$ and ${\frak a} \cap
\lambda_i \subseteq \dbcu_{\theta \in {\frak c}_i} {\frak
b}_\theta[{\frak a}]$ (otherwise we can find a 
$\sigma_1$-complete proper ideal $J$ on $\lambda_i \cap {\frak a}$ such that

$$
[\sigma < {\text{\rm pp\/}}^+_{\Gamma(\theta,\sigma_1)}(\lambda_i) 
\and \sigma \in  
{\text{\rm pcf\/}}(\lambda_i \cap {\frak a}) \Rightarrow 
{\frak b}_\sigma[{\frak a}] \cap (\lambda_i \cap {\frak a} \backslash 
\mu_i) \in J],
$$
\mn
a contradiction to $(*)_1$.
Note that ${\frak c}_i \subseteq \lambda_{i+1}$.

Let $S_0 = \{\delta < \kappa:{\text{\rm cf\/}}(\delta)\geqq \sigma_1\}$, 
so for some $i(*) < \kappa$ we have: 

$$
S_1 = \{ \delta \in S_0:{\frak c}_\delta \cap \lambda_\delta \subseteq 
\lambda_{i(*)}\}.
$$
\mn
is a stationary subset of $\kappa$.
\sn
By renaming, without loss of generality $i(*) = 0$ and so ${\frak c}_i
\subseteq (\lambda_i,\lambda_{i+1})$, and for $i \in S_1$ let  
$\langle(\theta_{i,\zeta},{\frak e}_{i,\zeta}):
\zeta < \zeta_i \rangle$ list $\{(\theta,{\frak b}_\theta[{\frak a}]):
\theta \in {\frak c}_i\}$; so: 
\mr
\item "{$(*)_1$}"   ${\frak a} \cap \lambda_i \subseteq 
\dbcu_{\zeta < \zeta_i}
{\frak e}_{i,\zeta}$,  max pcf$({\frak e}_{i,\zeta}) = 
\theta_{i,\zeta},\lambda_i < \theta_{i,
\zeta} < \lambda_{i+1}$ \nl
and $\theta_{i,\zeta} \in \text{ pcf}_{\sigma \text{-complete}} 
({\frak a} \cap \lambda_i)$.
\endroster
\enddemo
\bigskip

\demo{\stag{1.2A} Fact}  There are finite ${\frak d}_{i,\zeta} 
\subseteq \text{ pcf}({\frak e}_{i,\zeta}) \backslash 
\lambda_i$ for $i \in S_1,\zeta < \zeta_i$ and 
stationary $S_2 \subseteq S_1$ such that letting ${\frak d}_i = 
\dbcu_{\zeta < \zeta_i} d_{i,\zeta}$ we have: if $S \subseteq S_2,\kappa = 
\text{ sup}(S)$ then $\mu \in \text{ pcf}_{\sigma_2 \text{-complete}}
\left( \dbcu_{i \in S} {\frak d}_i \right)$.
\enddemo
\bigskip

\demo{Proof of 1.2 from the Fact}  Now the preliminary part of \scite{1.2} 
is easy; as ${\frak d}_i \subseteq (\lambda_i,\lambda_{i+1})$, and 
$|{\frak d}_i| < \sigma_1 < \kappa = \text{ cf}(\kappa)$, clearly 
${\frak d} =: \dbcu_{i \in S_1} {\frak d}_i$ has order 
type $\kappa$, and by \scite{1.2A} $\mu \in 
\text{ pcf}_{\sigma_2 \text{-complete}} \left(\dbcu_{i \in S_1} 
{\frak d}_i\right)$ and by $(*)_2$ above, for $j < \kappa$ implies
$\mu \notin \text{ pcf}_{\sigma_1 \text{-complete}} 
\left(\dbcu_{i \in S_1} {\frak d}_i \cap  
\lambda_j \right)$; also $\mu = \text{ max pcf}({\frak d})$.  So we are 
left with clauses $(a) + (b)$.  Let ${\frak d}_\delta 
= \{\lambda_{\delta,\zeta}:
\zeta < \zeta_\delta < \sigma_1\}$, as $\sigma_1 < \kappa$, clearly 
without loss of generality for some $\zeta(*)$

$$
S_3 = \{\delta \in S_2:\zeta_\delta = \zeta(*)\} \text{ is stationary}.
$$
\mn
For each $\zeta < \zeta(*)$ let

$$
{\Cal P}_\zeta =: \{S \subseteq S_2:\text{max pcf}\{\lambda_{\delta,\zeta}:
\delta \in  S\} < \mu\}.
$$
\mn
If for some $\zeta < \zeta(*)$, the normal ideal on $\kappa$ which  
${\Cal P}_\zeta$ generates is proper, we have finished.  If not, for each  
$\zeta < \zeta(*)$ there are members $S_{\zeta,i}(i < \kappa)$ of 
${\Cal P}_\zeta$ and club $C_\zeta$ of $\kappa$  such that:

$$
\delta \in S_3 \cap C_\zeta \Rightarrow \dsize \bigvee_{i < \delta} \delta \in
S_{\zeta,i}.
$$
\mn
Clearly $C = \dbca_{\zeta < \zeta(*)} C_\zeta$ is a club of $\kappa$,  
now remembering $S_3 \subseteq S_0$ we know that on $S_3 \cap C$ the
function $\delta 
\mapsto \text{ sup}\{\text{Min}\{i:\delta \in S_{\zeta,i}\}:\zeta < \zeta
(*)\}$ is a pressing down function, so for some stationary 
$S_4 \subseteq S_3$ and ordinal $j(*) < \text{ Min}(S_4) < \kappa$  
we have:

$$
\delta \in S_4 \Leftrightarrow \dsize \bigwedge_{\zeta < \zeta(*)} \,
\dsize \bigvee_{j < j(*)} \delta \in C \cap S_{\zeta,j}.
$$
\mn
But as cov$(|j(*)|,\sigma_2,\sigma_1,2) < \kappa$, there is $w \subseteq 
j(*)$ such that $|w| < \sigma_2$ and $S_5$ is a stationary subset of  
$\kappa$ where

$$
S_5 = \bigl\{ \delta \in S_3:\dsize \bigwedge_{\zeta < \zeta(*)} \,
\dsize \bigvee_{i \in w} \delta \in S_{\zeta,i} \bigr\}.
$$
\mn
Let ${\frak b} = \{\lambda_{i,\zeta}$: for some $j$ and $\zeta$ we have 
$i \in S_{\zeta,j},j \in w$ and $\zeta < \zeta(*)\}$, it is the union of  
$< \sigma_2$ sets $({\frak b}_{j,\zeta} = 
\{\lambda_{i,\zeta}:i \in S_{\zeta,j}\}$ 
for $j \in w,\zeta < \zeta(*))$, each with max pcf $< \mu$.

This contradicts the fact (\scite{1.2A}).
\enddemo
\bigskip

\demo{Proof of Fact 1.2A}  Similar to the proof of \cite[1.5]{Sh:371}.
\hfill$\square_{\scite{1.2}}$\margincite{1.2}
\enddemo
\newpage

\head {\S2 On measures of the size of ${\Cal D}_{< \kappa}(\lambda)$}
\endhead  \resetall 
\bn
SAHARON:  on 2.1 see 430? \nl
Improving a little \cite[5.9]{Sh:400}.
\proclaim{\stag{2.1} Claim}  Assume $\lambda > {\text{\rm cf\/}}(\lambda)
= \aleph_0,\lambda > 2^\theta$ and 
$[\lambda' < \lambda \and {\text{\rm cf\/}}(\lambda') \leqq  
\theta \Rightarrow {\text{\rm pp\/}}_\theta
(\lambda') < \lambda]$ and ${\text{\rm pp\/}}(\lambda) < 
{\text{\rm cov\/}}(\lambda,\lambda,\aleph_1,2)$.  \ub{Then} 
$\{\mu:\lambda < 
\mu = \aleph_\mu \leq {\text{\rm pp\/}}_\theta(\lambda)\}$ has order type  
$\geq \theta$. 

We shall return to this in \cite[\S1]{Sh:430} so we do not elaborate.
\endproclaim
\bigskip

\proclaim{\stag{2.2} Claim}  Suppose $\theta,\kappa$ are regular, $\aleph_0 < 
\theta \leq \kappa < \lambda \leq \lambda_1 \leq \lambda^*$, and  
$(\forall \mu)[\lambda \leq \mu \leq \lambda_1 \and {\text{\rm
cf\/}}(\mu) 
< \theta \Rightarrow {\text{\rm pp\/}}_{<\theta}(\mu) \leqq \lambda^*]$ 
and ${\text{\rm cov\/}}(\kappa,\theta,\theta,2) = \kappa$. 
\ub{Then} 
there is a family ${\Cal P}$ of $\leq \lambda^*$ subsets of $\lambda_1$,
each of cardinality $< \lambda$ such that: 
\mr
\item "{$(*)_1$}"  for every $Y \subseteq \lambda_1,|Y| \le \kappa$ there 
are $Z_n(n < \omega)$ such that:  \nl
$Y \subseteq \dbcu_{n < \omega} Z_n,|Z_n| = \kappa$ and for each $n$
$$
(\forall Z)[Z \subseteq Z_n \and |Z| < \theta \Rightarrow (\exists X \in  
{\Cal P})[Z \subseteq X]],
$$ 
\sn
\item "{$(*)_2$}"  for every $Y \subseteq \lambda_1,|Y| \ge \theta_1$
and $\aleph_1 < {\text{\rm cf\/}}(\theta_1) \and \theta_1 \le \theta$ 
there is $X \in {\Cal P}$ such that:  $X \cap Y$  
has cardinality $\geq \theta_1$.
\endroster
\endproclaim
\bigskip

\remark{Remark}   1) Here and later we can replace $\leq \lambda^*$ by $< 
\lambda^* = {\text{\rm cf\/}}(\lambda^*)$.  [Saharon: check.] \nl
2) See \cite[x.x.]{Sh:430}.
\endremark
\bigskip

\demo{Proof}  It suffices to prove $(*)_1$ as $(*)_2$ follows.  Let  
$\Theta = \{\mu:\lambda \leq \mu \leq \lambda_1$ and cf$(\mu) < \theta\}$.  
Clearly if $\Theta = \emptyset$, the conclusion is straightforward (by 
induction on $\lambda_1$). 

Without loss of generality $\lambda^* = \text{ sup}\{\text{pp}_{<\theta}(\mu):
\mu \in \Theta\}$.  Now each pp$_{<\theta}(\mu)$ (for $\mu \in \Theta$)  
has cofinality $\geq \theta$, and if $\delta < \theta,\langle \mu_i:i < 
\delta \rangle$ increasing, $\langle \text{ pp}_{<\theta}(\mu_i):i < \delta 
\rangle$ strictly increasing then pp$_{<\theta}(\dbcu_{i < \delta} \mu_i)  
> \dsize \sum_{i < \delta} \text{ pp}_{<\theta}(\mu_i)$, hence  
cf$(\lambda^*) \geq \theta$ and, by \cite[2.3]{Sh:355},  without loss of 
generality $\lambda^* = \lambda_1$.  Let $\chi$ be regular large enough and  
$B$ be an elementary submodel of 
$({\Cal H}(\chi),\in,<^*_\chi)$ of cardinality  
$\lambda^*$ such that $\lambda^* + 1 \subseteq B$. \nl 
Let

$$
{\Cal P} = \{X \in B:X \subseteq \lambda^* \text{ and } |X| < \lambda\}.
$$
\enddemo
\bn
Now repeat the proof of \cite[3.5]{Sh:400}, noting
\demo{\stag{2.2A} Observation}  Suppose: 
\mr
\widestnumber\item{$(iii)$}
\item "{$(i)$}"  ${\frak a}$ is 
a set of regular cardinals $|{\frak a}| \leq \kappa < \text{ Min}({\frak a})$;
\sn
\item "{$(ii)$}"  the function ${\frak b} \mapsto \langle f^{\frak
b}_\alpha:\alpha < \text{ max pcf}({\frak b}) \rangle$ (for ${\frak b}
\subseteq {\frak a}$) is as in \cite[\S1]{Sh:371}; 
i.e. satisfies:
{\roster
\itemitem{ $(*)_1$ }   $f^{\frak b}_\alpha \in \Pi {\frak b}$, we stipulate  
$f^b_\alpha \restriction ({\frak a} \backslash {\frak b}) = 
0_{{\frak a} \backslash {\frak b}}$ and  
$\langle f^{\frak b}_\alpha:\alpha < \theta \rangle$ is 
$<_{J_{<\text{max pcf}[{\frak b}]}}$-increasing and cofinal in 
$(\Pi{\frak b},<_{J_{<\text{max pcf}[b]}})$, and 
$(\forall g \in \Pi {\frak b})(\exists^\theta \alpha < 
\theta)[g < f^\theta_\alpha]$
\sn
\itemitem{ $(*)_2$ }  if $\delta < \theta \in \text{ pcf}({\frak a})$,
cf$(\delta) \in (|{\frak b}|,\text{Min } {\frak b}]$ 
and $\sigma \in {\frak a} \backslash \{{\text{\rm cf\/}}
(\delta)\}$ then $f^{\frak b}_\delta(\sigma) = 
\text{ Min}\{ \dbcu_{\alpha \in C} f^{\frak b}_\alpha(\sigma):C$ 
a club of $\delta\}$
\endroster}
\item "{$(iii)$}"  ${\frak a} \subseteq \lambda^*,
\lambda \leq \lambda^*$, and \nl  
$(\forall {\frak b})[{\frak b} \subseteq {\frak a} \and |{\frak b}| 
< \theta \and {\text{\rm sup\/}}({\frak b}) \geq \lambda
\and \text{ sup}({\frak b}) \notin {\frak b} \Rightarrow 
\text{ max pcf}_{J^{bd}_{\frak b}} {\frak b} \leq \lambda^*]$
\sn
\item "{$(iv)$}"   We let for $\theta \in \text{ pcf}({\frak a})$, 
and $\alpha < \theta:f^\theta_\alpha = 
f^{{\frak b}_\theta[{\frak a}]}_\alpha$.
\ermn
\ub{Then} for every $\{{\frak a}_i:i < \kappa\} 
\subseteq \{{\frak b}:{\frak b} \subseteq {\frak a},|{\frak b}| < 
\theta\},g_i \in \Pi {\frak a}_i(i < \kappa)$ we 
can find $g$ and $\lambda'_i < 
\lambda$ such that:  $g \in \Pi {\frak a}$, for 
each $i < \kappa,g_i \leq g$ and 
for every $i < \kappa$ there are $\lambda'_i < \lambda$ such that letting  
${\frak c}^i = {\frak a}^i \backslash \lambda'_i$ 
we have $g \restriction {\frak c}^i$ is Max of 
finitely many functions from $\{f^\mu_\alpha \restriction c^i:i < \kappa,
\alpha < \mu,\mu \in {\text{\rm pcf\/}}({\frak a})$ 
and $\mu \leq \lambda^*\}$.  
Moreover for some ${\frak d} \subseteq [\lambda,\lambda^*] \cap 
{\text{\rm pcf\/}}({\frak a})$, for every 
$\theta \in \text{ Reg } \cap [\lambda,\lambda^*] \cap {\frak d}$, 
for some $\lambda_\theta < \lambda,g \restriction ({\frak
b}_\theta[{\frak a}] \backslash \lambda_\theta)$ is 
(the suitable restriction of $a$) Max of finitely 
many functions from $\{f^\sigma_{g(\sigma)}:\sigma \in (\lambda,\lambda^*] 
\cap {\text{\rm pcf\/}}({\frak a} \backslash \lambda')$ 
for every $\lambda' < \lambda\}$.
\enddemo
\bigskip

\demo{Proof}  Without loss of generality $\kappa^+ < \text{ Min}({\frak a})$.
Use \cite[1.4]{Sh:371} with $\kappa^+,{\frak a},\langle {\frak a}_i:
i < \kappa \rangle$, here standing for $\delta,{\frak a},
\langle {\frak b}_i:i < \zeta^* \rangle$ there and we get  
$\langle < \lambda_{i,\ell},{\frak c}_{i,\ell}:
\ell \leq n(i)>:i < \kappa \rangle$.  
Let $\lambda'_i = \sup[\cup \{{\frak c}_{i,\ell}:{\frak c}_{i,\ell}$ a
bounded subset of $\lambda\}]$.  
\hfill$\square_{\scite{2.2A}} \,\, \square_{\scite{2.2}}$
\enddemo
\bn
Similarly
\proclaim{\stag{2.3} Claim}  1) In \scite{2.2}, if $\theta = \kappa$ and  
$(\forall \mu )[\kappa \leq \mu \leq \lambda_1 \and \, 
{\text{\rm cf\/}}(\mu) < \theta  
\Rightarrow$ {\rm pp}$_{<\theta}(\mu) \leq \lambda^*]$ (i.e.  
$\kappa = \lambda$) \ub{then} we can add 
\mr
\item "{$(*)_3$}"  for every $Y \subseteq \lambda,|Y| \leq \kappa$ there 
are $Z_n,\langle Z_{n,i}:i < \theta \rangle$ for $n < \omega$ such that  
$Y \subseteq \dbcu_{n < \omega} Z_n,Z_n = \dbcu_{i < \theta} Z_{n,i},
|Z_{n,i}| < \theta,\langle Z_{n,i}:i < \theta \rangle$ increasing continuous 
and each $Z_{n,i}$ belongs to ${\Cal P}$.
\endroster
\endproclaim
\bn
Hence
\demo{\stag{2.4} Conclusion}  If $\lambda > \text{ cf}(\lambda) = \aleph_0$ 
\ub{then} there is a family ${\Cal P}$ of cardinality 
$\le \sup\{pp(\mu):\mu \leq \lambda$,
cf$(\mu) = \lambda_0\}$ consisting of countable subsets of 
$\lambda$ such that: 
\mr
\item "{$(*)$}"  if $Y \subseteq \lambda,|Y| = \aleph_1$ \ub{then} for some  
$Z \in {\Cal P},Y \cap Z$ is infinite.   
Moreover, we can find $\alpha^n_i \, (n < \omega,i < \omega_1)$ such that  
$Y \subseteq \{\alpha^n_i:n < \omega,i < \omega_1\}$ and for each $n$  
for arbitrarily large $i < \omega_1,\{\alpha^n_j:j < i\} \in {\Cal P}$.
\endroster
\enddemo
\bigskip

\demo{\stag{2.5} Conclusion}  If $\lambda > \kappa \geq \text{ cf}(\lambda)$,  
and pp$_\kappa(\lambda) < \text{ cov}(\lambda,\lambda,\kappa^+,2)$ \ub{then}
pp$_\kappa(\lambda) < \text{ cf}_\kappa(\Pi(\lambda \cap \text{ Reg}
\backslash \kappa^+),<_{J^{bd}_\lambda})$.
\enddemo
\bigskip

\demo{Proof}  Let $\chi = \beth_3(\lambda)^+$, and for $\zeta \leq \kappa^+$ 
let ${\frak B}_\zeta$ be an elementary submodel of $({\Cal
H}(\chi),\in,<^*_\chi)$ of 
cardinality pp$_\kappa(\lambda)$ such that pp$_\kappa(\lambda) \subseteq  
{\frak B}_\zeta,{\frak B}_\zeta$ 
increasing continuous, and $\langle {\frak B}_\epsilon:\epsilon  
\leqq \zeta \rangle \in {\frak B}_{\zeta +1}$.  
Let ${\frak B} = {\frak B}_{\kappa^+}$.  Assume 
that the conclusion fail and we shall prove that cov$(\lambda,\lambda,
\kappa^+,2) \leqq$ pp$_\kappa(\lambda)$, in fact that ${\Cal P} =: 
{\frak B}_{\kappa^+} \cap {\Cal S}_{<\lambda}(\lambda)$ exemplify it.
Let $a \subseteq \lambda,|a| \leq \kappa$ and we shall find $A \in {\Cal P}$
such that $a \subseteq A$; this suffice.  Choose by induction on 
$\zeta < \omega_1,f_\zeta \in {\frak B} 
\cap \Pi(\text{Reg } \cap \lambda)$ such 
that letting $N_\zeta$ be the Skolem Hull of 
$a \cup \{f_\xi:\xi < \zeta\} \cup \{i:i \le \kappa\}$,  
we have: for every large enough $\sigma \in \text{ Reg } \cap \lambda,
[\sigma \in N_\zeta \Rightarrow \sup(\sigma \cap N_\zeta) 
< f_\zeta(\sigma)]$.
Now use \scite{2.4} + proof of \cite[3.5]{Sh:400}.  
\hfill$\square_{\scite{2.5}}$\margincite{2.5} 
\enddemo
\bn
\centerline {$* \qquad * \qquad *$}
\bn
We now return to the issue of strong covering (from 
\cite[Ch.XIII,\S1-\S4]{Sh:b}) (better version \cite[Ch.VII]{Sh:g}).  It
influenced the first proof of a bound on $\aleph^{\aleph_0}_\omega$, and is 
clearly related to computing Min$\{|S|:S \subseteq {\Cal S}_{\leqq \kappa}
(\lambda)$ is stationary$\}$.
\proclaim{\stag{2.6} Lemma}  Suppose $\bold W \subseteq \bold V$ 
is a transitive class of $\bold V$ 
including all the ordinals and is a model of ZFC.
\nl
1)  For every set $Y \in \bold V$ of ordinals of cardinality 
$< \kappa$ (in $\bold V)$ there are $Y_n \in \bold W \, n < \omega$ 
(so we know only that $\langle Y_n:n < \omega \rangle  
\in \bold V!),\bold W \models 
``Y_n$ a set of $< \kappa$ ordinals" such that  
$Y \subseteq \dbcu_{n < \omega} Y_n$ provided that: 
\mr
\item "{$(*)_\kappa\,\,(i)$}"  $\kappa$ is a regular uncountable cardinal 
in $\bold V$,
\sn
\item "{$(ii)$}" if $a \in \bold V$ is a subset of 
Reg$^{\bold W} \backslash \kappa,|a| < \kappa$ and $g \in (\Pi
a)^{\bold V}$ then     
{\roster
\itemitem{ $\otimes$ }   there is a function $h \in \bold W$ such 
that $\theta \in \text{ Dom}(g) \Rightarrow g(\theta) < h(\theta) < 
\theta$ (so Dom$(g) \subseteq \text{ Dom}(h)$)
\endroster}
\ermn
or even just
\mr
\item "{$(*)^-_\kappa$}"  like $(*)_\kappa$ but in (ii) we demand only:     
{\roster
\itemitem{ $\otimes^-$ }  there are functions 
$h_n \in \bold W$ (for $n < \omega$) such that 
$(\forall \theta \in \text{ Dom}(g))[\dsize \bigvee_{n < \omega}
g(\theta) < h_n(\theta) < \theta]$.
\endroster}
\ermn
2)  For every set $Y \in \bold V$ of ordinals of cardinality $<
\kappa$ (in $\bold V$) there is $Z \in \bold W$ satisfying 
$\bold W \models ``Z$ a set of $< \kappa$ ordinals" such that  
$Y \subseteq Z$ provided that:

$$
(*)_\kappa + \aleph^{\bold V}_2 \leq \kappa.
$$
\mn
3)  Assume that $\kappa = \aleph^{\bold V}_1$ and $(*)_\kappa$ holds and 
\mr
\item "{$\oplus_0$}"   $\bold V \models ``A$  a set of ordinals of power  
$\kappa" \Rightarrow (\exists B \in \bold W)
[A \cap B$ infinite $\and \bold W \models ``|B| < \kappa"]$.
\ermn
\ub{Then} the conclusion of part (2) holds. \nl
4)  Assume $\kappa = \aleph^{\bold V}_1,(*)_\kappa,(\kappa^+)^{\bold
V} = (\kappa^+)^{\bold W}$ and $(*)_{\kappa^+}$.  \ub{Then} 
the conclusion of part (2) holds. \nl
5)  Assume 
\mr
\item "{$(a)$}"  $(\kappa^+)^{\bold V} = 
(\kappa^+)^{\bold W},(*)_\kappa,(*)_{\kappa^+}$;
\sn
\item "{$(b)$}"  there is $\langle C_\delta:\delta \in (\lambda + 1)
\backslash (\kappa + 1),\aleph_0 \leq \text{ cf}^{\bold W} \delta \leq
\kappa \rangle \in \bold W$ satisfying $C_\delta$ a club of $\delta$
for each $\delta$ such that $[\alpha \in \text{ acc}(C_\delta) 
\Rightarrow C_\alpha = C_\delta \cap \alpha]$ and otp$(C_\delta) \leq  
\kappa$.  [nec?]
\ermn
\ub{Then} $(\bold W,\bold V)$ 
satisfies the $\kappa$-strong covering (see \cite[Ch.VII]{Sh:g}) which
means:
\mr
\item "{$\otimes$}"   for every ordinal $\alpha$ and model $M \in
\bold V$ with 
universe $\alpha$, with countable vocabulary, there is $N \prec M$ of power  
$< \kappa,N \cap \kappa$ an ordinal and the universe of $N$ belongs to
$\bold W$.
\ermn
6)  Moreover, in (5) 
\mr
\item "{$\oplus^+$}"  in the game where a play last $\kappa$ moves, in the 
ith move (for $i < \kappa$) the first and second players 
choose $a_i,b_i \in [\lambda]^{< \kappa}$, 
respectively preserving $\dbcu_{j<i} b_j 
\subseteq a_i \subseteq b_i$, the first player has a winning strategy 
where the first player winning a play means $\{\delta < \kappa:
\dbcu_{i < \delta} a_i \in \bold W\} \in {\Cal D}_\kappa$ 
(= the clubs filter on $\kappa$ (in $\bold V$)).
\endroster
\endproclaim
\bigskip

\remark{Remark}  1) Note that parts (3), (4) hold for $\kappa \geq
\aleph^{\bold V}_1$,  but this is covered already by part (2). \nl
2) Note that in part (9), $\aleph^{\bold V}_2 = \aleph^{\bold W}_2$.
\endremark
\bigskip

\demo{Proof}  Should be straightforward (if you read till here).  
[Originally we say only; for (1) imitate the proof of \cite[3.5]{Sh:400}, 
for (2) --- repeat the proof of \cite[3.5]{Sh:400} by doing the
induction for $i < \aleph_1$,  then use part (1).  
For (3) --- instead using part (1) in the end, use the assumption, 
for (4), (5) imitate the proof of \cite[3.6]{Sh:400}.]

For the proof of 1)-5), let 
$Y$ be a subset of the ordinal $\lambda$, a cardinal of $\bold V$  
(for part (5) $\lambda$ is given), and let $\chi =:
[(2^\lambda)^+]^{\bold W}$.  
Let $({\Cal H}(\chi)^{\bold W},\in,<^*_\chi) \in \bold W$. \nl
1)  In $\bold V$ we choose by 
induction on $n < \omega,N_n,\alpha_n,h_{n,\ell},g_n$ such that:
\mr
\item "{$(a)$}"  $N_n \prec ({\Cal H}(\chi)^{\bold W},\in,<^*_\chi)$
\sn
\item "{$(b)$}"  $\bold V 
\models \|N_n\| < \kappa$ and $N_n \cap \kappa = 
\alpha_n$;
\sn
\item "{$(c)$}"  $Y \subseteq N_0$ and $\{\kappa,\lambda\} \in N_0$ 
\sn
\item "{$(d)$}"  $N_n \prec N_{n+1}$ 
\sn
\item "{$(e)$}"  $h_{n,\ell} \in \bold W$ is a 
(partial) function (for $\ell < 
\omega)$ and from $\lambda$ to $\lambda$ 
\sn
\item "{$(f)$}"  $g_n$ is a function, Dom$(g_n) = 
(\lambda \cap \text{ Reg}^{\bold W} \backslash \kappa) \cap N_n,g_n(\theta) =: \sup(N_n \cap \theta)$
\sn
\item "{$(g)$}"  for every $\theta \in \text{ Dom}(g_n)$ for some $\ell,g_n
(\theta) < h_{n,\ell}(\theta)$
\sn
\item "{$(h)$}"  $h_{n,\ell} \in N_{n+1}$ for $\ell < \omega$.
\ermn  
There is no problem to carry the definition.  Let $N = N_\omega =: 
\dbcu_{n < \omega} N_n$ and $\alpha^*_\omega = \dbcu_{n < \omega} \alpha_n = 
N \cap \kappa$.  For $\alpha \leq \omega$ let $M_\alpha$ be defined as the 
Skolem Hull in $({\Cal H}(\chi)^{\bold N},\in,<^*_\chi)$ 
of $\{i:i < \alpha_n\} \cup \{h_{n,\ell }:
n,\ell < \alpha\} \cup \{\kappa,\lambda\}$.  
Clearly $[\alpha < \omega \Rightarrow M_\alpha \in \bold W]$, and
$[\alpha \le \omega \Rightarrow M_\alpha \prec N_\alpha]$ and $\bold V
\models ``\|M_\alpha\| 
< \kappa"$.  Now $M = M_\omega =  
\dbcu_{n < \omega} M_n \prec N,M \cap \kappa = 
\alpha_\omega = N \cap \kappa$ and for every $\theta \in (\lambda^+ \cap
\text{ Reg}^{\bold W} \backslash \kappa) \cap M$ 
and $n$ for some $m < \omega,\theta
\in N_m \and m > n$ so for some $\ell < \omega$ we have 
$\sup(N_n \cap \theta) < 
h_{m,\ell}(\theta) \in M$ hence $\sup(N \cap \theta) = 
\underset {n < \omega} {}\to \sup(N_n \cap \theta) \leq \sup(M \cap \theta) 
\leq \sup(N \cap \theta)$.  So by \cite[3.3A]{Sh:400}, $M \cap \lambda = 
N \cap \lambda$, so $Y_n =: \lambda \cap M_n$ for $n < \omega$ are as 
required.

2)-5)  The following will be used in proving during the proof of
2) - 5).  Let  
$\delta(*) \leq \kappa$ be given (i.e. we shall choose it for each part) 
and we assume $(*)_\kappa$.  We fix $({\Cal H}(\chi)^{\bold
W},\in,<^*_\chi) \in \bold W$ as in part (1).

In $\bold V$ we choose by induction on $i < \delta(*),N_i,\alpha_i,h_i,g_i$ 
such that: 
\mr
\item "{$(a)$}"   $N_i \prec ({\Cal H}(\chi)^{\bold W},\in,<^*_\chi)$
\sn
\item "{$(b)$}"  $\bold V \models 
``\|N_i\| < \kappa"$ and $N_i \cap \kappa = \alpha_i$ 
\sn
\item "{$(c)$}"   $Y \subseteq N_0,\{\kappa,\lambda\} \subseteq N_0$ 
\sn
\item "{$(d)$}"  $N_i$ is increasing continuous
\sn  
\item "{$(e)$}"  $h_i \in \bold W$ 
is a partial function from $\lambda + 1$ to  
$\lambda +1$
\sn
\item "{$(f)$}"   $g_i$ is the function with Dom$(g_i) = (\lambda^+ \cap  
\text{ Reg}^{\bold W} \backslash \kappa) \cap N_i$ such that
$$
g_i(\theta) = \sup(N_i \cap \theta)
$$
\item "{$(g)$}"  $\theta \in \text{ Dom}(g_i) \Rightarrow  g_i(\theta) < 
h_i(\theta) < \theta$
\sn  
\item "{$(h)$}"  $h_i \in N_{i+1}$.
\ermn
There is no problem to carry the construction.  We let $N = N_{\delta(*)}
= \dbcu_{i< \delta(*)} N_i$.
\enddemo
\bigskip

\demo{Proof of Part (2)}  Choose $\delta(*) = \aleph^{\bold V}_1$, 
so under the assumptions of part (2) we have 
$\delta(*) < \kappa$, hence $\bold V \models ``\|N\|
< \kappa"$.  Apply part (1) to $Y = \{h_i:i < \delta(*)\}$ (which is  
$\subseteq \bold W$), so we can find $Y_n \in \bold W$ 
(and if you like $Y_n \subseteq H =: \{h \in \bold W:h$ a 
partial function from $\lambda +1$ to $\lambda +1\}$),  
for $n < \omega$ such that $\bold V 
\models ``|Y_n| < \kappa"$ and $Y \subseteq
\dbcu_{n < \omega} Y_n$; (well, $Y$ is not a set of ordinals, but we can 
code it as one).
So for some $n = n(*)$, we hvae $\bold V 
\models ``|Y \cap Y_n| = \aleph_1"$.  Let $M$ be 
the Skolem Hull in $({\Cal H}(\chi)^{\bold W},\in,<^*)$ 
of $\{\alpha:\alpha < 
\dbcu_{i < \delta(*)} \alpha_i\} \cup Y_{n(*)} \cup \{\kappa,\lambda\}$;  
so $M \in \bold W,\bold V \models 
``\|M\| < \kappa"$  and so $\bold W \models ``\|M\| < 
\kappa"$.  Let $M' = M \cap N$; so clearly $M' \prec N,
\dbcu_{i < \delta(*)} \alpha_i \subseteq M' \cap \kappa \subseteq N \cap  
\kappa = \dbcu_{i < \delta(*)} \alpha_i$.  Lastly, if $i < \delta(*)$ and  
$\theta \in (\lambda^+ \cap \text{ Reg}^{\bold W} 
\backslash \kappa) \cap M'$  
then for some $j \in (i,\delta(*)),\theta \in N_j$ and for some $\epsilon$ 
we have $j < \epsilon < \delta(*) \and h_\epsilon \in Y_{n(*)}$, so $\theta  
\in \text{ Dom}(h_\epsilon)$ and $\sup(N_i \cap \theta) \leq \sup(N_j \cap  
\theta) \leq \sup(N_\epsilon \cap \theta) < h_\epsilon(\theta) \in M'$.
So by \cite[3.3A]{Sh:400} we know $M' \cap \lambda = N \cap \lambda$, so 
$Y \subseteq M' \cap \lambda \subseteq M \cap \lambda \in \bold W$ and
$\bold W \models ``|M \cap \lambda| \le \|M\| < \kappa$ as required.
[Saharon: fill details]
\enddemo
\bigskip

\demo{Proof of Part (3)}  As in part (2), we choose $\delta(*) = 
\aleph^{\bold V}_1$ and get 
$N_i,\alpha_i,g_i,h_i$ (for $i < \delta(*)$).  By $\oplus_0$ applied 
to $\{h_i:i < \delta(*)\}$ (again translating to a set of ordinals) there 
is a set $B \in \bold W$ such that 
$\bold W \models ``|B| < \kappa"$ and $A \cap B$ 
infinite, without loss of generality  $B \subseteq H$ 
(see proof of part (2)).
So there is a limit ordinal $\zeta < \delta(*)$ such that $\zeta = 
\sup \{i < \zeta:h_i \in B\}$.  Now let $M$ be the Skolem Hull in  
$({\Cal H}(\chi)^{\bold W},\in,<^*_\chi)$ of 
$\{\alpha:\alpha < \dbcu_{i < \zeta} \alpha_i\} \cup
B \cup \{\kappa,\lambda\}$, and let $M' = M \cap N_\zeta$ so clearly
$M' \cap \lambda \subseteq M \cap \lambda \in \bold W,|M' \cap
\lambda| \le |M \cap \lambda| \le \|M\| < \kappa$ and as 
above $M' \cap \lambda = N_\zeta \cap \lambda$, so  
$N_\zeta \cap \lambda = M' \cap \lambda \in \bold W$ but 
$Y \subseteq N_\zeta \cap \lambda = M' \cap \lambda \subseteq M \cap
\lambda$ and are finished.
\enddemo
\bigskip

\demo{Proof of Part (4)}  We let $\delta(*) = \kappa$ and get $N_i,\alpha_i,
g_i,h_i$ (for $i < \delta(*)$) be as before.  Now we apply part (1) with
$\kappa^+ ( = (\kappa^+)^{\bold V} = (\kappa^+)^{\bold W}),
\{h_i:i < \kappa\}$ here standing for 
$\kappa$ and $Y$ there, and get $\langle Y_n:n < \omega \rangle$.
So for some $n(*)$ the set $\{i:g_i \in Y_{n(*)}\}$ is unbounded in $\kappa$.
Let $M$ be the Skolem Hull in $({\Cal H}(\chi)^{\bold
W},\in,<^*_\chi)$ of $\{\alpha:\alpha < 
\dbcu_{i < \delta(*)} \alpha_i\} \cup Y_{n(*)} \cup \{\kappa,\lambda\}$.  
As before $N \subseteq M,M \in \bold W$ and $\bold W 
\models ``\|M\| < \kappa^+"$.  So $\bold W \models 
``\|M\| = \kappa"$  hence there is a one to one function $f \in W$
from $\{i:i < \kappa\}$ onto $M$,  so for some club $E \in \bold V$ 
of $\kappa$ (in $\bold V$)

$$
i \in  E \Rightarrow  N_i \subseteq \text{ Rang}(f \restriction i).
$$
\mn
So for each $i \in E,Y \subseteq \text{ Rang}(f \restriction i)$ hence we 
are done. \nl
5)  So we have already chosen $N_i,\alpha_i,g_i,h_i$ for $i < \kappa$.
By parts (2) + (4) we can assume \wilog \, that Dom$(h_i)$ has
cardinality $< \kappa$.  Let $N_\kappa = 
\dbcu_{i < \kappa} N_i$, by part (1) there is a sequence $\langle Y_n:
n < \omega \rangle$ satisfying $Y_n \in \bold W$ such that:

$$
\align
\{h_i:i < \kappa\} \subseteq \dbcu_{n < \omega} Y_n \subseteq H = \{f \in 
\bold W:&f \text{ a partial function from} \\
  &(\lambda +1) \cap \text{ Reg}^{\bold W} 
\backslash \kappa \text{ to } \lambda + 1\}
\endalign
$$
\mn
and for each $n,\bold W \models 
``|Y_n| < \kappa^+"$.  So for some $n(*),\{i < 
\kappa:h_i \in Y_{n(*)}\}$ is unbounded in $\kappa$.  So in $\bold W$ 
there is a list $\langle f_i:i < \kappa \rangle$ of $Y_{n(*)}$.  
In $\bold V$, for each $i < \kappa$ let $j_i < \kappa$ be minimal such that:
\mr
\item "{$(a)$}"   $h_i \in \{f_\zeta:\zeta < j_i\}$;
\sn
\item "{$(b)$}"  if for some $\zeta < \kappa,f_i \leq h_\zeta$ (i.e.  
Dom$(f_i) \subseteq \text{ Dom}(h_\zeta)$ and $(\forall \theta \in
\text{ Dom}(f_i))[f_i(\theta) \leq h_\zeta(\theta)]$ then there is such  
$\zeta < j_i$ 
\sn
\item "{$(c)$}"   $j_i \geq i + 1$.
\ermn
Let $E = \{\zeta < \kappa:\text{for every } i < \zeta,j_i < \zeta$, and  
$\zeta$ is a limit ordinal$\}$.  Now for each $\zeta \in E$ note that 
\mr
\item "{$(*)_1(\alpha)$}"   $\{h_i:i < \zeta\} \subseteq \{f_i:i < \zeta\}$;
\sn
\item "{$(\beta)$}"   $\dsize \bigwedge_{i < \zeta} 
[\dsize \bigvee_{j < \kappa} f_i \leq h_j \leftrightarrow
\dsize \bigvee_{j < \zeta} f_i \leq h_j]$.
\ermn
As $[i < j \Rightarrow h_i < h_j]$ clearly we get:
\mr
\item "{$(\beta)'$}"   $\dsize \bigwedge_{i < \zeta}
[\dsize \bigvee_{j < \zeta} f_i \leq h_j \Leftrightarrow f_i \leq h_\zeta]$.
\ermn
Define for $\zeta \leq \kappa$ a function $f^*_\zeta$ as follows:
\mr
\item "{$(*)_2$}"  $\text{Dom}(f^*_\zeta) = 
\bigcup\{\text{Dom}(f_i):i < \zeta \text{ and } f_i \leqq h_\zeta\}$
\nl
$f^*_\zeta(\theta) = \sup\{f_i(\theta):i < \zeta \text{ and } f_i \leqq 
h_\zeta\}$.
\ermn
So clearly
\mr
\item "{$(*)_3$}"  for $\zeta \in E$, Dom$(f^*_\zeta) =
\cup\{\text{Dom}(h_i):i < \zeta\},f^*_\zeta(\theta) =
\sup\{h_i(\theta):i < \zeta\}$ hence Dom$(f^*_\zeta) = N_{\zeta'} \cap
(\lambda +1) \cap \text{ Reg}^{\bold W} \backslash (\kappa +1)$
\sn
\item "{$(*)_4$}"  $f^*_\zeta \in \bold W$, and even for 
$\langle f^*_\zeta:\zeta < \kappa$ limit$\rangle \in \bold W$.  
Let us define for $\zeta \leq \kappa$ a function $g^*_\zeta$:
$$
\text{Dom}(g^*_\zeta) = \dbcu_{i < \zeta} \text{ Dom}(h_i) \subseteq 
(\lambda +1) \cap \text{ Reg}^{\bold W} \backslash \kappa,
$$

$$
g^*_\zeta(\theta) = \dbcu_{i < \zeta} g_i(\theta).
$$
\ermn
Clearly for $\zeta$ limit: 
\mr
\item "{$(*)_5$}"  $\text{Dom}(g^*_\zeta) = 
(\lambda + 1) \cap \text{ Reg}^{\bold W} \cap N_\zeta \backslash
\kappa$
\ermn
and if $\zeta \in E$ 

$$
\align
\theta \in \text{ Dom}(g^*_\zeta) \Rightarrow g^*_\zeta(\theta) &= \bigcup 
\{g_i(\theta):i < \zeta,\theta \in N_i\} = \bigcup\{h_i(\theta):
i < \zeta,\theta \in N_i\} \\
  &= \bigcup\{f_i(\theta):i < \zeta \text{ and } f_i \leq h_\zeta\} = 
f^*_\zeta(\theta).
\endalign
$$
\mn
So for $\zeta \in E$ we have $g^*_\zeta \subseteq f^*_\zeta$ but by
$(*)_3 + (*)_5$ they have the same domain hence $g^*_\zeta = f^*_\zeta$.

For every $\theta \in N_\kappa \cap ((\lambda +1) \cap \text{
Reg}^{\bold W} \backslash (\kappa +1)),
g_\kappa(\theta)$ is an ordinal $\in (\kappa,\theta) 
\subseteq (\kappa,\lambda)$ of cofinality $\kappa$, so $C_{g_\kappa(\theta)}$ 
is a set of ordinals of order type $\kappa$; let $C_{g_\kappa(\theta)} = 
\{\alpha^\theta_\zeta:\zeta < \kappa\}$, (increasing); it is strictly 
increasing continuous and has limit $g_\kappa(\theta)$; also 
$\langle g^*_\zeta(\theta):\zeta < \kappa \rangle$ is also strictly 
increasing continuous with limit $g^*_\kappa(\theta)$.  Clearly 
$E_\theta = \{\zeta < \kappa:\zeta$ limit ordinal and $\alpha^\theta_\zeta = 
g^*_\zeta(\theta)$ (so $\theta \in N_\zeta)\}$ is a club of $\kappa$  
(in $V$; as $V \models ``\kappa$ regular uncountable").  So for
$\zeta \in E_\theta,C_{\alpha^\theta_\zeta} = C_{g_\kappa(\theta)} \cap 
\alpha^\theta_\zeta \subseteq N_{\zeta +1}$, hence $[\zeta \in \text{ acc } 
E_\theta \Rightarrow C_{g_\kappa(\theta)} \cap \alpha^\theta_\zeta \subseteq
N_\zeta]$.  Let $E^* = \{\zeta \in E:(\forall i < \zeta)(\forall \theta \in  
N_i)[\theta \in (\lambda +1) \cap \text{ Reg}^W \backslash (\kappa +1) 
\Rightarrow \zeta \in \text{ acc } E_\theta\}$.  So for $\zeta \in E^*,
\zeta$ is a limit ordinal and:  

$$
\align
\theta \in N_\zeta \cap ((\lambda +1) \cap \text{ Reg}^W \backslash 
(\kappa +1)) \Rightarrow C_{f^*_\zeta(\theta)} &= C_{g^*_\zeta(\theta)} \\
  &= C_{g^*_\kappa(\theta)} \cap g^*_\zeta(\theta) \subseteq N_\zeta.
\endalign
$$ 
\mn
Now we shall show that for $\zeta \in E^*,N_\zeta \cap \lambda \in W$. 
For $\zeta \in E^*$ we define by induction on $n < \omega,M^n_\zeta$:

$$
M^0_\zeta \text{ is the Skolem Hull in } ({\Cal H}(\chi)^{\bold W},\in,<^*)
\text{ of } \{\alpha:\alpha < \alpha_\zeta\}.
$$
\mn
Let $M^{n+1}_\zeta$ be the Skolem Hull in $({\Cal H}(\chi)^{\bold W},
\in,<^*_\chi)$ of

$$
M^n_\zeta \cup \bigcup \{C_{f^*_\zeta(\theta)}:\theta \in (\lambda +1) 
\cap \text{ Reg}^{\bold W} \backslash 
(\kappa +1) \text{ and } \theta \in M^n_\zeta  
\cap \text{ Dom } f^*_\zeta\}.
$$
\mn
Let $M_\zeta = \dbcu_{n < \omega} M^n_\zeta$, clearly $M_\zeta \in
\bold W$ (as $f^*_\zeta \in \bold W$) we can 
prove by induction on $n$ that $M^n_\zeta \subseteq
N_\zeta$, hence:  
$M_\zeta \subseteq N_\zeta$.  Also $\alpha_\zeta \subseteq  
M_\zeta \cap \kappa \subseteq N_\zeta \cap \kappa = \alpha_\zeta$, and: 
\mr
\item "{${{}}$}"  $\theta \in M_\zeta \cap ((\lambda +1) \cap \text{ Reg}^W
\backslash (\kappa +1)) \Rightarrow  C_{g^*_\zeta(\theta)}$ is an unbounded 
of subsets of both $N_\zeta \cap \theta,M_\zeta \cap \theta$.
\ermn
So by \cite[5.1A(1)]{Sh:400} we get $N_\zeta \cap \lambda= M_\zeta \cap
\lambda$.  Alternatively, let for $\zeta \in E,M^*_\zeta$ is the
Skolem hull of $\alpha_\zeta \cup \{f^*_i:i < \zeta\}$ in $({\Cal
H}(\chi)^{\bold W},\in,<^*_\chi)$; so clearly $M^*_\zeta \subseteq
N_\zeta$ and again by \cite[5.1A]{Sh:400} we have $M^*_\zeta \cap
\lambda = N_\zeta \cap \lambda$. \nl
6)  The winning strategy for the second player is to choose ``on the side" 
also $N_i,\alpha_i,h_i,g_i$ as in the common part of the proof of parts 
(2)-(5) and guaranteeing:  $a_i$ include $\dbcu_{j < i} (N_j \cap \lambda) 
\cup \dbcu_{j < i} b_j,N_{i+1}$ include $b_i$ and $N_i \cap \lambda$ is 
the universe of an elementary submodel of $({\Cal H}(\chi)^{\bold
W},\in,<^*_\chi) \in \bold W$.  \hfill$\square_{\scite{2.6}}$\margincite{2.6}
\enddemo
\bigskip

\remark{\stag{2.6A} Remark}  1)  We can put $\lambda$ as a parameter of the 
Lemma \scite{2.6}, then in $(*)_\kappa,a \subseteq (\lambda + 1) \cap  
\text{ Reg}^{\bold W} \backslash \kappa$, etc., 
(so we may write $(*)_\kappa$) and  
$Y \subseteq \lambda$ (in parts (1)-(4)) and $\alpha \leq \lambda$  (in 
$\otimes$ of part (5)). \nl
2)  Note that $(*)^\lambda_\kappa$ follows easily from the relevant 
covering property in \cite[Ch.VII]{Sh:g}: 
\mr
\item "{$(*)$}"   if $a \in \bold V,a \subseteq \lambda,
\bold V \models ``|a| < \kappa"$  
then for some $b \in \bold W,a \subseteq b,\bold W \models ``|b| < \kappa"$.
\endroster
\endremark
\newpage

\head {\S3  pcf Inaccessibility and Characterizing the Existence      
of Non $<$-decreasing Sequences (for a Topological Problem)} \endhead  \resetall 
\bigskip

\definition{\stag{3.1} Definition}  1) We say $\lambda$ is $(\mu,\theta,
\sigma)$-inaccessible if $\lambda > \mu \geq \theta > \sigma$ and for any  
${\frak a} \subseteq$ Reg we have: if $|{\frak a}| < \theta$, 
Min$({\frak a}) > \mu$ and ${\frak a} \subseteq
\lambda$, even $\sup({\frak a}) < \lambda$ and $I$ is a 
$\sigma$-complete ideal on ${\frak a}$,
\ub{then} $\lambda \nleq \text{ tcf}(\Pi {\frak a}/I)$ 
(when  tcf  is well defined).
\nl
2)  If we write $*$ instead of $\mu$ we mean ``for some $\mu,\theta \leq 
\mu < \lambda$". \nl
3)  If we omit $\sigma$ we mean $\sigma = \aleph_0$. \nl
4)  ``accessible" is just the negation of ``inaccessible".  \nl
\enddefinition
\bn
We now rephrase various old results.
\proclaim{\stag{3.2} Claim}  1) For $\lambda$ regular, in the 
definition, ``and $\sigma$-complete $I,\lambda \nleq {\text{\rm tcf\/}}
(\Pi {\frak a}/I)"$  can be replaced by $``\lambda \notin  
{\text{\rm pcf\/}}_{\sigma \text{-complete}}({\frak a})"$ and 
also by $``\lambda \ne {\text{\rm tcf\/}}(\Pi {\frak a}/I)$ 
for any ${\frak a} \subseteq \text{ Reg } \cap  
(\mu,\lambda),|{\frak a}| < \theta,I$ being $\sigma$-complete"; also if  
{\rm cf}$(\lambda) \notin [\sigma,\theta)$ then $``\sup({\frak a}) 
< \lambda"$ is not necessary just $``\lambda \notin {\frak a}"$.  \nl
2)  Assume $\lambda > \mu \geq \theta > \sigma$ and {\rm cf}$(\lambda) \geq  
\sigma$.  \ub{Then} $\lambda$ is $(\mu,\theta,\sigma)$-inaccessible iff  
$[\lambda' \in (\mu,\lambda) \and \sigma \leq$ {\rm cf}$(\lambda') < \theta  
\Rightarrow$ {\rm pp}$_{\Gamma(\theta,\sigma)}(\lambda') < \lambda]$.
\nl
3)  If $\lambda = {\text{\rm cf\/}}(\lambda) > \mu \geq \theta = 
{\text{\rm cf\/}}(\theta) > \sigma,\lambda$ is 
$(\mu,\theta,\sigma)$-accessible \ub{then} there is a 
set ${\frak a} \subseteq \text{ Reg } \cap (\mu,\lambda)$ of 
$(\mu,\theta,\sigma)$-inaccessible cardinals each $> \mu,{\frak a}$ 
of cardinality $< \theta$ such that $\lambda \in$ 
{\rm pcf}$_{\sigma \text{-complete}}({\frak a})$.
\nl
4) If $\lambda = {\text{\rm cf\/}}(\lambda) \geq \kappa > \mu \geq \theta > 
\sigma = {\text{\rm cf\/}}(\sigma)$, and $(\exists {\frak a})
[{\frak a} \in \text{ Reg } \cap (\mu,\kappa) \and |{\frak a}| 
< \theta \and \lambda \in {\text{\rm pcf\/}}_{\sigma
\text{-complete}}({\frak a})]$ \ub{then} there is a set ${\frak a}$
of $(\mu,\theta,\sigma)$-inaccessible cardinals $\in (\mu,\kappa)$
with $|{\frak a}| < \theta$ such that 
$\lambda = {\text{\rm max pcf\/}}({\frak a})$, and $\lambda \in  
{\text{\rm pcf\/}}_{\sigma \text{-complete}}({\frak a})$.  If $\kappa$ is 
$(\mu,\theta,\sigma)$-inaccessible \ub{then} necessarily 
$\sup({\frak a}) = \kappa,J^{bd}_{\frak a} \subseteq
J_{<\lambda}[{\frak a}]$.
If $\sigma = \aleph_0$, also there is a tree of 
cardinality $\leqq \kappa$ and $\geq \lambda$ \, 
((cf$\kappa$)-branches if $\kappa$ is $(\mu,\theta,\sigma)$-inaccessible).
\nl
5)  If $\lambda = \text{ max pcf}({\frak a}),\kappa = |{\frak a}| 
\leq \mu < \text{ Min}({\frak a})$, each $\theta \in {\frak a}$ 
is $(\mu,\kappa^+,2)$-inaccessible \ub{then} there is a tree 
of cardinality sup$({\frak a})$ and $\geq \lambda 
{\text{\rm cf(otp)\/}}$-branches.  
If we have $\kappa$ pairwise disjoint subsets of ${\frak a}$ not in  
$J_{<\lambda}[{\frak a}],2^\kappa \geq \sup({\frak a})$ 
or on each $\theta \in {\frak a}$ there is 
an entangled linear order \ub{then} there is an entangled linear order of 
cardinality $\lambda$. \nl
6)  If $\mu < \lambda < \text{ pp}^+(\mu)$, \ub{then} 
there is a tree with $\leq \mu$ nodes and $\geq \lambda$ branches.  
If $\mu$ is $(*,(\text{cf}\mu)^+,
2)$-inaccessible we can demand  $``\geq \lambda \, 
(\text{cf}\mu)$-branches".
\endproclaim
\bigskip

\demo{Proof}  1) Easy (using pcf analysis and \cite[\S1]{Sh:355}). \nl
2)  Easy, too (use \cite[2.4]{Sh:355}).  \nl
3)  Prove by induction on $\lambda$ using \cite[1.10]{Sh:345a} (so in 
\cite[1.12]{Sh:345a} we can replace pcf by pcf$_{\sigma \text{-complete}}$).
\nl
4)  Similar to (3). \nl
5)  By \cite[\S4]{Sh:355}. \nl
6)  Easy, too.
\enddemo
\bn
We state some variants of \cite[\S2,\S4]{Sh:400}; specifically combining 
\cite[2.4,4.2]{Sh:400}:
\proclaim{\stag{3.3} Claim}  Suppose: 
\mr
\widestnumber\item{$(iii)$}
\item "{$(i)$}"  $\langle \lambda_\zeta:\zeta < \zeta(*) \rangle$ is a 
strictly increasing sequence of regular cardinals $> \sigma$
\sn
\item "{$(ii)$}"  for $\zeta$ limit, $\lambda_\zeta = 
(\dbcu_{\xi < \zeta} \lambda_\xi)^+$
\sn
\item "{$(iii)$}"  $\lambda_\zeta \in {\text{\rm pcf\/}}({\frak a}_\zeta)$ 
\sn
\item "{$(iv)$}"  ${\frak a}_\zeta \subseteq \text{ Reg } \cap 
(\sigma^+,\lambda_0)$ and $|{\frak a}_\zeta| < \sigma,\sigma$ regular.
\ermn
\ub{Then} $\zeta(*) < \sigma^{+3}$.
\endproclaim
\bn
Similarly combining \cite[2.4,4.2]{Sh:400}
\proclaim{\stag{3.4} Claim}  Suppose: 
\mr
\widestnumber\item{$(iii)$}
\item "{$(i)$}"   $\langle \lambda_\zeta:\zeta < \zeta(*) \rangle$ is a 
strictly increasing sequence of regular cardinals $> \kappa$ 
\sn
\item "{$(ii)$}"  for $\zeta$ limit, $\lambda_\zeta = 
(\dbcu_{\xi < \zeta} \lambda_\xi)^+$ 
sn
\item "{$(iii)$}"  $\lambda_\zeta \in \text{ pcf}_{\sigma \text{-complete}}
(a_\zeta)$ where $\sigma$ is regular 
\sn
\item "{$(iv)$}"  ${\frak a}_\zeta 
\subseteq \text{ Reg } \cap (\kappa^+,\lambda_0),
|{\frak a}_\zeta| < \kappa,\kappa$ regular
\sn
\item "{$(v)$}"  ${\frak a}_\zeta \subseteq {\frak b}$ and if $\langle
\mu_\zeta:\zeta \leqq \kappa^{+2} \rangle$ is strictly increasing
sequence of regular ${\frak a} \subseteq {\frak b},|{\frak a}| <
\kappa,\{\mu_\zeta:\zeta < \kappa^{+2}\} \subseteq \text{ pcf}({\frak
a}),\mu_{\kappa^{+2}} = \text{ max pcf}\{\mu_\zeta:\zeta <
\kappa^{+2}\}$ and $(\lambda_0 + \dsize \sum_{\zeta < \kappa^{+2}}
\mu_\zeta)^+ < \mu_{\kappa^{+2}} < \dbcu_{\zeta <
\xi}\{\lambda_\zeta:\zeta < \kappa^{+3}\}$, then there are $\mu'_\zeta
< \mu_\zeta$ regular, as for $\zeta \leqq \kappa^{+2}$ and ${\frak a}'
\subseteq {\frak b}$ with the same properties and
$\Sigma\{\mu_\zeta:\zeta < \kappa^{+2}\} = \Sigma\{\mu'_\zeta:\zeta <
\kappa^{+2}\}$.    
\ermn
\ub{Then} $\zeta(*) < \text{ Max}\{\kappa^{+3},\text{cov}
(\dbcu_{\zeta < \xi} {\frak a}_\zeta,\kappa,\kappa,\sigma)^+\}$.
\endproclaim
\bigskip

\proclaim{\stag{3.5} Claim}  The following is impossible: 
\mr
\widestnumber\item{$(iii)$}
\item "{$(i)$}"   $\sigma < \kappa < \theta < \mu$ are regular,  
$\kappa^+ < \theta$ 
\sn
\item "{$(ii)$}"   $\langle \lambda_\zeta:\zeta < \mu \rangle$ is a 
strictly increasing sequence of regular cardinals $> \mu$ 
\sn
\item "{$(iii)$}"  $S = \{\epsilon < \mu:{\text{\rm cf\/}}(\epsilon)
= \theta \text{ and for some club } C \text{ of } \epsilon,
{\text{\rm sup pcf}}_{\sigma \text{-complete}}\{\lambda_\zeta:
\zeta \in C\} < \dsize \sum_{\zeta < \mu} \lambda_\zeta\}$ is stationary;
\sn
\item "{$(iv)$}"  $(a) \quad$ if $\delta < \mu$, cf$(\delta) 
= \kappa$ \ub{then} for every 
club $C$ of $\delta$, there is $\alpha \in C$ such that    
$$
{\text{\rm sup pcf\/}}_{\sigma \text{-complete}} \{\lambda_\zeta:\zeta \in  
\alpha \cap C\} \geq \dbcu_{\zeta \in C} \lambda_\zeta
$$  
or
\sn
\item "{${{}}$}"  $(b) \quad \lambda_\zeta \in  
\text{pcf}_{\sigma \text{-complete}}({\frak a}_\zeta),|{\frak
a}_\zeta| < \kappa$ and  
$\mu > {\text{\rm cov\/}}(\dbcu_{\zeta < \mu} {\frak a}_\zeta,\kappa,\kappa,
\sigma)$.
\endroster
\endproclaim
\bigskip

\proclaim{\stag{3.6} Claim}   Assume $|{\frak a}| < {\text{\rm
Min\/}}({\frak a})$, \ub{then} 

$$
{\text{\rm cf\/}}_{\le |{\frak a}|}(\Pi \,{\text{\rm pcf\/}}({\frak
a})) \leq {\text{\rm max pcf\/}}({\frak a}).
$$
\endproclaim
\bigskip

\demo{Proof}  More is proved in \cite[\S3]{Sh:371}. 
\enddemo
\bn
\centerline {$* \qquad * \qquad *$}
\bigskip

The following answers a question of Gerlits, Hajnal and Szentmiklossy in
\cite{GHS}.  They dealt with  $``\kappa$-good topological spaces $X"$  (i.e. 
every subset is the union of $\leq \kappa$ compact sets) and 
``weakly $\kappa$-good spaces" (every $Y \subseteq X$ of cardinality  
$> \kappa$  contains a compact subset of cardinality $> \kappa$). 
\cite{GHS} has the easy implication. \nl
We return to this in \cite[\S6]{Sh:513}.
\bigskip

\proclaim{\stag{3.7} Theorem}  The following conditions on $\kappa < \sigma 
< \theta$  are equivalent:  ($\kappa$ is an infinite cardinal, $\sigma$ and  
$\theta$ are ordinals) 
\mr
\item "{$(A)_{\kappa,\sigma,\theta}$}"  there are functions  
$f_\alpha:\kappa \rightarrow \sigma$ for $\alpha < \theta$ such that:
$$
\alpha < \beta \Rightarrow \dsize \bigvee_{i < \kappa} f_\alpha(i) < 
f_\beta(i)
$$ 
\sn
\item "{$(B)_{\kappa,\sigma,\theta}$}"  $2^\kappa 
\geq |\theta|$ or for every 
regular $\mu_1 \leq \theta$ for some 
singular cardinal $\lambda^* \leq \sigma$ we have:          
$$
{\text{\rm cf\/}}(\lambda^*) \leq \kappa,\lambda^* > 2^\kappa,
{\text{\rm pp\/}}^+(\lambda^*) > \mu_1.
$$
\endroster
\endproclaim
\bigskip

\demo{Proof}  First note
\enddemo
\bn
\demo{\stag{3.8} Observation}   Let $\kappa < \sigma,\kappa$ an infinite 
cardinal, $\sigma,\theta$ are ordinals.  If for every regular $\theta_1,  
\sigma \leq \theta_1 < \theta$ the statement
$(A)_{\kappa,\sigma,\theta_1}$ holds and 
$\theta$ is singular (e.g. $\theta > |\theta|$), \ub{then}
$(A)_{\kappa,\sigma,\theta}$ holds.
\enddemo
\bigskip

\demo{Proof}  We 
prove this by induction on $\theta$; if $\theta \leq \sigma$ 
-- trivial: use the constant functions.  As $\theta$ is singular $\theta  
= \dsize \sum_{\alpha < \theta(*)} \theta_\alpha$ where $\theta(*) < \theta,
\theta_\alpha < \theta,\theta_\alpha$ increasing continuous,  
$\theta_0 = 0$.  By the assumption for each $\alpha < \theta(*)$, there is 
a sequence $\langle f^\alpha_i:i < \theta_\alpha \rangle$ as required in 
$(A)_{\kappa,\sigma,\theta_\alpha}$, [why?  if $\theta_\alpha$ is 
singular by the induction hypothesis, if $\theta_\alpha$ is regular by an 
assumption of \scite{3.8}].  Similarly there is $\langle f_i:i < \theta(*)
\rangle$ exemplifying $(A)_{\kappa,\sigma,\theta(*)}$. 

For $i < \theta$ let $i = \dsize \sum_{\beta < \alpha (i)} \theta_\beta 
+ j(i),j(i) < \theta_{\alpha(i)},\alpha(i) < \theta(*)$ and let $g_i:\kappa  
\rightarrow 
\lambda$ be $g_i(2\zeta) = f_{\alpha(i)}(\zeta),g_i(2 \zeta +1) = 
f^{\alpha(i)}_{j(i)}(\zeta)$.  \hfill$\square_{\scite{3.8}}$\margincite{3.8}
\enddemo
\bigskip

\demo{Continuation of the Proof of 3.7}  First we do the easy direction.
\mn
\ub{$(B) \Rightarrow (A)$}
\mn
\ub{First Case}:  $2^\kappa \geq |\theta|$. 

Let $\{A_\alpha:\alpha < \theta\}$ be a family of $|\theta|$ distinct 
subsets of $\kappa$, let $A'_\alpha = \{2i:i \in A_\alpha\} \cup \{2i +1:i 
\notin A_\alpha,i < \kappa\}$ and let $f_\alpha:\theta \rightarrow \{0,1\} 
\subseteq \lambda$ be

$$
f_\alpha (i) = \cases 0 &\text{ if } i \in A'_\alpha \\
1 &\text{ if } i \notin A'_\alpha. \endcases 
$$
\mn
\ub{Second Case}:  $\lambda^* \leqq \sigma$,cf$(\lambda^*) \leq 
\kappa,\lambda^* > 2^\kappa$, 
pp$^+(\lambda^*) > \theta,\theta$ a regular cardinal. 

So there are regular cardinals $\lambda_i < \lambda^*$ for $i < \kappa$  
(such that $\lambda_i > \kappa$) and an ideal $J$ on $\kappa,\kappa \notin J$
such that $\dsize \prod_{i < \kappa} \lambda_i/J$ has true cofinality  
$\theta_1 \geq \theta$.  So there is a sequence $\langle f_\alpha:\alpha < 
\theta_1 \rangle,f_\alpha \in \dsize \prod_{i < \kappa} \lambda_i$ and  
$\alpha < \beta \Rightarrow f_\alpha < f_\beta$ mod $J$.  Now  
$\langle f_\alpha:\alpha < \theta \rangle$ is a sequence as required. 

By \scite{3.8} those two cases suffice.
\mn
\ub{$(A) \Rightarrow (B)$}
\mn
Let $\langle f_\alpha:\alpha < \theta \rangle$ be as in (A).

We can assume that $(B)$ fails, $\sigma$ minimal for which this occurs (for a
given $\kappa$ for some $\theta$) and $\theta$ minimal for the given  
$\kappa$ and $\sigma$.  So $\theta \geq (2^\kappa)^+$.  By Observation 
\scite{3.8}, $\theta$ is a regular cardinal.  So $2^\kappa < \theta$ (hence  
$2^\kappa < \sigma$) and $[a \subseteq \text{ Reg } \cap \sigma^+
\backslash \kappa^+ \and |a| \leq \kappa \Rightarrow \text{ max pcf } a < 
\theta]$, and $\sigma < \theta$. 

As $\theta$ is a regular cardinal necessarily cf$(\sigma) \leq \kappa$  
(otherwise for some $\sigma_1 < \sigma$ the set $\{\alpha < \theta:
\text{Rang}(f_\alpha) \subseteq \sigma_1\}$ is unbounded in $\theta$,  
contradicting the minimality of $\sigma$).  Also $\sigma$ is a limit 
ordinal as $2^\kappa < \theta = \text{ cf}(\theta)$ (as if $\sigma = \beta 
+ 1$,  for some $A \subseteq \kappa,B = \{\alpha:\dsize \bigwedge_{i <
\kappa} [i \in A \Leftrightarrow f_\alpha(i) = \beta]\}$ has cardinality  
$\theta$, so $\{f_\alpha \restriction (\kappa \backslash A):\alpha \in B\}$  
essentially contradicts the minimality of $\sigma$). 

Let $\chi$ be regular large enough.  We choose by induction on $i < 
(2^\kappa)^+$, a model $N_i$ such that:   

$N_i \prec ({\Cal H}(\chi),\in,<^*_\chi)$   

$\|N_i\| = 2^\kappa$   

$2^\kappa \subseteq N_0$   

$\kappa,\sigma,\theta \in N_0,\langle f_\alpha:\alpha < \theta \rangle \in  
N_0$   

$i < j \Rightarrow  N_i \prec  N_j$   

$\langle N_j:j \leq i \rangle \in N_{i+1}$   

$N_i$ increasing continuous. 

Let $\delta_i =: \sup(\theta \cap N_i)$ so $\langle \delta_i:i < 
(2^\kappa)^+ \rangle$ is strictly increasing continuous (as $\theta$ is 
regular, $\theta > \sigma$ and $\sigma > 2^\kappa$, necessarily $\delta_i < 
\theta$).  We define for $i < (2^\kappa)^+$, a function $g_i \in {}^\kappa 
\sigma$ by

$$
g_i(\zeta ) = \text{ Min}(N_i \cap \sigma \backslash f_{\delta_i}(\zeta)) 
$$
\mn
(it is well defined as $\sigma \in N_0 \subseteq N_i$ and $N \cap \sigma$ is 
unbounded in $\sigma$ as cf $\sigma \leq \kappa$). 

Now $i < (2^\kappa)^+$, cf$(i) = \kappa^+$ implies $N_i = 
\dbcu_{j < i} N_j$ and Rang $g_i \subseteq \dbcu_{j < i} N_j$ hence  
$\dsize \bigvee_{j < i} [\text{Rang}(g_i) \subseteq N_j]$; but every subset 
of $N_j$ of cardinality $\leq \kappa$ belongs to $N_{j+1}$, hence  
$g_i \in \dbcu_{j < i} N_j$.  So by Fodor 
Lemma for some stationary subset $S$
of $\{i < (2^\kappa)^+:\text{cf}(i) = \kappa^+\}$ and some $g^*:\kappa  
\rightarrow \sigma$ and some $A \subseteq \kappa$ and some $i(*) < 
(2^\kappa)^+$ we have:  $[i \in  S \Rightarrow g_i = g^*],(\forall i \in S)
(\forall \zeta < \kappa)[f_{\delta_i}(\zeta) = g^*(\zeta) \Leftrightarrow  
\zeta \notin A]$ and $g^* \in N_{i(*)}$; note $A \in N_0 \subseteq N_{i(*)}$ 
as $A \subseteq \kappa$.

Clearly $i \in S \and \zeta \in A \Rightarrow \text{ cf}[g^*(\zeta)] > 
2^\kappa$ (otherwise $g^*(\zeta) = \sup(N_i \cap g^*(\zeta))$ (as  
$N_i \prec ({\Cal H}(\chi),\in,
<^*_\chi),2^\kappa + 1 \subseteq N_0 \subseteq N_i)$  
and easy contradiction).  Also, as the $f_\alpha$'s are pairwise distinct, 
clearly $A \ne \emptyset$.
\bn
\ub{Question}:  What is cf$[\dsize \prod_{\zeta \in A} \text{ cf}[g^*(\zeta)]
]$?

(I.e. cofinality of the partial ordered set). 

By \cite[3.1]{Sh:355} it is  max pcf$\{\text{cf}[g^*(\zeta)]:\zeta \in A\}$,  
which by an assumption is $< \theta$, so there is a family $G \subseteq
\dsize \prod_{\zeta \in A} g^*(\zeta)$ of cardinality $< \theta$ such that  
$(\forall f \in \dsize \prod_{\zeta \in A} g^*(\zeta)) (\exists g \in G)
[f < g]$.  As the parameters in the demands on $G$ belongs to $N_{i(*)}$,  
without loss of generality $G \in N_{i(*)}$. 

Now we can define a partial function $H$ from the family $G$ to $\theta$: 

\ub{if} $g \in G$ and for some $\alpha$ the condition $(*)$ below holds 
then $H(g)$ is such an ordinal 

\ub{if} $g \in G$ and for no $\alpha$ the condition $(*)$ below holds then
$H(g)$ is not defined where 
\mr
\item "{$(*)$}"  $\alpha < \theta,f_\alpha \restriction (\kappa \backslash A) 
= g^* \restriction (\kappa \backslash A)$ and $g = g \restriction A \leq  
f_\alpha \restriction A \leq g^* \restriction A$.
\ermn
Now we can choose an ordinal $j(*)$ such that

$$
i(*) < j(*) < (2^\kappa)^+, \,\,\, j(*) \in S
$$
\mn
(possible as $S$ is a stationary subset of $(2^\kappa)^+$).

We know that there is a function $h \in G$ such that  
$f_{\delta_{j(*)}} \restriction A < h$.
\bn
\ub{Question}:  Is $H(h)$ well defined?
\bn
\ub{Possibility A}:  The answer is yes. 

Then $H(h) < \cup \{H(g) + 1:g \in (\text{Dom}(H)) \subseteq G\}$.  This 
union is an ordinal $< \theta$ (as $|G| < \theta$ and Rang$(H) \subseteq 
\theta$ and $\theta$ is regular); also this union belongs to $N_{i(*)}$ 
(as $G,H \in N_{i(*)})$,  hence the union is an ordinal
$< \delta_{i(*)} < \delta_{j(*)}$.  So $H(h) < \delta_{j(*)}$. 

But (by the choice of $h$ for the first inequality, and definition of $H(h)$
for the second inequality) 
\mr
\item "{$\otimes_1$}"   $f_{\delta_{j(*)}} \restriction A \leq g_{j(*)}
\restriction A < h \restriction A \leq f_{H(h)} \restriction A$
\ermn
and (by the definition of $H$ for the first equality, choice of $g^*$ and  
$j(*) \in S$  for the second): 
\mr
\item "{$\otimes_2$}"   $f_{H(h)} \restriction (\kappa \backslash A) = 
g^* \restriction (\kappa \backslash A) = f_{\delta_{j(*)}} \restriction 
(\kappa \backslash A)$.
\ermn
Now $\otimes_1,\otimes_2$ together implies $f_{\delta_{j(*)}} \leq f_{H(h)}$,
but as said earlier $H(h) < \delta_{j(*)} < \theta$, together they 
contradict the choice of $\langle f_\alpha:\alpha < \theta \rangle$.
\bn
\ub{Possibility $B$}:  The answer is no. 

So $H(h)$ is not well defined and without loss of generality $h \in N_{j(*)
+1}$ (as all parameters in the requirements on it are in  
$N_{j(*)+1})$.  Choose $j \in S,j > j(*)$;  as $H(h)$ is not well defined, no
$\alpha < \theta$ satisfies the requirements in $(*)$.
But of the three demands on $\alpha,\delta_j$ trivially satisfy two and 
a half:  $``\alpha < \theta,f_\alpha \restriction (\kappa \backslash A) = 
g^* \restriction (\kappa \backslash A)$ and $f_\alpha \restriction A \leq  
g^* \restriction A"$;  so the remaining one should fail, i.e.  
$\neg[h \restriction A \leq f_{\delta_j} \restriction A]$.  So for some  
$\zeta \in A$ 
we have $h(\zeta) > f_{\delta_j}(\zeta)$; now $h \in N_{j_{j(*)}
+1} \subseteq N_j$ hence $h(\zeta) \in N_j$ hence $h(\zeta) \in N_j \cap  
\sigma \backslash f_{\delta_j}(\zeta)$,  hence (by the definition of $g_j$),
$g_j(\zeta) \leq h(\zeta)$ hence (as $j \in S$) we have $g^*(\zeta) \leq  
h(\zeta)$ but $h \in G \subseteq \dsize \prod_{\xi \in A} g^*(\xi)$, so  
$h(\zeta) < g^*(\zeta)$,  contradiction.  \hfill$\square_{\scite{3.7}}$\margincite{3.7}
\enddemo
\newpage

\head {\S4  Entangled Orders --- Narrow Order Boolean Algebra Revisited}
\endhead  \resetall 
\bigskip

\proclaim{\stag{4.1} Theorem}  1) If 
$\kappa^{+4} \leq {\text{\rm cf\/}}(\lambda) < \lambda 
\leq 2^\kappa$ \ub{then} there is an entangled linear order 
of cardinality $\lambda$.  \nl
2)  Moreover, if $\chi_0 < \lambda$ we can demand that the linear order has
density character $\geq \chi_0$ (in fact, in every interval of the linear 
order).
\endproclaim
\bigskip

\remark{Remark}  See more \cite{Sh:462} and \cite{Sh:666}.
\endremark
\bigskip

\demo{Proof}  Without loss of generality $\chi_0 > \text{ cf}(\lambda)
\ge \kappa^{+4}$.  By \cite[2.1]{Sh:355} there 
is an increasing continuous sequence $\langle \lambda_i:
i < \text{ cf}(\lambda)\rangle$ of singular cardinals with 
limit $\lambda$ such that  
tcf$(\dsize \prod_{i < \text{cf}(\lambda)} \lambda^+_i,
<_{J^{\text{bd}}_{\text{cf}(\lambda)}}) = \lambda^+$ and 
$\lambda_0 > \chi_0$.
The proof will be split to cases (one of them relies on the solution to 
others for smaller cardinals, so you may want to say we are proving 
\scite{4.1} by induction on $\lambda$).  Without loss of generality 
$\chi_0 > \text{ cf}(\lambda)$.
\enddemo
\bn
\ub{Case I}:  For $i < \text{ cf}(\lambda)$ we have 
max pcf$\{\lambda^+_j:j < i\} < \lambda$. 

So for some unbounded $A \subseteq \text{ cf}(\lambda)$ we have $i \in A 
\Rightarrow \lambda^+_i > \text{ max pcf}\{\lambda^+_j:j \in A \cap i\}$. 

So ${\frak a} = \{\lambda^+_i:i \in A\}$ 
is as required in \cite[4.12]{Sh:355} (with  
$\lambda^+$, cf$(\lambda)$ here standing for $\lambda,\kappa$ there, noting 
that $2^{\text{cf}(\lambda)} \geq 2^\kappa \geq \lambda$).

So we can assume:
\bn
\ub{Assumption --- not Case I}

So there is $\mu,\chi_0 < \mu < \lambda$, cf$(\mu) < \text{ cf}(\lambda)$,  
pp$_{<\text{cf}(\lambda)}(\mu) > \lambda$.  Choose a minimal such $\mu$, so 
by \scite{3.2}(2): 
\mr
\item "{$(*)$}"  ${\frak a} \subseteq 
\text{ Reg } \backslash \chi_0 \and \sup({\frak a}) < 
\mu \and |{\frak a}| < \text{ cf}(\lambda) \Rightarrow
\text{ max pcf}({\frak a}) < \lambda$.
\ermn
Clearly (by \cite[2.3]{Sh:355}) in $(*)$'s conclusion we can replace  
$``< \lambda"$ by $``< \mu"$ i.e. 
\mr
\item "{$(*)'$}"  ${\frak a} \subseteq 
\text{ Reg } \backslash \chi_0 \and \sup({\frak a}) < 
\mu \and |{\frak a}| < \text{ cf}(\lambda) \Rightarrow 
\text{ max pcf}({\frak a}) < \mu$.
\ermn
Let $\sigma =: \text{ cf}(\mu)$, so pp$(\mu) = 
\text{ pp}_{<\text{cf}(\lambda)}
(\mu)$ 
(by \cite[1.6(3)]{Sh:371}) and remember pp$_{< \text{cf}(\lambda)}(\mu) 
> \lambda$.
\bn
\ub{Case II}:  $\sigma \geq \kappa$ 
(and not Case I; actually $2^\sigma \geq \mu$ suffices). 

First assume $\sigma > \aleph_0$.  As said above pp$_{<\text{cf}(\lambda)}
(\mu) > \lambda$ and by \cite[1.7]{Sh:371} there is a strictly increasing 
sequence  $\langle \mu^*_i:i < \sigma \rangle$ of regular cardinals
satisfying $\mu 
= \dbcu_{i < \sigma} \mu^*_i$, and $\lambda^+ = \text{ max pcf}\{\mu^*_i:
i < \sigma\} = \text{ tcf } \dsize \prod_{i < \sigma} \mu^*_i/
J^{\text{bd}}_\delta$.  Now as we can replace $\langle \mu^*_i:i < \sigma 
\rangle$ by $\langle \mu^*_i:i \in A \rangle$ for any $A \subseteq \sigma$  
unbounded, by $(*)'$ without loss of generality 
$\mu^*_i > \text{ max pcf}\{\mu^*_j:
j < i\}$, so we can apply again \cite[4.12]{Sh:355} 
(or \scite{3.2}(5)). 

When $\sigma = \aleph_0$, \scite{4.1} follows from \cite[4.13(1)]{Sh:355}.
\bn
\ub{Case III}:  cf$(\lambda) = \kappa^{+4}$ and $\sigma < \kappa$ (instead  
$\sigma < \mu,2^\sigma < \lambda$ suffice). 

So $\sigma^{+4} \leq \text{ cf}(\lambda)$.  Let ${\Cal P} =: \{A:A
\subseteq \text{ cf}(\lambda)$, otp$(A) = \kappa^+,A$ is a closed
subset of $\sup(A)$ and max pcf$\{\lambda^+_i:i \in c\} =
\lambda^+_{\sup(A)}\}$.  For any $C \in {\Cal P}$, try
to choose by induction on $i < \kappa^+({\frak b}_i) = {\frak b}_i[C]$
and $\gamma_i = \gamma_i[C]$ such that:
\mr
\widestnumber\item{$(iii)$}
\item "{$(i)$}"  ${\frak b}_i \subseteq \text{ Reg } \cap \mu \backslash
\dbcu_{j < i} {\frak b}_j \backslash \chi_0$ 
\sn
\item "{$(ii)$}"   $\gamma_i \in C \backslash \dbcu_{j < i} (\gamma_j+1)$
\sn
\item "{$(iii)$}"  $\lambda^+_{\gamma_i} \in \text{ pcf}(b_i)$ 
\sn
\item "{$(iv)$}"  $|{\frak b}_i| \leq \kappa$ 
\sn
\item "{$(v)$}"  all members of $b_i$ are $(\chi_0,\kappa^+,
\aleph_0)$-inaccessible
\sn
\item "{$(vi)$}"  $\gamma_i$ is minimal under those requirements.
\endroster
\bn
\ub{Subcase IIIa}:  For some $j < \kappa^+$.

For every $C \in {\Cal P}$ such that Min$(C) \ge j$.  
For some $\varepsilon(*) < \text{ cf}(\lambda)$, we cannot define  
${\frak b}_i,\gamma_i$ are defined iff $\varepsilon(*)$. 

Let $C,i(*)$ be as above.
Let $\gamma^* = \dbcu_{i< i(*)} \gamma_i$, so $\gamma^* \in C$.  
Now if $\gamma \in C \backslash \gamma^*$ then 
(by \cite[1.5B]{Sh:355}) as pp$_\sigma(\mu) 
\geq \lambda^+ > \lambda^+_\gamma$,
there is ${\frak a}_\gamma \subseteq 
\text{ Reg } \cap (\chi_0,\mu),|{\frak a}_\gamma| \leq  
\sigma$ such that $\lambda^+_\gamma \in \text{ pcf}({\frak
a}_\gamma)$.  
By \scite{3.2}(3) there is ${\frak c}_\gamma 
\subseteq \text{ Reg } \cap (\chi_0,\mu)$  
of cardinality $\leq \kappa$ consisting of $(\chi_0,\kappa^+,
\aleph_0)$-inaccessible cardinals such that 
$\lambda^+_\gamma \in \text{ pcf}({\frak c}_\gamma)$.  Now $\gamma,c_\gamma 
\backslash \dbcu_{i < i(*)} {\frak b}_i$ 
cannot serve as $\gamma_{i(*)},{\frak b}_{i(*)}$ so 
necessarily $\lambda^+_\gamma \notin \text{ pcf}({\frak c}_\gamma \backslash 
\dbcu_{i < i(*)} {\frak b}_i)$ 
hence without loss of generality ${\frak c}_\gamma \subseteq 
\dbcu_{i < i(*)} {\frak b}_i$.  
\sn
\ub{Version 2}:  So $\{\lambda^+_\gamma:\gamma \in C \backslash i(*)\}
\subseteq \text{ pcf}\{ \dbcu_{i < i(*)} {\frak b}_i\}$ and $|\dbcu_{i
< i(*)} {\frak b}_i| \le \kappa$.  By the proof of \cite[4.2]{Sh:400}
we get a contradiction.
\bn
\ub{Subcase IIIb}:  For every $j < \text{ cf}(\lambda)$ there are $C
\in {\Cal P}$ with Min$(C) > j$ such that for $C$, the pair
$({\frak b}_i,\gamma_i)$ defined for every $i < \kappa^+$.

We shall now show 
\mr
\item "{$\otimes$}"   for every $i(*) < \text{ cf}(\lambda)$ 
there is $\lambda' \in \lambda 
\cap \text{ pcf}\{\lambda^+_j:j < \text{ cf}(\lambda)\} \backslash 
\lambda_{i(*)}$ such that Ens$(\lambda',\lambda')$ (exemplified by linear 
order which has density character $> \chi_0$ in every interval). 
\ermn
\ub{Why $\otimes$ is sufficient}:  We can for 
$i < \text{ cf}(\lambda)$, choose
$\mu^*_i,\lambda_i < \mu^*_i = \text{ cf}(\mu^*_i) \in \lambda \cap
\text{ pcf}\{\lambda^+_j:j < \text{ cf}(\lambda)\}$, as 
required in $\otimes$.  As  
$\dsize \bigwedge_i \mu^*_i < \lambda$ without loss of generality 
$\langle \mu^*_i:i < \text{ cf}(\lambda)\rangle$ is (strictly) increasing.  
We try to choose by $i$ induction on $\epsilon < \text{ cf}(\lambda),
i(\epsilon) < \text{ cf}(\lambda)$ strictly increasing such that $\mu^*_{i(\epsilon)} > 
\text{ max pcf}\{\mu^*_{i(\zeta)}:\zeta < \epsilon\}$.

Let $i(\epsilon)$ be defined iff $\epsilon < \epsilon(*)$.  So  
$\epsilon(*)$ is $\le \text{ cf}(\lambda)$ and is a limit ordinal
$\lambda$ and $\text{max pcf}\{\mu^*_{i(\epsilon)}:
\epsilon < \epsilon(*)\} \ge \lambda$ hence $\ge \lambda^+$, but
pcf$\{\mu^*_{i(\varepsilon)}:\varepsilon < \varepsilon(*)\} \subseteq
\text{ pcf}\{\lambda^+_j:j < \text{ cf}(\lambda^+)\} \subseteq
\lambda^+ +1$ hence $\lambda^+ = \text{ max
pcf}\{\mu^*_{i(\varepsilon)}:\varepsilon < \varepsilon(*)\}$.
Note that $\mu^*_{i(\epsilon)} > \text{ max pcf}\{\mu^*
_{i(\zeta)}:\zeta < \epsilon\},\mu^*_{i(\epsilon)}$ is strictly increasing, 
and Ens$(\mu^*_{i(\epsilon)},\mu^*_{i(\epsilon)})$ for $\varepsilon < \varepsilon(*)$.  So applying \cite[4.12]{Sh:355} we finish. 
\bn
\ub{Why $\otimes$ holds}:  Let $i(*) < \text{ cf}(\lambda)$ be given.
Choose $C \subseteq (i(*)$, cf$(\lambda))$ from ${\Cal P}$  
such that $({\frak b}_\varepsilon[C],\gamma_\varepsilon[C])$ is
defined iff $\varepsilon < \varepsilon[C]$ and $\varepsilon[C] <
\kappa^+$. By the definition of ${\Cal P}$ we have
max pcf$\{\lambda^+_\gamma:\gamma \in C\} < \lambda$.
Let ${\frak d} =: \{\lambda^+_\gamma:\gamma \in C\}$,  
let $\langle {\frak b}_\theta[{\frak d}]:\theta \in \text{ pcf}({\frak
d}) \rangle$ be as in \cite[2.6]{Sh:371}.  Let $\theta$ be 
minimal such that otp$({\frak b}_\theta[{\frak d}])$ is $\ge \kappa$.
We can find $B_\epsilon \subseteq C$ (for $\epsilon < \kappa$) 
such that $\{\lambda^+_\gamma:\gamma \in  
B_\epsilon \} \subseteq {\frak b}_\theta[{\frak d}]$,  
otp$|B_\epsilon| = \kappa$ and the $B_\epsilon$'s are pairwise disjoint.
Clearly max pcf$\{\lambda^+_\gamma:\gamma \in B_\epsilon\} = \theta$ as  
$\{\lambda^+_\gamma:\gamma \in B_\epsilon\}$ is $\subseteq {\frak
b}_\theta[{\frak d}]$,  but is not a subset of any finite union 
of ${\frak b}_{\theta'}[{\frak c}],\theta' < \theta$.
Now letting ${\frak a}^* =: \cup\{{\frak b}_\varepsilon[C]:\varepsilon
< \varepsilon[C]\}$, there is (by \cite[2.6]{Sh:371}) a 
subset ${\frak a}$ of ${\frak a}^*$ such that 
$\theta = \text{ max pcf}({\frak a})$ but $\theta 
\notin \text{ pcf}({\frak a}^* \backslash {\frak a})$.  
Now as $\theta \in \text{ pcf}
\{\lambda^+_\gamma:\gamma \in B_\epsilon\},
\lambda^+_\gamma \in \text{ pcf}({\frak b}_\gamma)$ 
we have (by \cite[1.12]{Sh:345a})  $\theta \in \text{ pcf}
(\dbcu_{\gamma \in B_\epsilon} {\frak b}_\gamma)$ 
hence by the previous sentence  
$\theta \in  \text{ pcf}({\frak a} \cap 
\dbcu_{\gamma \in B_\epsilon} {\frak b}_\gamma)$.  
Let ${\frak c}_\epsilon =: {\frak a} \cap \dbcu_{\gamma \in
B_\epsilon} 
{\frak b}_\gamma,\lambda' = \theta$,  
we can apply \cite[4.12]{Sh:355} and get that there is an entangled linear 
order of cardinality $\lambda'$  (which is more than required, see 
\cite{Sh:345b}); and, of course, 
$\lambda_{i(*)} < \lambda' \in \lambda \cap  
\text{ pcf}\{\lambda_j:j < \text{ cf}(\lambda)\}$.  The 
assumptions of \cite[4.12]{Sh:355} holds as the ${\frak c}_\epsilon$ 
are pairwise disjoint (by (i) above), 
$\theta \in \text{ pcf}\{\lambda^+_\gamma:\gamma \in B_\epsilon 
\{\subseteq \text{ pcf}(\dbcu_{\gamma \in B_\epsilon} {\frak b}_\gamma) = 
\text{ pcf}({\frak c}_\epsilon)$ and 
$[\theta_1 \in {\frak a} \Rightarrow \text{ max pcf}
({\frak a} \cap \theta_1) < \theta_1]$ as $\theta_1$ is $(\chi_0,\kappa^+,
\aleph_0)$-inaccessible and $\theta = \lambda' \geq \sup 
\{\lambda^+_\gamma:\gamma \in C\} > \lambda_{i(*)} > \chi_0$.  So $\otimes$ 
holds and we finish Subcase IIIb hence Case III.
\bn
\ub{Case IV}:  cf$(\lambda) > 
\kappa^{+4}$ and $\sigma \leq \kappa$ (and not Case I). 

For each $\delta < \text{ cf}(\lambda)$ of cofinality $\kappa^{+4}$ we can 
apply the previous cases (or the induction hypothesis on $\lambda$) and 
get an entangled linear order of power $\lambda^+_\delta$.  So $\otimes$ 
holds and we finish as in Subcase IIIb.  
\hfill$\square_{\scite{4.1}}$\margincite{4.1}
\bigskip

\proclaim{\stag{4.2} Claim}  Assume $\kappa^{+4} \leq \theta = \text{ cf}
(\theta),\langle \lambda_i:i < \theta \rangle$ is a strictly increasing 
sequence of regular cardinals, $\theta < \lambda_i \leq 2^\kappa$ and  
$\lambda_\theta = \text{ tcf}(\dsize \prod_{i < \theta} \lambda_i/
J^{\text{bd}}_\theta)$. \nl
1) If $\underset {i < \theta} {}\to \sup \lambda_i \leqq 2^\kappa$ 
\ub{then} there 
is an entangled linear order of cardinality $\lambda_\theta$. \nl
2) {\rm Ens} $(\lambda_\theta,2^\kappa)$.
\endproclaim
\bigskip

\remark{Remark}  Remember that if there is an entangled linear order 
in $\lambda$ then Ens$(\lambda,\lambda)$ (so \cite[7(5)]{Sh:345b}).
\endremark
\bigskip

\demo{Proof}  Same proof as \scite{4.1}.
\enddemo
\bigskip

\proclaim{\stag{4.3} Claim}  Assume 
\mr
\widestnumber\item{$(iii)$}
\item "{$(i)$}"   $\lambda$ is regular, uncountable 
\sn
\item "{$(ii)$}"   $\kappa < \lambda \Rightarrow 2^\kappa < 2^\lambda$ 
\sn
\item "{$(iii)$}"   for some regular $\chi \leqq 2^\lambda$ there is no 
linear order of cardinality $\lambda$ with $\geqq \chi$ Dedekind cuts or 
even no tree of cardinality $\lambda \geqq \chi\,\,\lambda$-branches.
\ermn
\ub{Then} $(2^{<\lambda} < 2^\lambda$ and) for some $\mu$
\mr
\item "{$(a)$}"   for every regular $\chi$ in
$(2^{<\lambda},2^\lambda]$ (or even $(\mu,2^\lambda])$  
there is an entangled linear order of cardinality $\chi$ and density
$\mu$,
\sn
\item "{$(b)$}"  $\mu \in (\lambda,2^{<\lambda}]$, cf$(\mu) = \lambda$, 
pp$_{\Gamma(\lambda)}(\mu) = 2^\lambda,\mu$ is $(\lambda,
\lambda^+,2)$-inaccessible 
\ermn
(the linear order is $(T,<_{\ell x}), 
T \subseteq {}^{\mu >}2$ has $\leqq \mu$ nodes and $\geqq \chi \,
\lambda$-branches).
\sn
[Saharon: see also \cite[\S3]{Sh:430}.
\endproclaim
\bigskip

\demo{Proof}  Note: 
$2^{<\lambda} < 2^\lambda$ [if $(\exists \theta < \lambda)
(2^\theta = 2^{<\lambda})$  by (ii), otherwise cf$(2^{<\lambda}) =
\lambda$ and by classical cardinal arithmetic, cf$(2^\lambda) > \lambda$,  
hence $2^{<\lambda} < 2^\lambda$].  By \cite[Lemma 5.11]{Sh:355} if the 
conclusion fails then for every regular $\chi$ in $(2^{<\lambda},
2^\lambda]$ there is $\mu,\lambda = \text{ cf}(\mu) < 
\mu \leqq 2^{<\lambda}$,  
pp$_{\Gamma(\mu)}(\mu) \geqq \chi$.  Choose a minimal $\mu$ such that  
$\lambda < \mu \leqq 2^{<\lambda}$, cf$(\mu) \leqq \text{ cf}(\lambda)$ and  
pp$^+(\mu) > \chi$ (note:  $\mu$ does not depend on $\chi$, by 
\cite[2.3]{Sh:355}).  So necessarily $\mu$ is $(\lambda,
\lambda^+,2)$-inaccessible.  Let $\chi \in (\mu,2^\lambda]$ be regular.  
As $(2^{<\lambda})^{<\lambda} = 2^{<\lambda}$ necessarily 
cf$(\mu) = \lambda$, so by 
\cite[1.6(3)]{Sh:371} there is a strictly increasing sequence  
$\langle \mu_i:i < \text{ cf}(\mu) = \lambda \rangle$ of regular cardinals,  
$\lambda < \mu_i < \mu,\mu = \dsize \sum_i \mu_i$ and $\chi = 
\text{ tcf}(\Pi \mu_i/J^{\text{bd}}_{\text{cf}\lambda})$.  As $\mu$ is 
$(\lambda,\lambda^+,2)$-inaccessible without loss of generality $\mu_i > 
\text{ max pcf}\{\mu_j:j < i\}$.  So by \cite[4.12]{Sh:355} we finish.      
\hfill$\square_{\scite{4.3}}$\margincite{4.3}
\enddemo
\bigskip

\demo{\stag{4.4} Conclusion}  1) For a class of cardinals $\mu$, there is 
an entangled linear order of cardinality $\mu^+$. \nl
2)  Assume $\lambda$ is strong limit singular.  \ub{Then} for some successor 
cardinal in $(\lambda,2^\lambda]$ there is an entangled linear order.
\enddemo
\bigskip

\demo{Proof}  1) By part (2). \nl
2)  If $\aleph_{\lambda^{+4}} < 2^\lambda$ then apply theorem \scite{4.1} 
(with $\lambda,\aleph_{\lambda^{+4}}$ here standing for $\kappa,\lambda$
there) so there is an entangled linear order of cardinality  
$\aleph_{\lambda^{+4}+1} (\leqq 2^\lambda)$,  which is as required.  So 
assume $2^\lambda \leqq \aleph_{\lambda^{+4}}$.  We know that there is a 
linear order of cardinality $2^\lambda$ and density character $\lambda$;  
hence (see \cite[AP,\S1]{Sh:g}) 
there is an entangled linear order of cardinality  
cf$(2^\lambda)$.  But as $2^\lambda \leqq \aleph_{\lambda^{+4}}$ 
necessarily cf$(2^\lambda)$ is a successor cardinal.  
\enddemo
\newpage

\head {\S5  prd: Measuring $\Pi f(i)$ by a family $\Gamma$ of ideals     
and family sequences $\langle B_i:i < \kappa \rangle,|B_i| < \mu_i$}
\endhead  \resetall 
\bn
In \cite[\S4]{Sh:371}, and here in \S1 we have dealt with 
generalizations of the measuring $\dsize \prod_{i < \kappa} f(i)/I$, i.e. 
whereas defining cov$(\lambda,\lambda,\theta,\sigma)$ we cover 
a set $a \in [\lambda]^{< \theta}$ by $< \sigma$ subsets of cardinality  
$< \lambda$; there we ask that $\kappa$ belongs to the closure to a normal 
ideal of $J$ union the family of $A \subseteq \kappa$ for which we succeed to 
cover.  Here we replace ``normal" by an abstract property $\Gamma$ (and 
phrase the required properties).  We also generalize normality to ideals on  
${\Cal Y}$ with $\iota:{\Cal Y} \rightarrow \kappa$, a generalization
used in \cite[\S4]{Sh:420}, \cite{Sh:430}.
\bigskip

\demo{\stag{5.0} Context}  1)  $\kappa$ is a regular uncountable cardinal,  
${\Cal Y}$ a set, $\iota$ a function from ${\Cal Y}$ onto $\kappa,  
{\Cal Y}_i = \iota^{-1}(\{i\})$.  Here $I,J$ vary 
on ideals on ${\Cal Y},\Gamma$ a 
family of proper ideals on ${\Cal Y}$.
\enddemo
\bigskip

\definition{\stag{5.1} Definition}  1)  $\Gamma_{{\Cal Y},\kappa,\sigma} = 
\{J:J$ a $\sigma$-complete ideal on ${\Cal Y}\}$ (if $({\Cal Y},\iota) = 
(\kappa,\text{id}_\kappa)$ this $\Gamma_{\kappa,\sigma}$ is essentially  
$\Gamma(\kappa^+,\sigma))$. \nl
2)  $\Gamma_{{\Cal Y},\kappa} = \Gamma^{\text{nor}}_{{\Cal Y},\kappa} =: 
\{J:J$ a normal ideal on ${\Cal Y}\}$ (normal --- see \scite{5.2}(0) below).
If ${\Cal Y} = \kappa,\iota = \text{ id}$ we write 
$\Gamma^{\text{nor}}_\kappa$.
\enddefinition
\bigskip

\definition{\stag{5.2} Definition}  0)  An ideal $I$ on ${\Cal Y}$ is 
normal \ub{if}: for any club $C$ of $\kappa$ the set $\dbcu_{i \notin C}
{\Cal Y}_i$ belong to $I$ 
and for any sequence $\langle A_i:i < \kappa \rangle$  
of sets from $I$ the set $\triangledown_i A_i =: \{x \in {\Cal Y}:x \in  
\dbcu_{j < \iota(x)} A_j\}$ belongs to $I$.  (So normal implies 
$\kappa$-complete). \nl
1) We say $\Gamma$ is $\sigma$-complete \ub{if} every $J \in \Gamma$ is 
$\sigma$-complete. \nl
2)  We say $\Gamma$ is normal \ub{if} every $J \in \Gamma$ is normal.
\nl
3)  We say $\Gamma$ is restriction closed \ub{when}:  
$J \in \Gamma,A \subseteq  
\kappa,A \ne \emptyset \text{ mod } J$ implies there is $I \in \Gamma,
J \cup \{\kappa \backslash A\} \subseteq I$. \nl
4)  We say $\Gamma$ is closed \ub{if} for every ${\Cal P} \subseteq {\Cal P}
({\Cal Y}),c \ell_\Gamma({\Cal P})$ is well defined where 
$c \ell_\Gamma({\Cal P})$ is the minimal member of $\Gamma \cup 
\{{\Cal P}(\kappa)\}$ which include it, i.e. $(\forall I \in \Gamma)
[{\Cal P} \subseteq I \Leftrightarrow  c \ell_\Gamma({\Cal P}) \subseteq I]$.

Note:  $c \ell_\Gamma$ for $a$ not 
necessarily closed $\Gamma$, is a partial function. \nl
5)  We say $\Gamma$ has character $\leqq \mu$ \ub{when}: ${\Cal P}(\kappa) =
c \ell_\Gamma({\Cal P})$ where ${\Cal P} \subseteq {\Cal P}(\kappa)$ 
implies that for some ${\Cal P}' \subseteq {\Cal P}$ of cardinality  
$\leqq \mu$, we have ${\Cal P}(\kappa) = c \ell_\Gamma({\Cal P})$. \nl
6)  The character of $\Gamma$ is the minimal cardinal $\mu$ such that  
$\Gamma$ has character $\leqq \mu$.
\enddefinition
\bigskip

\definition{\stag{5.3} Definition}  1)  We say $\Gamma$ is suitable if it 
is $\theta$-suitable for every $\theta$?  We say that $\Gamma$ 
is $\theta$-suitable when: for every ideal 
$J \in \Gamma$ on ${\Cal Y}$, if
\ub{then} ${\Cal P}(\kappa) = c \ell_\Gamma(J \cup \{f^{-1}_\eta(\{0\}):
\eta \in T\})$ where
\mr
\widestnumber\item{$(iii)$}
\item "{$(i)$}"   $T$ is a (non-empty) set of finite sequences of 
ordinals $< \theta$ closed under initial segments 
\sn
\item "{$(ii)$}"   $A_\eta \subseteq {\Cal Y}$ for $\eta \in T$ and
$A_{<>} = {\Cal Y}$ 
\sn
\item "{$(iii)$}"   for each $\eta \in T$ of length $n$,
$$
{\Cal P}(\kappa) = 
c \ell_\Gamma (J \cup \{\kappa \backslash A_{\eta \restriction 0},\dotsc,
\kappa \backslash A_{\eta \restriction n}\} \cup \{A_{\eta \char94 <i>}:
\eta \char94 <i> \in T\})
$$ 
\item "{$(iv)$}"  $\eta \triangleleft \nu \Rightarrow A_\eta \supseteq
A_\nu$ 
\item "{$(v)$}"   $f_\eta:A_\eta \rightarrow$ Ord 
\sn
\item "{$(vi)$}"  if $\eta \char94 < \zeta > \in T$ and $y \in  
A_{\eta \char94 <\zeta >},f_\eta (y) \ne 0$ then $f_\eta(y) > 
f_{\eta \char94 <\zeta >}(y)$.
\ermn
\enddefinition
\bigskip

\remark{\stag{5.3A} Remark}  1) Clearly for 
$\theta_1 \leqq \theta_2$ if $\Gamma$ is $\theta_2$-suitable then 
$\Gamma$ is $\theta_1$-suitable. \nl
2) If $\theta \geqq 2^{|{\Cal Y}|}$ its value 
is immaterial, so we can omit it.
\endremark
\bigskip

\proclaim{\stag{5.4} Claim}  1) $\Gamma_{{\Cal Y},\kappa,\sigma}$, if 
$\kappa \geqq \sigma = {\text{\rm cf\/}}(\sigma) > \aleph_0$ is closed, 
restriction closed, of character $\sigma$ and is suitable.  For ${\Cal P} 
\subseteq {\Cal P}({\Cal Y})$,

$$
c \ell_{\Gamma_{{\Cal Y},\kappa,\sigma}}({\Cal P}) = 
\bigl\{ Z:\text{ for some } \alpha < \sigma,A_i \in {\Cal P} \text{ for }  
i < \alpha \text{ we have } Z \subseteq \dbcu_{i < \alpha} A_i \bigr\}.
$$
\mn
2) If ${\Cal P}$ is a family of subsets of ${\Cal Y},\Gamma = 
\Gamma^{\text{nor}}_{{\Cal Y},\kappa}$ \ub{then} 

$$
\align
c \ell_\Gamma({\Cal P}) = \biggl\{Z:&Z \subseteq {\Cal Y} \text{ and for some 
club } C \text{ of } \kappa \text{ and sequence} \\      
  &\langle A_i:i < \kappa \rangle \text{ of member of } {\Cal P}
\text{ we have}: \\
  &Z \subseteq \{x \in {\Cal Y}:\iota(x) \in C \text{ and } x \in 
\dbcu_{j < \iota(x)} A_j \biggr\}.
\endalign
$$
\mn
3)  $\Gamma^{\text{nor}}_{{\Cal Y},\kappa}$ (remember $\kappa = 
\text{ cf}(\kappa) > \aleph_0)$ is 
closed, restriction closed, suitable and of character $\kappa$. \nl
4)  If $\Gamma$ is $\theta$-suitable and has character $\leqq \theta$ \ub{then}
it is suitable (we shall use this freely).
\endproclaim
\bigskip

\demo{Proof}  1) Let us check suitability leaving the rest to the 
reader; so let $J,\theta,\langle A_\eta,f_\eta:\eta \in T \rangle$ be as in 
Definition \scite{5.3}.  If the required conclusion in Definition \scite{5.3}
fails then there is a $\sigma$-complete filter $I$ on $\kappa$ containing  
$J \cup \{f^{-1}_\eta(\{0\}):\eta \in T\}$.  For each $\eta$, by condition 
(iii), for some set $w_\eta$ of ordinals, $|w_\eta| < \sigma$ and  
$B_\eta \in J$ we have 
\mr
\item "{$(*)$}"   ${\Cal Y} = B_\eta \cup \{\kappa \backslash A_{\eta 
\restriction 0} \backslash \ldots \backslash A_{\eta \restriction \ell g
\eta}\} \cup \{A_{\eta \char94 <\zeta >}:\zeta \in w_\eta\}$ \nl
$\qquad$ (and $\zeta \in w_\eta \Rightarrow  \eta \char94 <\zeta > \in T)$.
\ermn
Let $T^* =: \{ \eta \in T:\ell < \ell g(\eta) \Rightarrow \eta(\ell) \in  
w_{\eta \restriction \ell}\}$, as $\sigma$ is regular uncountable, clearly  
$|T^*| < \sigma$.  Let

$$
B =: \cup \{B_\eta:\eta \in T^*\} \cup \{A_\eta:\eta \in T^* 
\text{ and } A_\eta \in I\} \cup \{f^{-1}_\eta(\{0\}):\eta \in T^*\}.
$$ 

Now $B$ is the union of $< \sigma$ members of $I$ and $I$ is 
$\sigma$-complete, so $B \in I$.  Choose  $y \in \kappa \backslash B$.

We now choose by induction on $n,\eta_n$ such that:

$$
\eta_0 = <>
$$

$$
\eta_n \in T^*,
$$

$$
\ell g(\eta_n) = n
$$

$$
\eta_n \trianglelefteq \eta_{n+1}
$$

$$
y \in A_{\eta_n}.
$$
\mn
For $n = 0,\eta_0 = <>$ so $y \in {\Cal Y} = A_{<>}$.  
For $n + 1,y \in A_{\eta_n \restriction n},\eta_n \in T^*$ and 
$y \notin B$ hence $y \notin B_{\eta_n}$, so by 
$(*)$ there is $\zeta_n \in w_{\eta_n}$ such that 
$y \in A_{\eta_n \char94 <\zeta_n>}$ so 
$\eta_{n+1} = \eta_n \char94 \langle \zeta \rangle$ is as required. 

In the end for 
each $n$ we have $f_{\eta_n}(y) > 0$ as $f^{-1}_{\eta_n}(\{0\})
\subseteq B$, (remember $\eta_n \in T^*$) hence by condition (vi) from 
\scite{5.3}, $\langle f_{\eta_n}(y):n < \omega \rangle$ is a strictly 
decreasing sequence of ordinals, contradiction. \nl
2)  Left to the reader. \nl
3)  Again we leave the proof of restriction closed and closed and having 
character $\kappa$ to the reader and prove suitability.  The proof is 
similar but use diagonal union.  So by part (2) and condition (iii) of 
\scite{5.3} for each $\eta \in T$ for some $B_\eta \in J$ and function  
$h_\eta$ from $\kappa$ to ordinals such that $\eta \char94 \langle 
h_\eta(i) \rangle \in T$ for $i < \kappa$,  we have 
\mr
\item "{$(*)$}"   ${\Cal Y} = B_\eta \cup \{\kappa \backslash 
A_{\eta \restriction 0} \backslash A_{\eta \restriction 1} \backslash \ldots
\backslash A_{\eta \restriction \ell g \eta}\}\cup \{x \in {\Cal Y}:
x \in \dbcu_{i< \iota(x)} A_{\eta \char94 \langle h_\eta (i)\rangle}\}$.
\ermn
By renaming, without loss of generality $T \subseteq {}^{\omega >}\kappa$ and
each $h_\eta$ the identity function.

Let $I$ be the normal ideal on ${\Cal Y}$ generated by $J \cup 
\{f^{-1}_n(\{0\}):\eta \in T\}$, so we assume $I$ is a normal proper ideal 
and we shall get a contradiction.

Now define ${\Cal Y}^*$  

$$
\align
{\Cal Y}^* = \{x \in {\Cal Y}:\,&(a) \quad \iota(x) \text{ a limit ordinal }
 < \kappa \\        
  &(b) \quad \text{ if } \eta \in T \cap {}^{\omega >} \iota(x) \text{ then }
x \notin B_\eta \\
  &(c) \quad \text{ if } \eta \in T \cap {}^{\omega >} \iota(x) \text{ and }
A_\eta \in I \text{ then } x \notin A_\eta\}.
\endalign
$$ 
\mn
Clearly ${\Cal Y}^*  \equiv {\Cal Y}$ mod $I$,  hence we can find  
$x(*) \in {\Cal Y}^*$.  Now we choose by induction on $n \eta_n \in
{}^{\omega >} \iota(x(*))$ as in the proof of part (1) and get similar 
contradiction. \nl
4)  Left to the reader.     \hfill$\square_{\scite{5.4}}$\margincite{5.4}
\enddemo
\bigskip

\definition{\stag{5.5} Definition}  1) For $\bar \mu = \langle \mu_i:
i < \kappa \rangle,J$ an ideal on ${\Cal Y}$ and $f:{\Cal Y} \rightarrow$ ord
we define 

$$
\align
\text{prd}^\Gamma_J(f,\bar \mu) = \text{ Min} \bigl\{ |{\Cal P}|:&{\Cal P}
\text{ is a family of sequences of the form } \langle B_x:x \in {\Cal Y}
\rangle, \\
  &\text{ each } B_x \text{ a set of ordinals}, B_x \text{ of cardinality} \\
  &< \mu_{\iota(x)} \text{ such that for every } g \in 
{}^x \text{ord satisfying} \\
  &g \leqq f, \text{ we have}: \\
  &{\Cal P}(\kappa) = c \ell_\Gamma[J \cup \{\{x \in {\Cal Y}:g(x) \in B_x\}:
\langle B_x:x \in {\Cal Y} \rangle \in {\Cal P}\}] \bigr\}.
\endalign
$$
\mn
2)  If above $J \subseteq I,I$ an ideal on $\kappa$   

$$
\align
\text{prd}^\Gamma_{J,I}(f,\bar \mu) = \text{ Min} \bigl\{ |{\Cal P}|:&{\Cal P}
\text{ is a family of sequences } \langle B_x:x \in {\Cal Y} \rangle, \\
  &\text{ each } B_x \text{ a set of ordinals of cardinality } < 
\mu_{\iota(x)} \\ 
  &\text{ and for every } g \in {}^x \text{ord satisfying } q \leqq f, \\
  &\text{ we have}: \\
  &I \subseteq c \ell_\Gamma(J \cup \{\{x \in {\Cal Y}:g(x) \in B_x\}:
\langle B_x:x \in X \rangle \in {\Cal P}\}) \bigr\}.
\endalign
$$
\mn
3) If $\bar \mu$ is constantly $\mu$, we may write $\mu$ instead.  We
can use also $\bar \mu = \langle \mu_x:x \in {\Cal Y}\rangle$, (but usually 
do not) with the obvious meaning.
\enddefinition
\bigskip

\proclaim{\stag{5.6} Claim}  1) If $\{x \in {\Cal Y}:f(x) < \omega\} \in J$  
we can in \scite{5.5} demand $g <_J f$, and if in addition 
$\dsize \bigwedge_x f(x) \ne 0$ we can demand $g < f$ (without changing 
the values). \nl
2)  If $\{x:f(x) \geqq \mu_{\iota(x)}\} \in J$ \ub{then}
{\rm prd}$^\Gamma_J(f,\bar \mu) = 1$. \nl
3)  If $\dsize \bigwedge \mu_i = \mu, {\text{\rm cf\/}}(\mu_ > 
|{\Cal Y}|,\{x:f(x) > \mu\} \in J$ \ub{then}
{\rm prd}$^\Gamma_J(f,\bar \mu) \leqq {\text{\rm cf\/}}(\mu)$. \nl
4)  If $f_1 \leqq_J f_2$ or just $\{x:|f_1(x)| > |f_2(x)|\} \in J$ 
\ub{then}  
{\rm prd}$^\Gamma_J(f_1,\bar \mu) \leqq {\text{\rm prd\/}}^\Gamma_J(f_2,
\bar \mu)$  (and the other obvious monotonicity properties). \nl
5)  If $\mu_i = \mu$, {\rm cf}$(\mu) > |{\Cal Y}|$  we can in 
Definition \scite{5.5} demand  $\dsize \bigwedge_{x \in {\Cal Y}} B_x = B_0$ 
for $\bar B \in {\Cal P}$ (i.e. without changing the value). \nl
6)  If $\Gamma$ is $\sigma$-complete and restriction closed, $\epsilon(*) 
< \sigma,\langle A_\epsilon:\epsilon < \epsilon(*) \rangle$ is a partition 
of ${\Cal Y} \, (J,\bar \mu$ as in \scite{5.5}) \ub{then}

$$
\underset\ \epsilon {}\to \sup \, \text{ prd}^\Gamma_{J+(\kappa \backslash 
A_\epsilon)}(f,\bar \mu) \leqq \, {\text{\rm prd\/}}^\Gamma_J
(f,\bar \mu) \leqq \dsize \sum_\epsilon \, 
{\text{\rm prd\/}}^\Gamma_{J+(\kappa \backslash A_\epsilon)}(f,\bar \mu).
$$
\mn
7)  If $\Gamma$ is normal and restriction closed, $A_\epsilon \subseteq  
\{x \in {\Cal Y}:\iota(x) > \epsilon\}$ for $\epsilon < \kappa,
\langle A_\epsilon:\epsilon < \kappa \rangle$ a partition of ${\Cal Y}
\backslash \iota^{-1}(\{0\})$ then

$$
\underset {\epsilon < \kappa} {}\to \sup \, {\text {\rm prd\/}}
^\Gamma_{J+({\Cal Y} \backslash A_\epsilon)} \leqq {\text{\rm prd\/}}
^\Gamma_J(f,\bar \mu) \leqq \dsize \sum_{\epsilon < \kappa}
{\text{\rm prd\/}}^\Gamma_{J+({\Cal Y} \backslash A_\epsilon)}
(f,\bar \mu).
$$
\mn
8)  {\rm prd}$^{\Gamma^{\text{nor}}_\kappa}_J(f,\bar \mu) = 
{\text{\rm prc\/}}_J(f,\bar \mu)$ if $J \in \Gamma^{\text{nor}}_\kappa$. \nl
9)  Assume $\Gamma$ is normal, $\bar \mu^\ell = \langle \mu^\ell_i:i < \kappa 
\rangle$ for $\ell = 1,2,\bar \mu^1$ increasing continuous and for each  
$i$, {\rm cf}$(\mu^1_i) < 
\mu^1_i \and \mu^2_1 = (\mu^1_i)^+$.  \ub{Then} for (any  
$\Gamma,J \in \Gamma$ and $f \in {}^{\Cal Y}$Ord) we have 
{\rm prd}$^\Gamma_J(f,\bar \mu^1) = {\text{\rm prd\/}}^\Gamma_J
(f,\bar \mu^2)$. \nl
10)  If $\bar \mu = \langle \mu_i:i < \kappa \rangle$ is increasing 
continuous with limit $\mu$, and $\Gamma$ is normal, $J \in \Gamma$ then  
{\rm prd}$^\Gamma_J(\bar \mu,\bar \mu) \leqq \mu$ and 
even {\rm prd}$^\Gamma_J(\bar \mu,\bar \mu) = {\text{\rm cf\/}}(\mu)$. \nl
11)  If $\dsize \bigwedge_{i < \kappa} \mu_i = \mu$, {\rm cf}$(\mu) > 
|{\Cal Y}|,\Gamma = \Gamma_{{\Cal Y},\kappa,\sigma}$ then for any $\alpha$,  
{\rm prd}$^\Gamma_{\{\emptyset\}}(\alpha,\bar \mu) = {\text{\rm cov\/}}
(|\alpha|,\mu,|{\Cal Y}|^+,\sigma)$.
\endproclaim
\bigskip

\demo{Proof}  E.g. \nl
6)  The first inequality should be clear.  Also the second: assume it fails, 
let  
$\lambda_\epsilon = \text{ prd}^\Gamma_{J+(\kappa \backslash A_\epsilon)}
(f,\bar \mu),\dsize \sum_{\epsilon < \epsilon(*)} \lambda_\epsilon < 
\text{ prd}^\Gamma_J(f,\bar \mu)$, let ${\Cal P}_\epsilon$ exemplify the 
definition of $\lambda_\epsilon$ and ${\Cal P}$ be 
$\dbcu_{\epsilon < \epsilon(\ast)}{\Cal P}_\epsilon$.  As 
$\dsize \sum_\epsilon \lambda_\epsilon < \text{ prd}^\Gamma_J(f,\bar \mu),
{\Cal P}$ cannot exemplify $\dsize \sum_\epsilon \lambda_\epsilon \geqq  
\text{ prd}^\Gamma_J(f,\bar \mu)$, so there is a function $g \in
{}^\kappa$ ord exemplifying this, so there is a proper ideal $I \in \Gamma$  
extending

$$
J \cup \bigl\{ \{x \in {\Cal Y}:g(x) \in B_x\}:\langle B_x:x \in {\Cal Y}
\rangle \in {\Cal P} \bigr\}.
$$
\mn
As $\Gamma$ is $\sigma$-complete also $I$ is $\sigma$-complete so for some
$\epsilon < \epsilon(*)$ we have $A_\epsilon \notin I$;  but  
$\Gamma$ is restriction closed so there is $I_1 \in \Gamma,I \cup \{{\Cal Y}
\backslash A_\epsilon \} \in I_1$.  So $I_1 \in \Gamma$ extend

$$
[J \cup \{{\Cal Y} \backslash A_\epsilon \}] \cup \{\{x \in {\Cal Y}:g(x) \in
B_x\}:\langle B_x:x \in {\Cal Y} \rangle \in {\Cal P}_\epsilon \} 
$$
\mn
contradicting the choice of ${\Cal P}_\epsilon$.  
\hfill$\square_{\scite{5.6}}$\margincite{5.6}
\enddemo
\bigskip

\proclaim{\stag{5.7} Lemma}  1) Suppose
\mr
\item "{$(*)$}"  $\mu_i = \mu = {\text{\rm cf\/}}(\mu) > |{\Cal
Y}|,\Gamma$ is suitable, restriction closed, $f \in {}^{\Cal Y}$Ord 
and $J$ an ideal on  
$\kappa$.
\ermn
\ub{Then}: 

$$
\align
\mu + {\text{\rm prd\/}}^\Gamma_J(f,\bar \mu) = \mu + \sup \bigl\{
{\text{\rm tcf\/}} \dsize \prod_{x \in {\Cal Y}} \lambda_x/I:&I
\text{ an ideal on } {\Cal Y} \text{ in } \Gamma \tag"${\otimes}$" \\       
  &\text{extending } J \text{ such that } \\
  &\mu < \lambda_x = \text{ cf}(\lambda_x) \leqq f(x) \bigr\}.
\endalign
$$
\mn
2)  If 
$\Gamma$ is normal, $\bar \mu = \langle \mu_i:i < \kappa \rangle,\mu_i 
= \theta^+_i,\langle \theta_i:i < \kappa \rangle$ is increasing continuous,  
$\theta_i > \kappa,\mu_i > |{\Cal Y}_i|,\mu = \dbcu_{i < \kappa}
\mu_i,\Gamma$ suitable and restriction closed and $f \in {}^{\Cal
Y}$Ord, $J$ an ideal on ${\Cal Y}$ then $(\otimes)$ above holds.
\endproclaim
\bigskip

\demo{Proof}  Like the proof of \scite{1.1}.
\enddemo
\bigskip

\demo{\stag{5.8} Conclusion}  1) Suppose $\mu,\bar \mu,\Gamma,J$ are as in 
\scite{5.7}, $f \in {}^{\Cal Y}$ ord and $\{x:\text{cf}[f(x)] = f(x) \geqq  
\mu_{\iota(x)}\} \in J$  \ub{then}

$$
\align
\mu + \text{ prd}^\Gamma_J(f,\bar \mu) = \mu &+ \sup \{\text{prd}^\Gamma_J
(g,\bar \mu):g <_J f\} \\
  &+ \sup \{ \text{tcf}_x \dsize \prod_{x \in {\Cal Y}} f(x)/I:J \subseteq I 
\in \Gamma \text{ (and the {\rm tcf} well defined)}\}.
\endalign
$$
\mn
2)  If in addition $G \subseteq \{g:g \in {}^{\Cal Y}$Ord, $g <_Jf\}$  
is cofinal or at least ${\Cal P}(\kappa) = c \ell_\Gamma [I \cup \{\{i:
h(i) < g(i)\}:g \in G\}]$  for every  $h <_J f$,  (and e.g. 
$\dsize \bigwedge_i f(i) \geqq \theta$) \ub{then}

$$
\align
\mu + \text{ prd}^\Gamma_J(f,\bar \mu) = \mu &+ \sup \{\text{prd}^\Gamma_J
(g,\bar \mu):g \in G\} \\
  &+ \sup \{ \text{tcf} \dsize \prod_{x \in {\Cal Y}} f(x)/I:J \subseteq I 
\in \Gamma\}.
\endalign
$$
\enddemo
\bigskip

\proclaim{\stag{5.9} Claim}  Suppose $\mu,\bar \mu,\Gamma,J$ are as in 
\scite{5.7}((1) or (2)), $g:{\Cal Y} \rightarrow$ card and 
$f(x) = g(x)^+ \geqq \mu^+_{\iota(x)}$.

\ub{Then} {\rm prd}$^\Gamma_J(f,\bar \mu) \leqq [{\text{\rm prd}}^\Gamma_J
(g,\bar \mu)]^+ + \mu$.
\endproclaim
\bigskip

\demo{Proof}  Let ${\Cal P}$ exemplify the value of prd$^\Gamma_J
(g,\bar \mu)$, say $|{\Cal P}| = \chi$.  So for every $h <_J f$, clearly  
$\{x \in {\Cal Y}:|h(x)| \leqq g(x)\} = {\Cal Y}$ mod $J$,  hence there is  
${\Cal P}_h \subseteq \dsize \prod_{x \in {\Cal Y}} h(x)$ exemplifying  
prd$^\Gamma_J(h,\bar \mu) \leqq \text{ prd}^\Gamma_J(g,\bar \mu) = \chi$.  
Assume prd$_J(f,\bar \mu) > \chi^+ + \mu$; by \scite{5.7} there is  
$f':{\Cal Y} \rightarrow$ ord, each $f'(x)$ a regular cardinal $\geqq \mu_i,
f' \leqq  f$ mod $I$ where 
$I \in \Gamma$ an ideal extending $J$ such that $\chi^+ 
< \text{ tcf } \dsize \prod_{x \in {\Cal Y}} f'(x)/I$.  Let $\langle h_\zeta:
\zeta < \chi'\rangle$ be $<_I$-increasing cofinal in $\dsize \prod_{x \in
{\Cal Y}} f'(x)/I$.  As in \cite[1.5]{Sh:355} without loss of generality 
for some $\zeta(*) < \text{ tcf } \dsize \prod_{x \in {\Cal Y}} f'(x)/I$  
of cofinality $\chi^+$ we have: $\langle h_\zeta:\xi < \zeta(*)\rangle$ has 
a $<_I$-lub $h'$ such that: for \scite{5.7}(1) $\{x \in {\Cal Y}:\text{cf}
[f(x)] \leqq \mu\} \in I$ and for \scite{5.7}(2) $\{x \in {\Cal Y}:\text{cf}
[f(x)] < \mu_{\iota(x)}\} \in I$;  without loss of generality it is 
$h_{\zeta(*)}$ and $\dsize \bigwedge_{x \in {\Cal Y}} \text{ cf}[h_{\zeta(*)}
(x)] \geqq \mu_{\iota(x)}$,  and without loss of generality:  $\xi < \zeta(*) 
\Rightarrow h_\xi < h_{\zeta(*)}$.  For each $\bar B = \langle B_x:x \in  
{\Cal Y} \rangle \in {\Cal P}_{h_{\zeta(*)}}$ define a function  
$f_{\bar B}:f_{\bar B}(x) = \sup(h_{\zeta(*)}(x) \cap B_x)$.  
So $f_{\bar B} < h_{\zeta(*)}$ hence for some $\xi(\bar B) < \zeta(*)$  
we have $f_{\bar B} < f_{\xi(\bar B)}$ mod $I$.  Let $\dbcu_{\bar B}
\xi(\bar B) < \xi < \zeta(*)$ --- possible as the number of $\bar B$'s is  
$\leqq |{\Cal P}_{h_{\zeta(*)}}| \leqq \chi < \chi^+ = \text{ cf}(\zeta(*))$.
So for $\bar B \in {\Cal P}_{h_{\zeta(*)}}$ we have
$\{x \in {\Cal Y}:f_\xi(x) \in B_x\} \in I$.  But $f_\xi < f_{\zeta(*)}$ 
so we get contradiction to the choice of ${\Cal P}_{h_{\zeta(*)}}$.
\hfill$\square_{\scite{5.9}}$\margincite{5.9}
\enddemo
\bigskip

\definition{\stag{5.10} Definition}  Let $1-\text{cf}^\Gamma_J(\langle a_x:
x \in {\Cal Y} \rangle)$ be  

$$
\sup \bigl\{\text{tcf } \dsize \prod_{x \in {\Cal Y}} \lambda_x/I:J 
\subseteq I \in \Gamma \text{ and } \lambda_x \in a_x \text{ for } x \in  
{\Cal Y} \bigr\}.
$$
\enddefinition
\bigskip

\proclaim{\stag{5.11} Claim}  $\bar \mu = \langle \mu_i:i < \kappa \rangle$ is 
non-decreasing, $\Gamma$ is a suitable restriction close family of ideals on
${\Cal Y},J \in \Gamma,f \in {}^{\Cal Y}$ ord and 
$\dsize \bigwedge_{x \in {\Cal Y}} f(x) \geqq \mu_{\iota(x)} = \text{ cf }
\mu_{\iota(x)}$. 
\sn
1)  If $\lambda \leqq \text{ prd}^\Gamma_J(f,\bar \mu)$ is regular, then 
for some $\langle {\frak a}_x:x \in {\Cal Y} \rangle$  we have:  
\mr
\widestnumber\item{$(iii)$}
\item "{$(i)$}"   ${\frak a}_x \subseteq \text{ Reg } \cap f(x)^+ \backslash 
\mu_{\iota(x)}$  
\sn
\item "{$(ii)$}"   $|\dbcu_{x \in {\Cal Y}} {\frak a}_x| < \mu$ if $(*)$ of 
\scite{5.7}(1) and $|{\frak a}_x|^+ < \mu_{\iota(x)}$ when \scite{5.7}(2)'s 
assumptions hold
\sn
\item "{$(iii)$}"  $\lambda = 1-\text{cf}^\Gamma_J(\langle a_x:x \in {\Cal Y}
\rangle)$.
\ermn
2)  If $\lambda$ is inaccessible, $\Gamma,\sigma$-complete and $[\chi < \mu  
\Rightarrow  \text{ cov}(\chi,\chi_0,\kappa,\sigma) < \lambda]$ \ub{then} 
without loss of generality $|\dbcu_{x \in {\Cal Y}} a_x| < \chi_0$.
\endproclaim
\bigskip

\demo{Proof}  1) Like the proof of 1.1.
\nl
2)  Straight.
\enddemo
\bigskip

\proclaim{\stag{5.12} Claim}  Assume the hypothesis of \scite{5.7}. 

If $g \in {}^{\Cal Y}$Ord, and each $g(s)$ is an ordinal $\geqq  
\mu_{\iota(x)}$ and $f(i) = \aleph_{g(i)}$ and let $\lambda = 
{\text{\rm prd\/}}^\Gamma_J(g,\bar \mu) + \dbcu_{\alpha < \mu} |\alpha|^{+3} 
+ |{\Cal Y}|$ \ub{then} 
{\rm prd}$^\Gamma_J(f,\bar \mu) \leqq \aleph_{\lambda^+}$.
\endproclaim
\bigskip

\demo{Proof}  Assume not, so prd$^\Gamma_J(f,\bar \mu) \geqq 
\aleph_{\lambda^+ +1}$ hence by \scite{5.7} there is $I \in \Gamma,J 
\subseteq I$ and $f^* \leqq f$ such that: each $f^*(x)$ is a regular cardinal
$\geqq \mu_{\iota(x)}$ and $\aleph_{\lambda^+ +1} \leqq \text{ tcf}
(\dsize \prod_{x \in {\Cal Y}} f^*(x)/I)$.  By \cite[1.5]{Sh:355} for each  
$\alpha < \lambda^+$ such that $\aleph_\alpha \geqq |{\Cal Y}|$ (e.g.  
$\alpha \geqq \lambda$) there is $f_\alpha <_I f^*,f_\alpha \leqq f^*$,  
each $f_\alpha(x)$ a regular cardinal $> \mu_{\iota(x)}$ and  
tcf$(\dsize \prod_{x \in {\Cal Y}} f^*(x)/I) = \aleph_{\alpha +1}$.  
Clearly $\alpha \ne \beta \Rightarrow f_\alpha \ne_I f_\beta$.  Now let  
$f_\alpha = \aleph_{g_\alpha}$, so $g_\alpha <_I g,g_\alpha \leqq g$.  Let  
${\Cal P}$ exemplify prd$^\Gamma_J(g,\bar \mu) \leqq \lambda$, so for each  
$\alpha < \lambda^+$ for some $\bar B^\alpha \in {\Cal P},\{y:f_\alpha(y) 
\in \{\aleph_j:j \in B^\alpha_y\}\} \notin I$ (and, by \scite{5.6}(9), and 
normality without loss of generality 
$\sup\{|B^\alpha_y|:y \in {\Cal Y}\} < \dsize \sum_{i<\kappa}
\mu_i)$.  By \scite{3.3} we get contradiction.       
\hfill$\square_{\scite{5.12}}$\margincite{5.12}
\enddemo
\newpage

\head {\S6  The Existence of Strongly Almost Disjoint Families} \endhead  \resetall \sectno=6
\bn
See \cite[\S0]{Sh:355} on the history of the subject.
\proclaim{\stag{6.1} Theorem}  Assume $J$ is an ideal on $\kappa,\kappa$ not 
the union of $\aleph_0$ members of $J,\mu > \kappa^{<\sigma}$ where 
\mr
\item "{$\otimes$}"   $\sigma = \kappa^+$ or at least $\forall A \in J^+
\exists B \in J^+[B \subseteq A \and |B| < \sigma]$
\ermn
and $\sigma = \text{ cf}(\sigma) > \aleph_0$.

Then 

$$
T^3_J(\mu) = T^2_J(\mu) \leqq T^1_J(\mu) \leqq T^4_J(\mu)
$$
\mn
where 

$$
\align
T^1_J(\mu) = T_J(\mu) = \sup\{|{\Cal F}|:&{\Cal F} 
\text{ is a family of functions from } \kappa  \\
  &\text{ to } \mu \text{ such that for } f \ne g \text{ from } {\Cal F}
\text{ we have } f \ne_J g\}.
\endalign
$$

$$
\align
T^2_J(\mu) = \sup \bigl\{\lambda:&\text{ there are } n_i < \omega
\text{ for } i < \kappa \text{ and regular} \\
  &\lambda_{i,\ell} > \kappa^{<\sigma} \text{ for } i < \kappa,\ell < n_i 
\text{ such that}: \lambda_{i,\ell} \leqq \mu \\
  &\text{ and } \lambda \leqq \text{ max pcf}\{\lambda_{i,\ell}:i < \kappa,
\ell < n_i\}; \text{ moreover if } A \in J^+ \\
  &(= {\Cal P}(\kappa) \backslash J) \text{ \ub{then} } \lambda \leqq
\text{ max pcf}\{\lambda_{i,\ell}:i \in A,\ell < n_i\}\}.
\endalign
$$ 
\mn
$T^3_J(\mu)$ is defined similarly but for $A \in J^+$ we demand:      

$$
\lambda = {\text{\rm max pcf\/}}\{\lambda_{i,\ell}:i \in A,\ell < n_i\}
$$ 

$$
T^4_J(\mu) = {\text{\rm Min\/}}
\{\underset s {}\to \sup T^2_{J+(\kappa \backslash 
A_n)}(\lambda):A_n \subseteq A_{n+1} \subseteq \kappa = 
\dbcu_{n < \omega} A_n,A_n \notin J\}.
$$
\endproclaim
\bigskip

\remark{\stag{6.1A} Remark}   1) Note that usually the four terms in 
the conclusion of the theorem are equal. 
\mr
\item "{$(\alpha)$}"  If $J$ is $\aleph_1$-complete then $T^2_J(\mu) = 
T^4_J(\mu)$ hence all are equal
\sn
\item "{$(\beta)$}"  all terms are equal if for $\langle A_n:n < \omega 
\rangle$ such that $A_n \subseteq \kappa,A_n \notin J,A_n \subseteq 
A_{n+1},\kappa = \dbcu_n A_n$ we have: for some $n$ and $B \subseteq A_n$ 
we have $(\kappa,J),(B,J \cap {\Cal P}(B))$ are isomorphic. 
\ermn
2)  The supremum in the definition of $T^1_J(\mu)$ is always obtained.
\nl
[Why?  If ${\Cal F}_1,{\Cal F}_2$ are as there, ${\Cal F}_1$ 
maximal $|{\Cal F}_1| < |{\Cal F}_2|$ then for every 
$f \in F_2$ there is $g_f \in {\Cal F}$, etc.].
\endremark
\bigskip

\demo{Proof}  \ub{$T^2_J(\mu) = T^3_J(\mu)$}.

Trivially $T^3_J(\mu) \leqq T^2_J(\mu)$; for the other direction let  
$\lambda$ appear in the  sup  defining $T^2_J(\mu)$, as exemplified by  
$\langle < \lambda_{i,\ell}:\ell < n_i>:i < \kappa \rangle$; as 
max pcf$\{\lambda_{i,\ell}:i \in A,\ell < n_i\}$ is always regular, without 
loss of generality $\lambda$ is regular. 

By \cite[\S1]{Sh:355}; more elaborately, for some ${\frak a} \subseteq
{\frak a}^* =: \{\lambda_{i,\ell}:\ell < n_i,i < \kappa\}$ we have 
$[{\frak a} \ne \emptyset \Rightarrow \lambda = 
\text{ max pcf}({\frak a})]$ and $\lambda \notin  
\text{ pcf}[{\frak a}^* \backslash {\frak a}]$ and ${\frak a} 
\ne {\frak a}^* \Rightarrow \lambda <  \text{ max pcf}({\frak a}^* 
\backslash {\frak a})$.  Define $\lambda'_{i,\ell}$ to be  
$\lambda_{i,\ell}$ if $\lambda_{i,\ell} \in a^*$.  If 
${\frak a} \ne {\frak a}^*$ let $u^* = \{(i,\ell):
\lambda_{i,\ell} \in {\frak a}^* \backslash {\frak a})\},J = \{u \subseteq  
u^*:\text{max pcf}\{\lambda_{(i,\ell)}:(i,\ell) \in u\} < \lambda\}$.  By 
\cite[1.5]{Sh:355} we can find regular $\lambda'_{i,\ell} < \lambda_{i,\ell}$ 
for $(i,\ell) \in u^*$ such that $\lambda = \text{ tcf }
\dsize \prod_{(i,\ell) \in u^*} \lambda'_{i,\ell}/J$.  Now 
$\langle \lambda'_{i,\ell}:\ell < n_i,i < \kappa \rangle$ exemplifies  
$\lambda \leqq T^3_J(\mu)$.
\bn
\ub{$T^3_J(\mu) \leqq T^1_J(\mu)$}.  Very easy; of course, instead 
$F \subseteq {}^\kappa \mu$ we can 
have $F \subseteq {}^\kappa Y$ as long as  
$|Y| \leqq \mu$.  For $\lambda,\langle < \lambda_{i,\ell}:\ell < n_i>:
i < \kappa \rangle$ as in the definition of $T^3_J(\mu)$, let 
$a = \{\lambda_{i,\ell}:i < \kappa,\ell < n_i\}$ and $\langle f_\alpha:
\alpha < \lambda \rangle$ be a sequence of members of $\Pi a$ which is 
$<_{J_{<\lambda}[a]}$-increasing and cofinal.  
Now let $Y = {}^{\omega >}(\mu + 1)$ and for each $\alpha < \lambda$ we 
define $g_\alpha \in {}^\kappa Y$ as follows:

$$
g_\alpha(i) =: \langle \lambda_{i,\ell}:\ell < n_i \rangle \char94  
\langle f_\alpha(\lambda_{i,\ell}):\ell < n_i \rangle .
$$
\mn
We leave the checking to the reader.  We now turn to the main case.
\bn
\ub{$T^1_J(\mu) \leqq T^4_J(\mu)$}.  
Let $\lambda$ be the right side expression  
$-T^4_J(\mu)$ (so clearly $\lambda \geqq \mu),\chi =: \beth_3(\lambda)^+$ 
and for $\zeta \leqq \omega +1$ let: $M^*_\zeta \prec (H(\chi),\in,<^*_\chi),
\|M^*_\zeta\| = \lambda,\lambda + 1 \subseteq M^*_\zeta,\zeta < \xi \leqq
\omega +1 
\Rightarrow M^*_\zeta \in M^*_\xi$.  Suppose $F \subseteq {}^\kappa 
\mu$ exemplify $T^1_J(\mu) > \lambda$ and we shall get a contradiction, 
without loss of generality $F \in M^*_0$.  Clearly for every $f \in {}^\kappa 
\mu$ we have: $\{g \in F:\neg g \ne_J f\}$ has cardinality $\leqq 
\kappa^{<\sigma}$ 
(remember $\otimes$), hence necessarily there is $f^* \in F$
such that for every $g \in {}^\kappa \mu \cap M^*_{\omega +1}$ (e.g. 
$g \in F \cap M^*_{\omega +1})$ we have $f^* \ne_J g$.  Moreover, if  
$A \subseteq \kappa,A \notin J,B \subseteq \mu,|B| \leqq \kappa^{<\sigma}$ 
then $\{f \in F:\{\alpha \in A:f(\alpha) \in B\} \notin J\}$ has cardinality
$\leqq \kappa^{<\sigma}$ (again, remember $\otimes$), so if in addition
$A,B \in M^*_\zeta$ then $[f \in F \and \{\alpha \in A:f(\alpha) \in B\} 
\notin J \Rightarrow f \in M^*_\zeta]$. 

We define by induction on $k < \omega,N^a_k,N^b_k,a_k,f^k$ such that: 
\mr
\item "{$(a)$}"  $N^a_k \prec M^*_0,N^b_k \prec (H(\chi),\in,<^*_\chi)$ 
\sn
\item "{$(b)$}"  $N^b_0$ is the Skolem Hull of $\{f^*\} \cup \{i:i \leqq  
\kappa^{<\sigma}\}$ in $({\Cal H}(\chi),\in,<^*_\chi)$ 
\sn
\item "{$(c)$}"  $N^a_0$ is the Skolem Hull of $\{i:i \leqq 
\kappa^{<\sigma}\}$ (in $({\Cal H}(\chi),
\in,<^*_\chi)$, equivalently in $M^*_0)$ 
\sn
\item "{$(d)$}"   $N^b_{n+1}$ is the Skolem Hull of $N^b_n \cup \{f^n(\theta):
\theta \in a_n\}$
\sn
\item "{$(e)$}"   $N^a_{n+1}$ is the Skolem Hull of $N^a_n \cup \{f^n(\theta):
\theta \in a_n\}$ 
\sn
\item "{$(f)$}"  ${\frak a}_n = 
N^a_n \cap \lambda^+ \cap \text{ Reg } \backslash 
(\kappa^{<\sigma})^+$
\sn
\item "{$(g)$}"  $f^n \in \Pi {\frak a}_n$ 
and for each $\theta \in a_n,f^n(\theta) > 
\text{ sup}(\theta \cap N^b_n)$ 
\sn
\item "{$(h)$}"  if ${\frak b} \subseteq {\frak a}_n,
\text{ max pcf}({\frak b}) \leqq \lambda$ and  
$|{\frak b}| < \sigma$ \ub{then} \nl
$f^n \restriction {\frak b} \in \{\text{Max}\{f^{c_\ell}
_{\lambda_\ell,\alpha_\ell} \restriction {\frak b}:
\ell < n\}:n < \omega,\alpha_\ell < \text{ max pcf}({\frak b}) \text{
and}$ 
\nl

\hskip50pt $\lambda_\ell \in \text{ pcf}({\frak a}_n) \text{ and } 
{\frak c}_\ell \in \{{\frak b}_\theta[{\frak a}_n]:
\theta \in \text{ pcf}({\frak a}_n)\}\}$
\ermn
where ${\frak b} \mapsto \langle f^{\frak b}_\alpha:
\alpha < \text{ max pcf}({\frak b}) \rangle$ for ${\frak b} 
\subseteq {\frak a}_n$ is a definable function (in $({\Cal H}(\chi),
\in,<^*_\chi),\langle f^b_\alpha:\alpha < \text{ max pcf}({\frak b})
\rangle$ as in \cite[\S1]{Sh:371}.

By \scite{2.2A} (i.e. \cite[\S1]{Sh:371}) there is no problem to do it,  
$N^a_n \prec N^b_n,N^a_n \prec M^*_0,N^a_n \prec N^a_{n+1},N^b_n \prec  
N^b_{n+1}$ and (as in \cite[3.3A,5.1A]{Sh:400}) we get $\dbcu_n N^a_n =
\dbcu_n N^b_n$, hence Rang$(f^*) \subseteq \dbcu_{n < \omega} N^a_n$.  Now 
for each $i < \kappa$ let $m(i) = \text{ min}\{m:f^*(i) \in N^a_m\}$, and we 
can find finite ${\frak e}(i) 
\subseteq \dbcu_{\ell < m(i)} {\frak a}_\ell,y(i) \subseteq  
\kappa^{<\sigma} + 1$  such that $f^*(i) \in M^{e(i),y(i)}_{m(i)}$ 
where for any ${\frak e} \subseteq \dbcu_n {\frak a}_k$ 
and $y \subseteq \kappa^{<\sigma} + 1$ (we 
define by induction of $\ell$):
\mr
\item "{${{}}$}"  $M^{e,y}_0$ is the Skolem Hull of $y$ in $({\Cal H}
(\chi),\in,<^*_\chi)$,
\sn
\item "{${{}}$}"  $M^{e,y}_{\ell +1} =$ Skolem Hull of $M^e_\ell \cup 
\{f^m(\theta):m \leqq \ell$ and $\theta \in e \cap M^e_\ell\}$.
\ermn
Clearly:  $[{\frak e} \subseteq \dbcu_m {\frak a}_m 
\and {\frak e} \in M^*_\ell \and \text{ max pcf}({\frak e}) 
\leqq \lambda \and \ell < \omega \Rightarrow M^{e,y}_\ell \in 
M^*_{\ell +1}]$, and $[{\frak e} \subseteq {\frak d} \and y 
\subseteq z \Rightarrow M^{{\frak e},y}_\ell \subseteq M^{d,z}_\ell]$ and 
$M^{e,y}_\ell \subseteq N^a_\ell$.

Let $A_n = \{i < \kappa:m(i) \leqq n\}$.  Clearly $A_n \subseteq A_{n+1}, 
\kappa = 
\dbcu_{n < \omega} A_n$, but $\kappa$ is not the union of $\aleph_0$ 
members of $J$, so for some $n(*) < \omega,[n \geqq n(*) \Rightarrow  A_n 
\notin J]$.  It suffices to prove: 
\mr
\item "{$(*)$}"  if $m(*) < \omega,A \subseteq \kappa,\dsize
\bigwedge_{i \in A} m(i) \leqq m(*),A \notin J$ then max pcf 
$(\dbcu_{i \in A} {\frak e}(i) > \lambda$ 
\ermn
[as this means $n(*) \leqq n < \omega \Rightarrow T^2_{J+(\kappa \backslash 
A_n)}(\mu) > \lambda$,  hence $\langle A_{n(*)+\ell}:\ell < \omega \rangle$  
and $\langle e(i):i < \kappa \rangle$ exemplified $T^4_J(\mu) > \lambda$  
contradiction].

We can replace $A$ by any subset which is not in $J$.

By the assumption $\otimes$ without loss of generality $|A| < \sigma$, and 
suppose $A$ contradicts $(*)$.  Let ${\frak e}^* = 
\dbcu_{i \in A} {\frak e}(i),y^* = 
\dbcu_{i \in A} y(i)$.  As $\sigma$ is infinite, clearly $|{\frak
e}^*| < \sigma,  
|y^*| < \sigma,\dsize \bigwedge_n \|M^{{\frak e}^*,y^*}_n\| 
< \sigma$ (remember $\sigma > \aleph_0)$.

Prove by induction on $\ell$ (suffice for $\ell \leqq m(*))$ that  
${\frak e}^* \cap M^{{\frak e}^*,y^*}_n \in 
M^*_{\ell +1},y^* \cap M^{{\frak e}^*,y^*} \in 
M^*_{\ell +1}$ and $M^{{\frak e}^*,y^*}_\ell 
\subseteq M^{{\frak e}^*,y^*}_{\ell +1}$. 
For $\ell = 0$ this holds as $M^{{\frak e}^*,y^*}_0 \prec 
N^a_0,\|M^{{\frak e}^*,y^*}_0\| 
< \sigma$ and ${\frak e}^* \cap M^{{\frak e}^*}_0 = {\frak e}^* 
\cap N^a_0$ is a subset of  
$N^a_0$ of cardinality $< \sigma,N^a_0 \in M^*_1,\|N^a_0\|^{<\sigma} 
= \|N^a_0\| \leqq \kappa^{<\sigma}$.  Similarly for $y^*$.  For $\ell + 1$,  
as we know $M^{{\frak e}^*,y^*}_\ell \in 
M^*_{\ell +1}$ and $f^\ell \restriction
({\frak e}^* \cap M^{{\frak e}^*,y^*}_\ell) \in 
M^*_{\ell +1}$ by (h) as max pcf$({\frak e}^*) \leqq  
\lambda$  by an assumption hence  
$M^{{\frak e}^*,y^*}_{\ell +1} \in M^* _{\ell +2}$.  As  
$\|M^{{\frak e}^*,y^*}_{\ell +1}\| = \kappa^{<\sigma},|{\frak e}^*| 
< \sigma$ and $\kappa^{<\sigma} + 1 
\subseteq M^*_0$ necessarily $e^* \cap M^{e^*,y^*}
_{\ell +1} \in M^*_{\ell +2}$.  So  
$M^{{\frak e}^*,y^*}_{m(*)} \in M^*_{m(*)+1}$ so Rang$(f^* \restriction A) 
\subseteq M^{{\frak e}^*,y^*}_{m(*)} \in 
M^*_{m(*)+1}$, so by the choice of $f^*,
A \in I$.
\enddemo
\bigskip

\demo{Proof of 6.1A(1)}  $T^4_J(\mu) \leqq T^1_J(\mu)$.  Let $A_n \subseteq  
\kappa,A_n \subseteq A_{n+1},\kappa = \dbcu_{n < \omega} A_n,A_n \notin J$  
and $T^2_{J+(\kappa \backslash A_n)}(\mu) \geqq \lambda$.  For each $n$, by 
earlier parts of the proof, there is ${\Cal F}_n 
\subseteq {}^\kappa \mu$ such that  
$|{\Cal F}_n| \geqq \lambda$ and $[f \ne g \in F_n \Rightarrow f 
\ne_{J+(\kappa \backslash A_n)},g]$.

Let $F_n = \{f^n_\alpha:\alpha < \alpha_n\},\alpha_n \geqq \lambda$  
exemplify this.  Now define $f_\alpha \in {}^\kappa \mu$ for 
$\alpha < \lambda$ as follows: for $\zeta < \kappa$ let $n(\zeta) = 
\text{ Min}\{n:\zeta \in A_n\}$ and $f_\alpha(\zeta) = \omega 
f^{n(\zeta)}_\alpha(\zeta) + n(\zeta)$.
\hfill$\square_{\scite{6.1}}$\margincite{6.1}
\enddemo
\bigskip

\demo{\stag{6.2} Conclusion}  Suppose cf$(\kappa) > \aleph_0,\kappa > \sigma  
\geqq \aleph_0$ and $I = [\kappa]^{< \sigma},\mu > \kappa^\sigma$. 
\ub{Then} $T^4_I(\mu)$ is $T^2_I(\mu)$ hence is $T^1_J(\mu)$.
\enddemo
\bigskip

\demo{Proof}  Apply \scite{6.1} ($\sigma^+$ here corresponds to $\sigma$  
there), more exactly by \scite{6.1}(A)(2).
\enddemo
\bigskip

\remark{\stag{6.3} Remark}  Asking on almost disjoint 
sets is an inessential change.
\endremark

\newpage
    
REFERENCES.  
\bibliographystyle{lit-plain}
\bibliography{lista,listb,listx,listf,liste}

\shlhetal

\enddocument

\enddocument